\documentclass{smfart}

\usepackage[english, francais]{babel}
\usepackage[latin1]{inputenc}
\usepackage[T1]{fontenc}
\usepackage{graphicx}

\usepackage[version=0.96]{pgf}
\usepackage{tikz}
\usetikzlibrary{arrows,shapes,decorations,automata,backgrounds,petri}

\usepackage{algorithm2e}

\usepackage{amsmath,amssymb,amsfonts,bbm}

\usepackage{float}

\usepackage[justification=centering]{caption}

\newtheorem{thm}{Th\'eor\`eme}[section]
\newtheorem{prop}[thm]{Proposition}
\newtheorem{lemme}[thm]{Lemme}

\newtheorem{props}[thm]{Propri\'et\'es}
\newtheorem{define}[thm]{D\'efinition}
\newtheorem{def-prop}[thm]{D\'efinition-proposition}
\newtheorem{cor}[thm]{Corollaire}
\newtheorem{notation}{Notation}

\newtheorem{ex}[thm]{Exemple}
\newtheorem{rem}[thm]{Remarque}

\newtheorem{quest}{Question}

\newcommand\Part{\mathcal P}
\newcommand\A{\mathcal A}
\newcommand\Adele{\mathbb A}
\newcommand\Ring{\mathcal R}
\newcommand\ane{\mathcal O}

\newcommand\N{\mathbb N}
\newcommand\Z{\mathbb Z}
\newcommand\Q{\mathbb Q}
\newcommand\R{\mathbb R}
\newcommand\C{\mathbb C}
\newcommand\D{\mathbb D}

\newcommand\abs[1]{\left|#1\right|}
\newcommand\norm[1]{\|#1\|}
\newcommand\eq{\Leftrightarrow}

\newcommand\F{\mathbb{F}}

\newcommand\defi[1]{\textbf{#1}}

\newcommand\Frac{\operatorname{Frac}}
\newcommand\pgcd{\operatorname{pgcd}}
\newcommand\Rat{\operatorname{Rat}}
\newcommand\Arel{\A^{\operatorname{rel}}}
\newcommand\Ared{\A^{\operatorname{red}}}
\newcommand\Lrel{L^{\operatorname{rel}}}
\newcommand\Lred{L^{\operatorname{red}}}
\newcommand\Lnonred{L^{\operatorname{non red}}}

\usepackage{hyperref}

\author{Paul MERCAT}
\address{Université Paris-Sud 11 \\ 91405 Orsay}
\email{paul.mercat@math.u-psud.fr}
\urladdr{http://www.math.u-psud.fr/~mercat}
\title{Semi-groupes fortement automatiques}
\alttitle{Strongly automatic semigroups}

\begin{document}
\frontmatter

\begin{abstract}
	Dans cet article, nous introduisons la notion de semi-groupe fortement automatique, qui entraîne la notion d'automaticité des semi-groupes usuelle.
	On s'intéresse particulièrement aux semi-groupes de développements en base $\beta$, pour lesquels on obtient un critère de forte automaticité.
\end{abstract}

\begin{altabstract}
	In this paper, we introduce the notion of strongly automatic semigroup, which implies the usual notion of automaticity.
	We focus on semigroups of $\beta$-adics developpements, for which we obtain a criterion of strong automaticity.
\end{altabstract}

\altkeywords{}

\date{\today}
\maketitle



\tableofcontents

\mainmatter

\section{Introduction} 
\subsection{Organisation de l'article}
Dans cet article, nous introduisons la notion de semi-groupe fortement automatique, qui consiste à avoir un ensemble de relations qui soit un langage rationnel, c'est-à-dire reconnaissable par un automate fini.
On démontre que la forte automaticité entraîne l'automaticité au sens usuel.

Pour les semi-groupes correspondant aux développements en base $\beta$, on démontre le résultat suivant :

\begin{thm} 
Définissons un semi-groupe $\Gamma$ engendré par les transformations :
$$ x \mapsto \beta x + t $$
pout $t \in A \subset \C$, où $A$ est une partie finie de $\C$, et $\beta$ est un nombre complexe.

Si le nombre complexe $\beta$ est transcendant, ou bien algébrique mais sans conjugué de module 1, alors pour toute partie $A \subset \C$ finie, le semi-groupe $\Gamma$ est fortement automatique.

Réciproquement, si le nombre complexe $\beta$ est algébrique et a au moins un conjugué de module 1, alors il existe une partie $A \subset \C$ finie telle que le semi-groupe $\Gamma$ n'est pas fortement automatique.
\end{thm}

On commencera par faire des rappels sur les automates (voir partie \ref{rappels_automates}). Puis dans la partie \ref{smgafa}, on définira ce qu'est un semi-groupe fortement automatique, et l'on donnera quelques propriétés.
On rappellera ensuite ce qu'est un semi-groupe automatique, et l'on montrera que les semi-groupes fortement automatiques sont automatiques.
On s'intéressera ensuite aux semi-groupes correspondants aux développements en base $\beta$ dans la partie \ref{smgdb} où l'on démontrera le théorème annoncé (voir théorème \ref{cm1} et proposition \ref{rcm1}).
Enfin, la partie \ref{Exemples} est consacrée à des exemples, et rappelle des travaux en lien avec cet article.

\subsection{Motivation}
Ces notions d'automaticité et de forte automaticité fournissent d'une part une façon de représenter un semi-groupe infini avec une quantité finie de données (et donc cela permet de manipuler facilement ce semi-groupe sur ordinateur), et d'autre part donnent des informations combinatoires sur le semi-groupe (on obtient par exemple une asymptotique très précise du nombre d'éléments du semi-groupe).
La notion de forte automaticité que l'on introduit, bien que plus simple que la notion d'automaticité, ne semble pas avoir été étudiée parce-qu'elle n'est pas intéressante pour les groupes. Mais il existe des exemples de semi-groupes fortement automatiques intéressants : \\

Voici un exemple simple de semi-groupe fortement automatique.
Considérons le semi-groupe engendré par les trois transformations affines :
$$ \left\{ \begin{matrix} 0 :& x & \mapsto &  x/3, & \\ 1 :& x & \mapsto & x/3 & + & 1, \\ 3 :& x & \mapsto & x/3 & + & 3. \end{matrix} \right. $$
Par définition, c'est l'ensemble des composées de ces trois applications. \\
Voici quelques questions que l'on peut se poser :
\begin{itemize}
	\item Quel est l'asymptotique du nombre d'éléments pour la longueur des mots ?
	\item	Comment peut-on déterminer si deux mots en les générateurs représentent le même élément du semi-groupe ?
	\item	Y a-t'il une façon de représenter les éléments du semi-groupes par des mots uniques particuliers (que l'on appelera mots réduits) ?
	\item Le semi-groupe a-t'il une présentation finie ?
\end{itemize}
La réponse à ces questions est donnée par la structure automatique du semi-groupe. Celle-ci est donnée par des automates tels que l'on peut en voir sur les figures suivantes (voir la partie \ref{rappels_automates} pour des rappels sur les automates) :

\begin{figure}[H]
\centering
\caption{Automate reconnaissant un ensemble de mots réduits du semi-groupe. Les mots réduits sont ici les mots minimaux pour l'ordre lexicographique inverse, avec $0 < 1 < 3$.} \label{fig_intro_red}
\includegraphics[scale=.8]{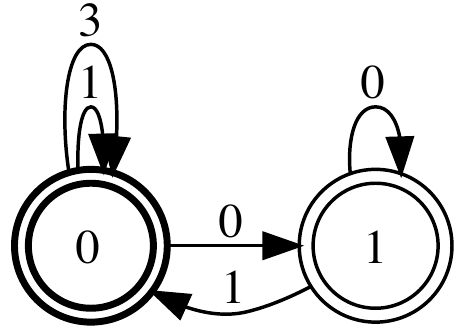}
\end{figure}

On appelle \defi{mots réduits} un choix de représentants uniques pour les éléments du semi-groupe par des mots en les générateurs.
On voit sur l'automate de la figure \ref{fig_intro_red} que les mots réduits sont ici exactement les mots ne contenant pas le mot 03.

\begin{figure}[H]
\centering
\caption{Automate reconnaissant les relations du semi-groupe.} \label{fig_intro_rel}
\includegraphics[scale=.8]{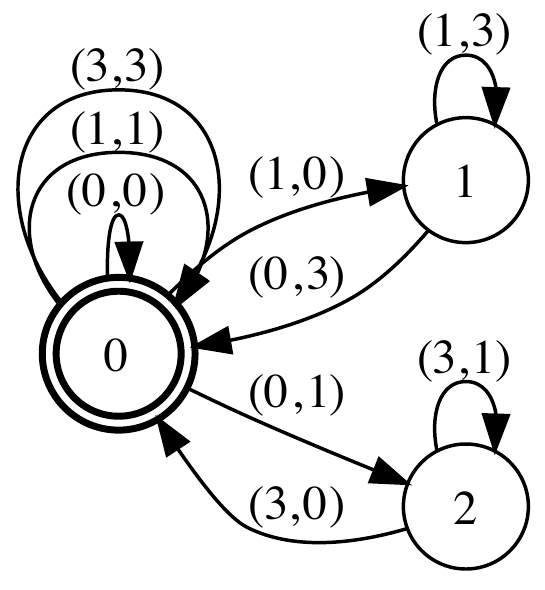}
\end{figure}

On voit sur l'automate de la figure \ref{fig_intro_rel} que les relations du semi-groupe $\Gamma$ s'obtiennent toutes à partir des relations $11^n0 = 03^n3$, par concaténation.


\begin{ex}
Le mot $(1,0)(1,3)(0,3)$ est reconnu par l'automate de la figure \ref{fig_intro_rel}, et on a en effet la relation $1\circ 1 \circ 0 = 0 \circ 3 \circ 3$, puisque l'on a l'égalité
$$ \frac{\frac{\frac{x}{3}}{3}+1}{3} + 1 = \frac{\frac{\frac{x}{3} + 3}{3} + 3}{3}. $$
\end{ex}

Deux mots $u_1...u_n$ et $v_1...v_n$ en les générateurs $\{0, 1, 3\}$ représentent le même élément du semi-groupe si et seulement si le mot $(u_1, v_1) ... (u_n, v_n)$ est reconnu par l'automate de la figure \ref{fig_intro_rel}.

L'automate de la figure \ref{fig_intro_red} fournit un moyen de connaître le nombre d'éléments du semi-groupe de longueur $n$ donnée : celui-ci est en effet égal au nombre de chemins de longueur $n$ de l'état initial $0$ vers les états finaux $0$ et $1$. Ceci est donné par la somme des deux premiers coefficients des puissances de la matrice d'adjacence du graphe :
$$
	\begin{pmatrix}
		2 & 1 \\
		1 & 1
	\end{pmatrix},
$$
dont les valeurs propres sont $\frac{\sqrt{5}+3}{2}$ et $\frac{3 - \sqrt{5}}{2}$.
Ainsi, on voit que le nombre d'éléments du semi-groupe de longueur $n$ est exactement $f_{2n+2}$, où $(f_n)_{n \in \N}$ est la suite de Fibonnacci :
\begin{eqnarray}
	& f_0 &= 0, \nonumber \\
	& f_1 &= 1, \nonumber \\
	& f_{n+2} &= f_{n+1} + f_n. \nonumber
\end{eqnarray}
En particulier, le nombre d'éléments du semi-groupe de longueur $n$ est asymptotiquement
$$ c  \left( \frac{\sqrt{5}+3}{2} \right)^n + O \left( \left(\frac{3 - \sqrt{5}}{2}\right)^n \right) $$
pour une constante $c > 0$.

L'existence de cette structure automatique assure que le semi-groupe admet une présentation finie.
Ici, celle-ci est
$$ <0, 1, 3 | 1 \circ 0 = 0 \circ 3 >. $$
Autrement dit, toutes les relations de ce semi-groupe se déduisent de la relation $1 \circ 0 = 0 \circ 3$.
Par exemple, on en déduit la relation $1 \circ 1 \circ 0 = 1 \circ 0 \circ 3 = 0 \circ 3 \circ 3$. \\

\section{Rappels sur les automates et les langages rationnels} \label{rappels_automates}
Dans cette partie, je donne des rappels sur les automates qui seront utiles dans la suite.
Les automates sont en quelques sortes des machines qui peuvent réaliser tous les calculs en temps linéaire ne nécessitant qu'une mémoire finie.
Pour plus de détails, voir par exemple \cite{oc}, 1.5.2.

\begin{define}
	On appelle \defi{automate} un quintuplet $\A := (\Sigma, Q, T, I, F)$, où
	\begin{enumerate}
		\item $\Sigma$ est un ensemble fini appelé \defi{alphabet},
		\item $Q$ est un ensemble fini d'\defi{états},
		\item $T \subseteq Q \times \Sigma \times Q$ est l'ensemble des \defi{transitions},
		\item $I \subseteq Q$ est l'ensemble des \defi{états initiaux},
		\item $F \subseteq Q$ est l'ensemble des \defi{états finaux}.
	\end{enumerate}
	On dira que l'automate est \defi{déterministe} si l'on a $\#I = 1$ et 
	$$ \left[ (p, a, q) \in T \text{ et } (p, a, r) \in T \right] \text{ implique } q = r.$$
	Autrement dit, quand l'automate $\A$ est déterministe, $T$ est le graphe d'une fonction partielle de transition $Q \times \Sigma \rightarrow Q$, et il n'y a qu'un seul état initial.
\end{define}

On considèrera parfois des automates infinis, c'est-à-dire des automates pour lesquels l'ensemble d'états $Q$ est infini.

\begin{notation}
On notera $p \xrightarrow{a} q$ si $(p, a, q) \in T$.
\end{notation}

\begin{notation}
\'Etant donné un alphabet $\Sigma$, on notera $\Sigma^* := \Sigma^{(\N)}$ l'ensemble des mots finis. Pour $u \in \Sigma^*$, on notera également $u^* := \cup_{n \in \N} \{u^n\}=\{u\}^*$.
\end{notation}

\paragraph{Représentation graphique}
On représente les automates comme des graphes dont les arêtes sont étiquetées par des lettres de l'alphabet. Sur les dessins de cet article, l'état initial est en gras, et les états finaux sont les ronds dessinés avec un trait double.

\begin{figure}[H]
   \begin{minipage}[b]{0.45\linewidth}
	\caption{Automate ayant pour états \{0, 1, 2, 3, 4\}, pour alphabet \{(0,0), (0,1), (1,0), (1,1)\}, pour ensemble d'états initiaux \{0\} et pour ensemble d'états finaux \{0\}.}
      \centering \includegraphics[scale=0.6]{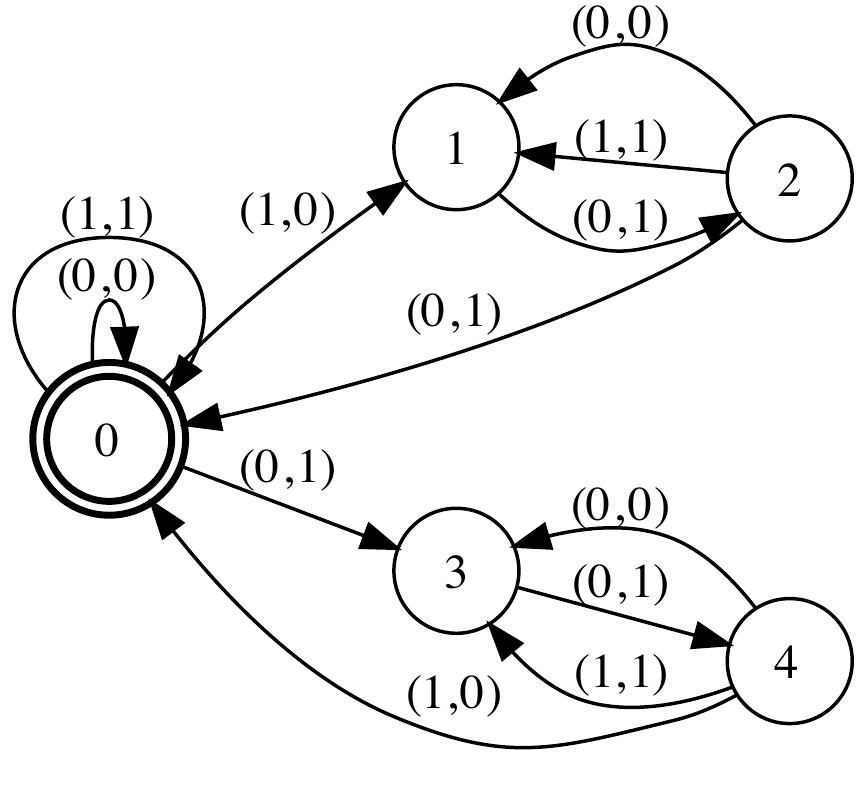}
   \end{minipage}\hfill
   \begin{minipage}[b]{0.45\linewidth}   
 	\caption{Automate ayant pour états \{0, 1, 2, 3, 4\}, pour alphabet \{0, 1, *\}, pour ensemble d'états initiaux \{0\} et pour ensemble d'états finaux \{0\}.}
      \centering \includegraphics[scale=0.6]{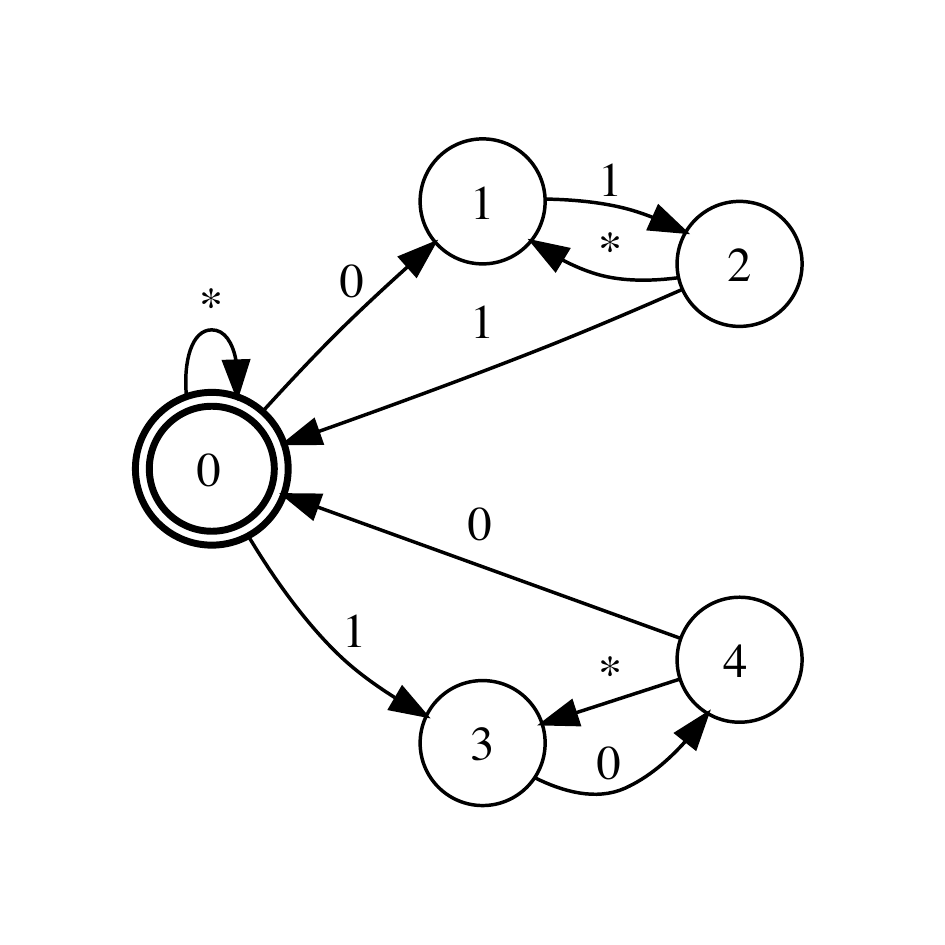}
   \end{minipage}
\end{figure}

\begin{figure}[H]
\centering
	\caption{Automate ayant pour états \{0, 1, 2, 3, 4, 5, 6\}, pour alphabet \{(0,0), (0,1), (1,0), (1,1)\}, pour ensemble d'états initiaux \{0\} et pour ensemble d'états finaux \{0, 1, 2\}.} \label{fig_ex3}
	\includegraphics[scale=.5]{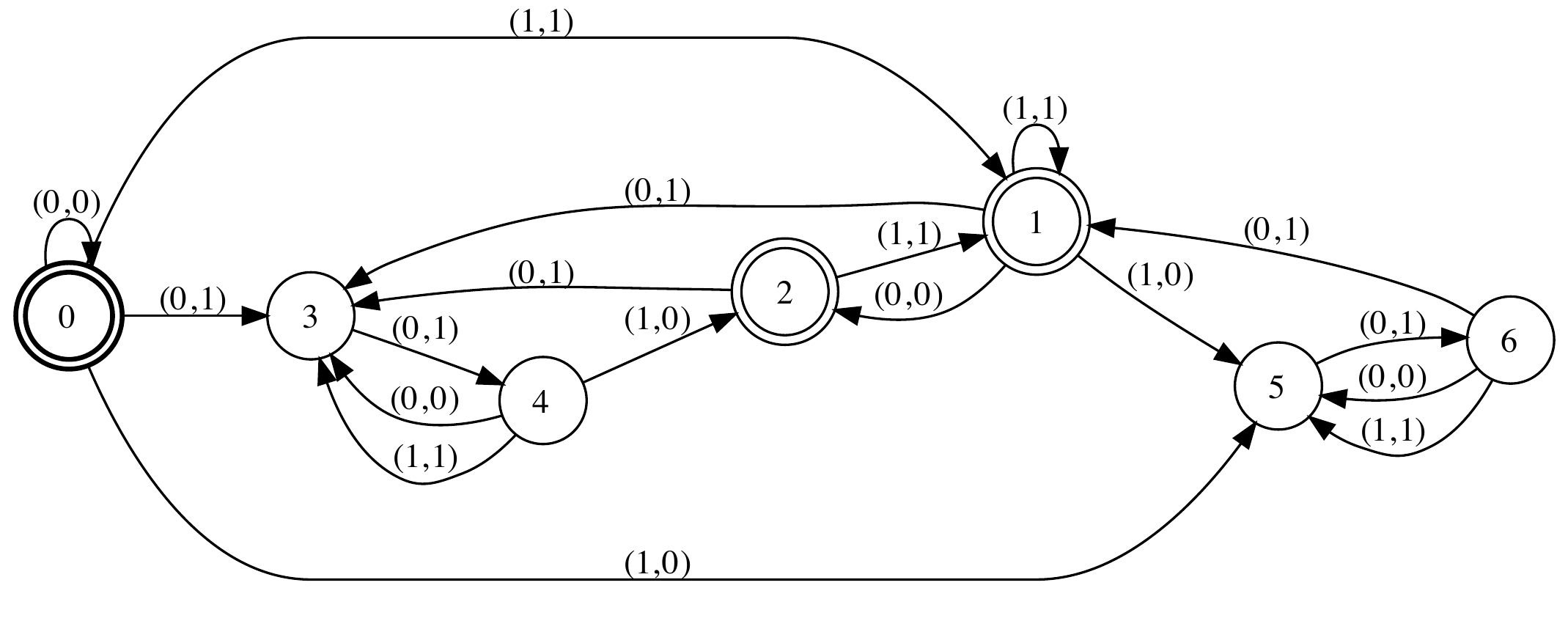}
\end{figure}

\begin{define}

On appelle \defi{langage reconnu par un automate} $\A = (\Sigma, Q, T, I, F)$ l'ensemble $L_\A$ des mots $a_1 ... a_n \in \Sigma^*$ tels qu'il existe un chemin
$$ I \ni q_0 \xrightarrow{a_1} q_1 \xrightarrow{a_2} ... ... \xrightarrow{a_{n-1}} q_{n-1} \xrightarrow{a_n} q_n \in F $$
dans l'automate $A$, d'un état initial vers un état final. \\
On dit qu'un mot $u \in \Sigma^*$ est \defi{reconnu} par l'automate $\A$ si l'on a $u \in L_\A$.
\end{define}
Un mot $a_1...a_n$ est donc reconnu par l'automate $\A$ s'il existe un chemin dans le graphe, étiqueté par $a_1, a_2, ..., a_n$, partant d'un état initial et aboutissant à un état final.
\begin{rem}
Si l'automate est déterministe, un tel chemin est unique.
\end{rem}

\begin{figure}[!h]
\centering
\caption{Automate reconnaissant l'ensemble des mots de la forme $a(baa)^n$.} \label{fig_0(101)*}
\includegraphics[scale=.8]{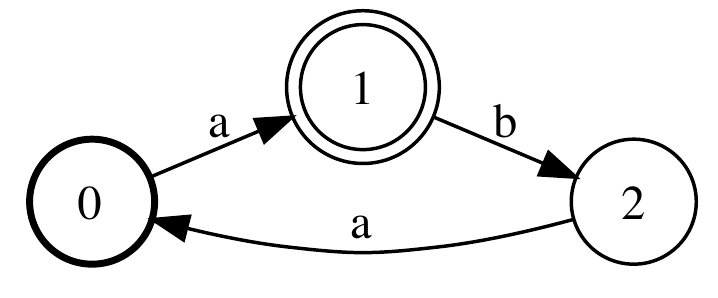}
\end{figure}

\begin{figure}[!h]
\centering
\caption{Automate non déterministe reconnaissant l'ensemble de mots \{lapin, laitue\}.}  \label{fig_lapin_laitue}
\includegraphics[scale=.6]{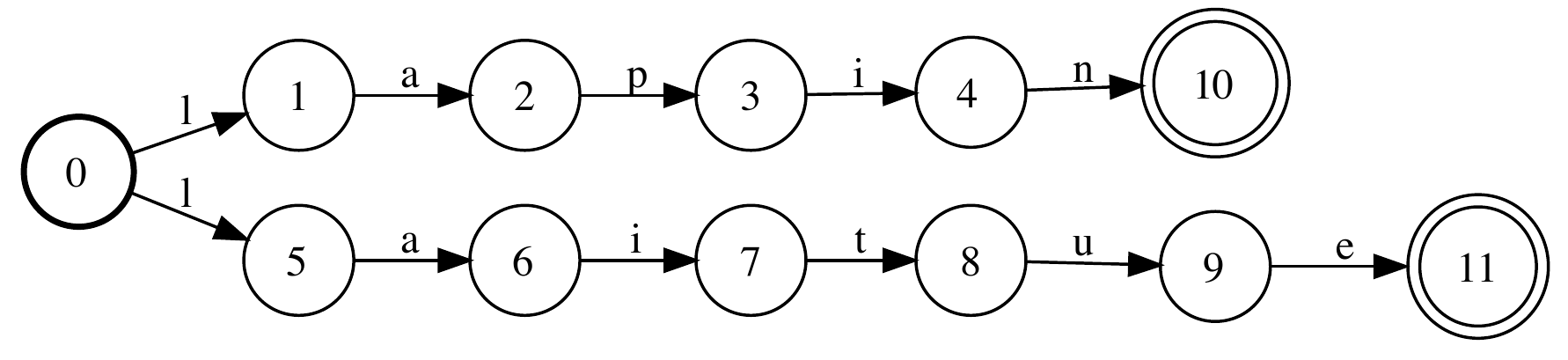}
\end{figure}

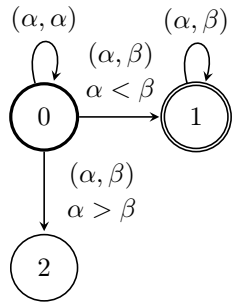
\begin{figure}[H]
	\centering
	\caption{Automate reconnaissant les couples $(u,v)$ de mots avec $u$ strictement inférieur à $v$ dans l'ordre lexicographique.} \label{fig_lex}	
	\begin{tikzpicture}[->,>=stealth,shorten >=1pt,auto,node distance=2cm,semithick]
		\tikzstyle{every state}=[fill=white,draw=black,text=black]
		\node[state][very thick]	(0)					{0};
		\node[state][double]		(1) [right of =0]	{1};
		\node[state]		(2) [below of =0]	{2};
		
		\path	(0)	[loop above]		edge	 node {$	(\alpha , \alpha)
										$} (0);
		\path(0)					edge	 node {$	\begin{array}{c}
												(\alpha, \beta) \\
												\alpha < \beta
											\end{array}
										$} (1)
								edge node {$	\begin{array}{c}
												(\alpha, \beta) \\
												\alpha > \beta
											\end{array}
										$} (2);
		\path	(1)	[loop above]		edge	 node {$	(\alpha , \beta)
										$} (1);
	\end{tikzpicture}
\end{figure}

\begin{ex}
Voir figures \ref{fig_ex3}, \ref{fig_0(101)*}, \ref{fig_lapin_laitue} et \ref{fig_lex}.
\end{ex}

\begin{figure}[!h]
\centering
\caption{Automate reconnaissant l'ensemble des nombres écrits en binaires qui sont divisibles par 3.} \label{fig_A1min} 
\includegraphics[scale=.8]{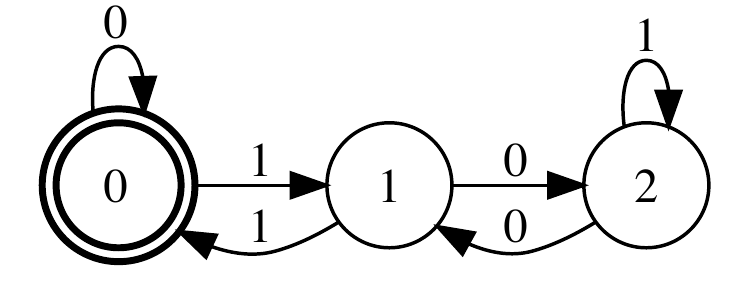}
\end{figure}

Dans la figure \ref{fig_lex}, dire qu'un couple $(u, v)$ de mots de $\Sigma^*$ est reconnu signifie qu'un mot $(u_1, v_1) ... (u_n, v_n) \in (\Sigma \times \Sigma)^*$ est reconnu, avec $u = u_1 ... u_n$ et $v = v_1 ... v_n$. On suppose que l'alphabet $\Sigma$ est muni d'une relation d'ordre totale.

\begin{define}
On dit que deux automates $\A$ et $\A'$ sont \defi{équivalents} s'ils reconnaissent le même langage : $L_\A = L_{\A'}$.
\end{define}

\begin{figure}[!h]
\centering
\caption{Automate déterministe équivalent à celui de la figure \ref{fig_lapin_laitue}.}  \label{fig_lapin_laitue_det}
\includegraphics[scale=.6]{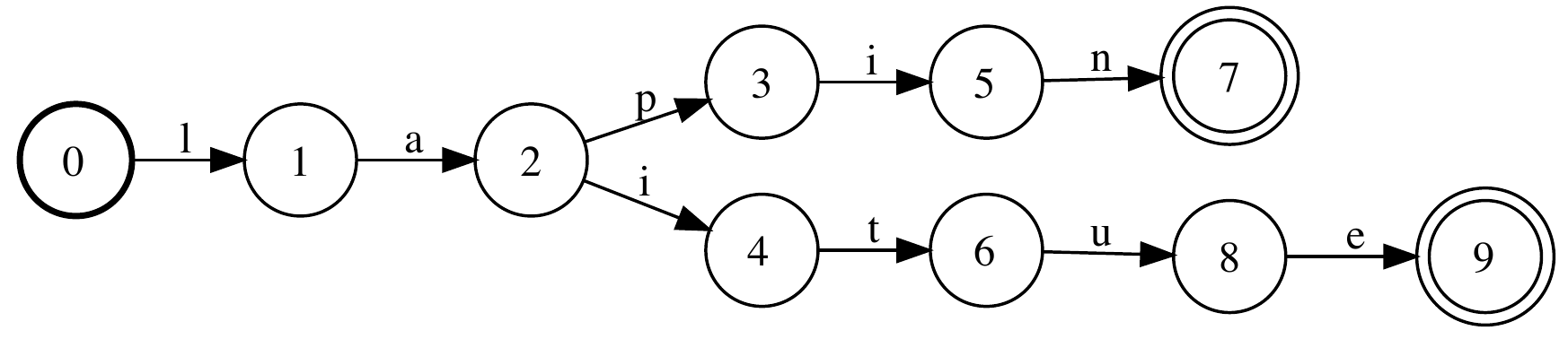}
\end{figure}

On a la proposition suivante :

\begin{prop}
Tout automate est équivalent à un automate déterministe.
\end{prop}

\begin{define}
On appelle \defi{automate minimal} d'un automate $\A$, un automate $\A'$ déterministe, équivalent à $\A$, et ayant un nombre minimal de sommets pour ces propriétés.
\end{define}

\begin{prop}
L'automate minimal d'un automate $\A$ est unique. De plus si l'automate $\A$ est déterministe, alors l'automate minimal s'obtient comme le quotient de l'automate $\A$ par une relation d'équivalence consistant à identifier des sommets entre eux.
\end{prop}

\begin{ex}
L'automate de la figure \ref{fig_lapin_laitue_det} est minimal.
\end{ex}

\begin{define}
On appelle \defi{transposée} d'un automate $\A = (\Sigma, Q, T, I, F)$ l'automate
$$ \A^t := (\Sigma, Q, T^t, F, I) $$
où $T^t := \{ (p, a, q) \in Q \times \Sigma \times Q | (q, a, p) \in T \}$.
\end{define}

\begin{rem}
Le langage reconnu par l'automate transposé $\A^t$ est la transposée du langage reconnu par l'automate initial $\A$.
\end{rem}

\begin{define}
On appelle \defi{émondé} d'un automate, l'automate restreint aux sommets par lesquels il passe un chemin d'un état initial à un état final.
On dit qu'un automate est \defi{émondé} s'il est égal à son émondé.
\end{define}
Autrement dit, un automate est émondé s'il n'existe pas de sommet qui ne sert à rien !

\begin{prop} \label{dt-min}
Un automate émondé, déterministe, et de transposée déterministe est minimal.
\end{prop}

\begin{figure}[H]
\centering
\caption{Automate minimal car vérifiant les conditions de la proposition \ref{dt-min}.}  \label{fig_dt-min}
\includegraphics[scale=.5]{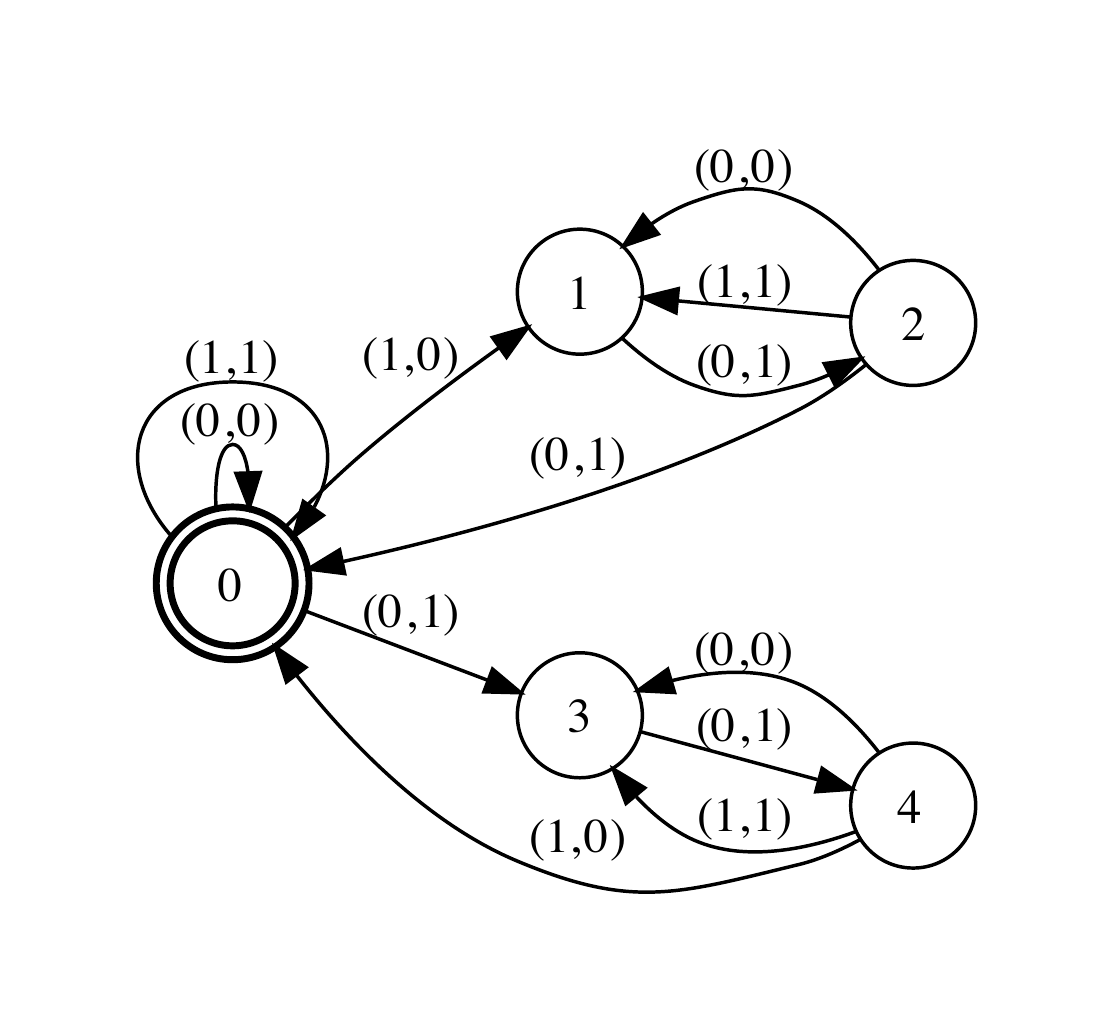}
\end{figure}

\begin{notation}
Dans tout l'article, $\epsilon$ dénote le \defi{mot vide}, c'est-à-dire le mot ayant $0$ lettres.
\end{notation}

\begin{define}
On définit l'ensemble des \defi{langages rationnels} comme étant la plus petite partie $\Rat \subset \Part(\Sigma^*)$ de l'ensemble des langages sur l'alphabet $\Sigma$ vérifiant :
\begin{enumerate}
	\item	$\emptyset \in \Rat$,
	\item $\{ \epsilon \} \in \Rat$ (où $\epsilon$ est le mot vide),
	\item	$\{ a \}  \in \Rat$ pour tout $a \in \Sigma$,
	\item $\Rat$ est stable par union : $L, L' \in \Rat \text{ implique } L \cup L'  \in \Rat$.
	\item $\Rat$ est stable par concaténation : $L, L' \in \Rat \text{ implique } L L'  \in \Rat$.
	\item $\Rat$ est stable par complémentaire : $L \in \Rat \text{ implique }  \Sigma^* \backslash L \in \Rat$.
	\item $\Rat$ est stable par étoile : $L \in \Rat \text{ implique }  L^* \in \Rat$.
\end{enumerate}
\end{define}

On a le théorème suivant :

\begin{thm}[Kleene]
Un langage $L \subseteq \Sigma^*$ est rationnel si et seulement si c'est le langage d'un automate.
\end{thm}

Voir par exemple \cite{oc}, théorème 1.59, page 36.

\begin{define}
\'Etant donnés deux langages $L \subseteq \Sigma^*$ et $K \subseteq \Lambda^*$, respectivement sur les alphabets $\Sigma$ et $\Lambda$, on appelle \defi{produit} des deux langages, le langage noté $L \times K$ sur l'alphabet $\Sigma \times \Lambda$ définit par :
$$ L \times K := \{ (a_1, b_1)...(a_n, b_n) \in (\Sigma \times \Lambda)^* | a_1 ... a_n \in L \text{ et } b_1 ... b_n \in K\}. $$
\end{define}

\begin{prop} \label{prod_rat}
Le produit de deux langages rationnels est un langage rationnel.
\end{prop}

\begin{define}
\'Etant donné un langage $L \subset (\Sigma \times \Lambda)^*$, on définit les projetés $p_1(L) \subset \Sigma^*$ et $p_2(L) \subset \Lambda^*$ du langage $L$ par
$$ p_1(L) := \{ u \in \Sigma^* | \text{Il existe } v \in \Lambda^* \text{ tel que } (u,v) \in L \}, $$
$$ p_2(L) := \{ v \in \Lambda^* | \text{Il existe } u \in \Sigma^* \text{ tel que } (u,v) \in L \}. $$
\end{define}

\begin{prop} \label{proj_rat}
Un projeté d'un langage rationnel est un langage rationnel.
\end{prop}

Voir \cite{oc}, Proposition 1.95.

\begin{lemme}[de l'étoile] \label{le}
Si un langage $L$ est rationnel, alors il existe une constante $N > 0$ telle que pour tout mot $u_1u_2u_3 \in \Sigma^*$ avec $\abs{u_2} > N$, il existe trois mots $v_1, v_2$ et $v_3$ avec $\abs{v_2} > 0$ tels que l'on ait $u_2 = v_1 v_2 v_3$ et
$$ \text{ pour tout entier } n \in \N, \quad u_1 v_1 {v_2}^n v_3 u_3 \in L. $$
\end{lemme}
(Ici, $\abs{u}$ dénote la longueur du mot $u$)

\section{Semi-groupe automatique et fortement automatique} \label{smgafa}
Dans cette partie, nous allons définir ce que sont les structures automatique et fortement automatique, et nous allons voir comment déterminer la structure fortement automatique d'un sous-semi-groupe de type fini d'un groupe.

Dans toute la suite, $\Gamma$ est un sous-semi-groupe d'un groupe $G$, et $\Sigma$ est une partie génératrice finie de $\Gamma$.
On notera $e$ l'élément neutre du groupe $G$.

\subsection{Semi-groupe fortement automatique}

\begin{define}
	On dit que le semi-groupe $\Gamma$ est \defi{fortement automatique} pour la partie génératrice $\Sigma$ si l'ensemble des relations
	$$\Lrel := \{ (u_1, v_1)...(u_n, v_n) \in (\Sigma \times (\Sigma \cup \{ e \}))^* | u_1...u_n = v_1...v_n \text{ dans } \Gamma \}$$
	est un langage rationnel. 
	On dira qu'un semi-groupe est fortement automatique s'il existe une partie génératrice pour laquelle il est fortement automatique. \\
	On appelera \defi{automate des relations} du semi-groupe $\Gamma$ l'automate minimal reconnaissant le langage $\Lrel$.
\end{define}

\begin{ex}
 Le semi-groupe engendré par les trois transformations
 $$\left\{ \begin{matrix} 0 :& x & \mapsto &  x/3, & \\ 1 :& x & \mapsto & x/3 & + & 1, \\ 3 :& x & \mapsto & x/3 & + & 3. \end{matrix} \right.$$
 est fortement automatique : voir figure \ref{fig_intro_rel}.
\end{ex}


\begin{rem}
L'automate minimal des relations ne dépend pas du plongement du semi-groupe dans un groupe, puisque le langage $\Lrel$ n'en dépend pas.
\end{rem}

\paragraph{Obtention de la structure fortement automatique d'un semi-groupe}
Dans ce paragraphe, nous allons voir comment obtenir la structure fortement automatique d'un sous-semi-groupe d'un groupe.

\begin{define}
	On définit un automate $\A = (\Sigma_\A, Q, T, I, F)$ (éventuellement infini) de la façon suivante :
	\begin{enumerate}
		\item $\Sigma_\A = \Sigma \times (\Sigma \cup \{ e\})$,
		\item $Q = \Gamma^{-1} (\Gamma \cup \{e\}) = \Gamma^{-1} \Gamma \cup \Gamma^{-1} \subseteq G$,
		\item	$I = \{ e \}$,
		\item $F = \{ e\}$,
		\item $T$ est définit par :
			$(p, (g,h), q) \in T  \text{ si et seulement si } q = g^{-1} p h$.
	\end{enumerate}
\end{define}

\begin{prop} \label{osfa}
L'émondé de l'automate $\A$ est l'automate des relations du semi-groupe $\Gamma$.
\end{prop}

\emph{Preuve :}
Montrons que l'automate $\A$ reconnait bien le langage $\Lrel$. \\
Si $(a_1, b_1)...(a_n, b_n)$ est un mot reconnu par l'automate, alors par définition on a $a_n^{-1}...a_1^{-1} e b_1 ... b_n = e$ dans $G$. Donc on a bien $a_1 ... a_n = b_1 ... b_n$ dans $\Gamma$. Réciproquement, la relation $a_1 ... a_n = b_1 ... b_n$ dans $\Gamma$ donne un chemin
$$e \xrightarrow{(a_1, b_1)} a_1^{-1} b_1 \rightarrow ... \xrightarrow{(a_n, b_n)} a_n^{-1} ... a_1^{-1} b_1 ... b_n = e$$
dans l'automate $\A$.

Montrons maintenant que l'automate émondé est minimal.
D'après la proposition \ref{dt-min}, il suffit de montrer qu'il est déterministe et de transposée déterministe.
L'automate est clairement déterministe, et la transposée s'obtient en remplaçant l'ensemble des transitions par les $p \xrightarrow{(g,h)} q$ pour $q = g p h^{-1}$, ce qui donne bien un automate déterministe.
$\Box$

\begin{props} \label{pptsfa}
	On a les propriétés :
	\begin{enumerate}
		\item	Un groupe est fortement automatique si et seulement s'il est fini.
		\item Un semi-groupe est libre (pour un système de générateurs donné) si et seulement si son automate des relations est trivial (c'est-à-dire réduit à un seul état).
		\item Si un semi-groupe fortement automatique contient une relation entre deux éléments de longueurs distinctes, alors c'est un groupe fini.
		\item Il y a unicité de la partie génératrice pour laquelle un semi-groupe ne contient pas de relation entre deux éléments de longueurs distinctes.
	\end{enumerate}
\end{props}

De ces propriétés, on déduit qu'un semi-groupe possédant une partie génératrice pour laquelle il n'y a pas de relation entre deux éléments de longueurs distinctes est automatique si et seulement si il l'est pour cette partie génératrice. \\

\emph{Preuve des propriétés :}
\begin{enumerate}
	\item	Si $\Gamma$ est un groupe, alors l'automate $\A$ est déjà émondé et a pour ensemble de sommets $\Gamma$, d'où la première propriété.

	\item Si le semi-groupe $\Gamma$ n'est pas libre pour la partie génératrice $\Sigma$, alors il existe une relation $a_1...a_n = b_1...b_n$ pour des éléments $a_i \in \Sigma$ et $b_i \in \Sigma \cup \{e\}$, avec $n \geq 1$ et $a_1 \neq b_1$. Le mot $(a_1, b_1) ... (a_n, b_n)$ est reconnu par l'automate des relations, et donc il existe une arête de $e$ à ${a_1}^{-1}b_1 \neq e$ dans l'automate des relations.
Réciproquement, si l'automate des relations n'est pas trivial, alors il existe un chemin de $e$ vers un état $g \neq e$, et de $g$ vers $e$ puisque l'automate est émondé. Le chemin de $e$ vers $e$ obtenu en concaténant ces deux chemins fournit alors une relation non triviale dans le semi-groupe $\Gamma$ (puisque les relations triviales de $\Gamma$ correspondent à des chemins qui ne passent que par l'état $e$ dans l'automate des relations).

	\item 
	Supposons qu'il existe une relation $u=v$ dans le semi-groupe fortement automatique $\Gamma$, pour deux mots $u$ et $v \in \Sigma^*$, avec $u$ de longueur strictement supérieure à $v$ : $\abs{u} > \abs{v}$.
Considérons alors le mot $w_n \in \Sigma \times (\Sigma \cup \{e\})$ correspondant au couple $(u^n,v^n)$. Celui-ci termine par au moins $n$ fois une lettre de la forme $(a, e)$ pour des lettres $a \in \Sigma$.
On a la relation $u^n = v^n$ dans le semi-groupe $\Sigma$, et donc le mot $w_n$ est reconnu par l'automate $\Arel$.
En utilisant le lemme de l'étoile (voir \ref{le}), avec $u_2$ de la forme $(a_1, e)...(a_k,e)$, on obtient, pour un entier $k$ assez grand (et donc pour un entier $n$ assez grand), une relation de la forme
$$ u_1 u_2^k u_3 =  v_1 e^{\alpha k + \beta} = v_1 \text{ dans } \Gamma, \text{ pour tout entier } k \in \N,$$ 
avec $\alpha > 0$, et $u_2$ de longueur $\alpha$.
On a donc $u_2^k = u_1^{-1} v_1 u_3^{-1} = cste$, d'où $u_2 = e$ dans le semi-groupe $\Gamma$.
De ceci, on déduit l'existence d'un générateur $a \in \Sigma$ qui est inversible dans $\Gamma$ (i.e. $a^{-1} \in \Gamma$). Posons $a' \in \Sigma^*$ tel que $a' = a^{-1}$ dans le semi-groupe $\Gamma$.
En considérant maintenant la relation $a^n(a')^{n}b^n = b^n$ dans le semi-groupe $\Gamma$, pour un générateur $b \in \Sigma$, le lemme de l'étoile nous donne, en prenant $n$ assez grand, l'égalité
$a^{n}(a')^nb^{n + k p} = b^n$ dans le semi-groupe $\Gamma$, pour un entier $p > 0$, et pour tout entier $k \in \N$. On en déduit que le générateur $b$ est inversible.
Tous les générateurs du semi-groupe étant ainsi inversibles, on en déduit que le semi-groupe $\Gamma$ est un groupe.
Par la première propriété, c'est un groupe fini.

	\item Supposons que l'on ait deux partie génératrices $\Sigma$ et $\Sigma'$ du même semi-groupe $\Gamma$. Alors un générateur $a \in \Sigma$ s'écrit comme produit d'éléments de $\Sigma'$ qui eux-même s'écrivent chacun comme produit d'éléments de $\Sigma$. Mais la longueur en la partie génératrice $\Sigma$ du dernier produit obtenu doit être 1 puisque sinon on obtiendrait une relation entre des éléments de longueurs distinctes. On en déduit que l'on a $a \in \Sigma'$. Ainsi, on a l'inclusion $\Sigma \subseteq \Sigma'$, et par symétrie $\Sigma = \Sigma'$.
	
\end{enumerate}
$\Box$

\begin{ex}
Le sous-semi-groupe $\Z_{\geq 1}$ est fortement automatique, tandis que le sous-semi-groupe $\Z_{\geq 2}$ ne l'est pas.
\end{ex}

\begin{ex}
Le sous-semi-groupe de $SL(2,\Z)$ engendré par les trois matrices $\begin{pmatrix} 1 & 2 \\ 0 & 1 \end{pmatrix}$, $\begin{pmatrix} 2 & -1 \\ 1 & 0 \end{pmatrix}$ et $\begin{pmatrix} 1 & 0 \\ -2 & 1 \end{pmatrix}$ est fortement automatique.
\end{ex}

\noindent Voir la partie \ref{Exemples} pour plus d'exemples.

\begin{rem}
	Déterminer si une partie finie d'un groupe engendre un semi-groupe libre (et donc déterminer si l'automate des relations correspondant est trivial) est décidable pour les sous-semi-groupes de type fini de $GL(2, \Z)$ mais est indécidable pour les sous-semi-groupes de $SL(3,\N)$ de type $\geq 13$. C'est une question ouverte pour les sous-semi-groupes de type fini de $SL(2, \Q)$. Voir \cite{CN} pour plus de détails.
	Ainsi, il ne peut pas exister d'algorithme pour déterminer la structure fortement automatique des sous-semi-groupes de type fini de $SL(3,\N)$.
\end{rem}

\begin{quest}
Le problème de déterminer si un sous-semi-groupe de type fini de $PSL(2,\Z) \simeq \Z/{2\Z} * \Z/{3\Z}$ est fortement automatique est-il décidable ? 
\end{quest}
(Ici, $\Z/{2\Z} * \Z/{3\Z}$ est le produit libre des deux groupes $\Z/{2\Z}$ et $\Z/{3\Z}$)

\subsection{Semi-groupe automatique}

La structure fortement automatique est utile pour déterminer facilement si deux mots représentent le même élément du semi-groupe, et nous allons voir qu'elle permet également d'obtenir d'autres informations sur le semi-groupe puisqu'elle impliquera la structure automatique usuelle :

\begin{define}
	On dit que le semi-groupe $\Gamma$ est \defi{automatique} s'il existe une partie génératrice $\Sigma$, un automate $\Ared$ appelé \defi{automate des mots réduits} et une famille d'automates $(\A^{g})_{g \in \Sigma}$ appelés \defi{automates de multiplication} vérifiants les propriétés :
	\begin{enumerate}
		\item	$\Ared$ a pour alphabet $\Sigma$,
		\item Le langage $L_{\Ared}$ est un ensemble de mots réduits. C'est-à-dire que l'on a : 
			$$\text{ pour tout } \gamma \in \Gamma, \text{ il existe un unique } u \in L_{\Ared} \text{ tel que } \gamma = u \text{ dans } \Gamma.$$
		
		Autrement dit, l'application de $\Sigma^*$ dans $\Gamma$ induit une bijection de $L_{\Ared}$ dans $\Gamma$ : on peut identifier un mot réduit à un élément de $\Gamma$.
		
		\item Pour tout $g \in \Sigma$, $\A^g$ a pour alphabet $(\Sigma \cup \{e\}) \times (\Sigma \cup \{e\})$,
		\item Pour tout $g \in \Sigma$, l'automate $\A^g$ reconnait si un mot réduit de $\Gamma$ s'obtient à partir d'un autre par multiplication à droite par $g$. Plus précisément, le langage $L_{\A^g}$ est l'ensemble des mots $(u_1, v_1) ... (u_n, v_n)$ vérifiant :
		\begin{eqnarray}
			&(u_1, v_1) ... (u_n, v_n) \in (L_{\Ared}g \times L_{\Ared}e^*) \cup (L_{\Ared}ge^* \times L_{\Ared}), \\
			&u_1...u_n = v_1...v_n \text{ dans } \Gamma.
		\end{eqnarray}
	\end{enumerate}
\end{define}

\begin{rem}
Il existe des variantes de la notion de structure automatique, qui autorise par exemple une notion de mots réduits un peu plus souple, ou encore qui reconnait si un mot réduit de $\Gamma$ s'obtient à partir d'un autre par multiplication à droite par $g$ d'une façon différente. Voir \cite{cain} pour plus de détails.
\end{rem}

\begin{figure}[H]
	\centering
	\caption{Structure automatique de $\Gamma=\Z$ avec $\Sigma=\{-1, 1\}$} \label{fig_ex_Z}
	\begin{tabular}{ccc}
	\includegraphics[scale=.4]{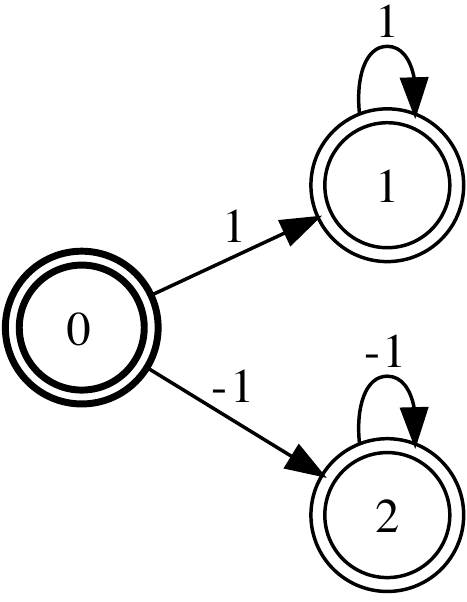} & 
	\includegraphics[scale=.4]{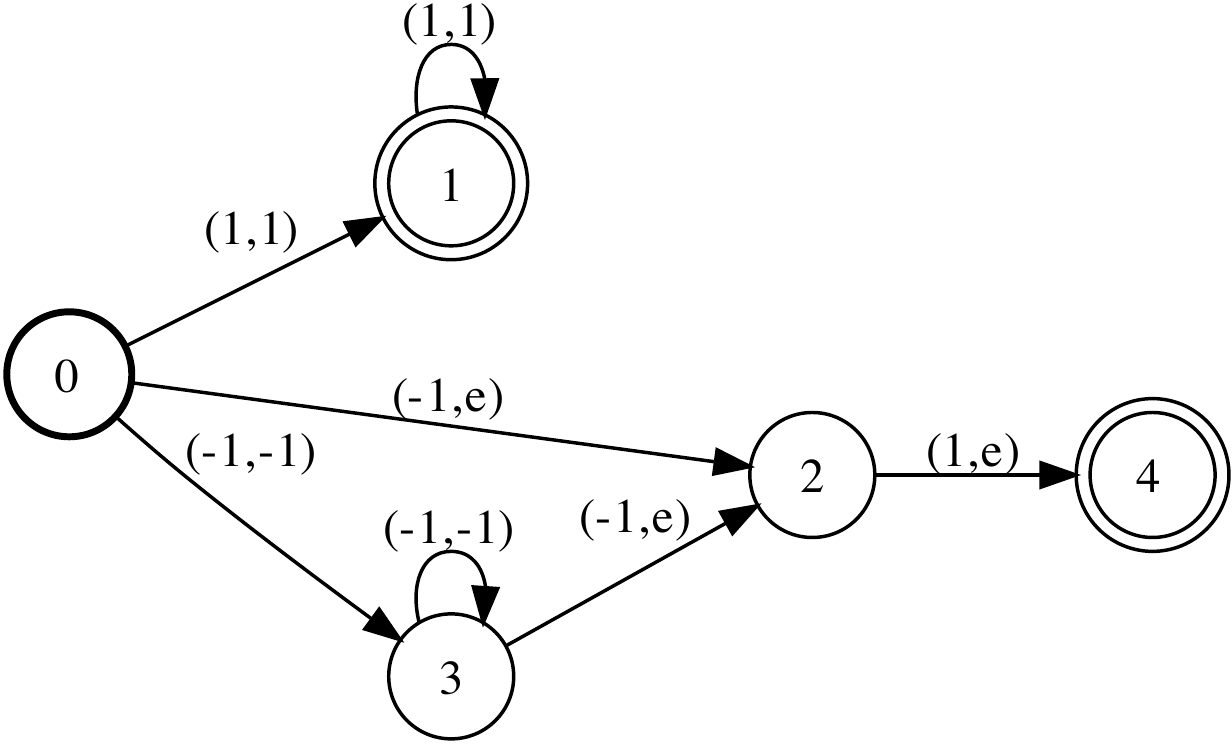} &
	\includegraphics[scale=.4]{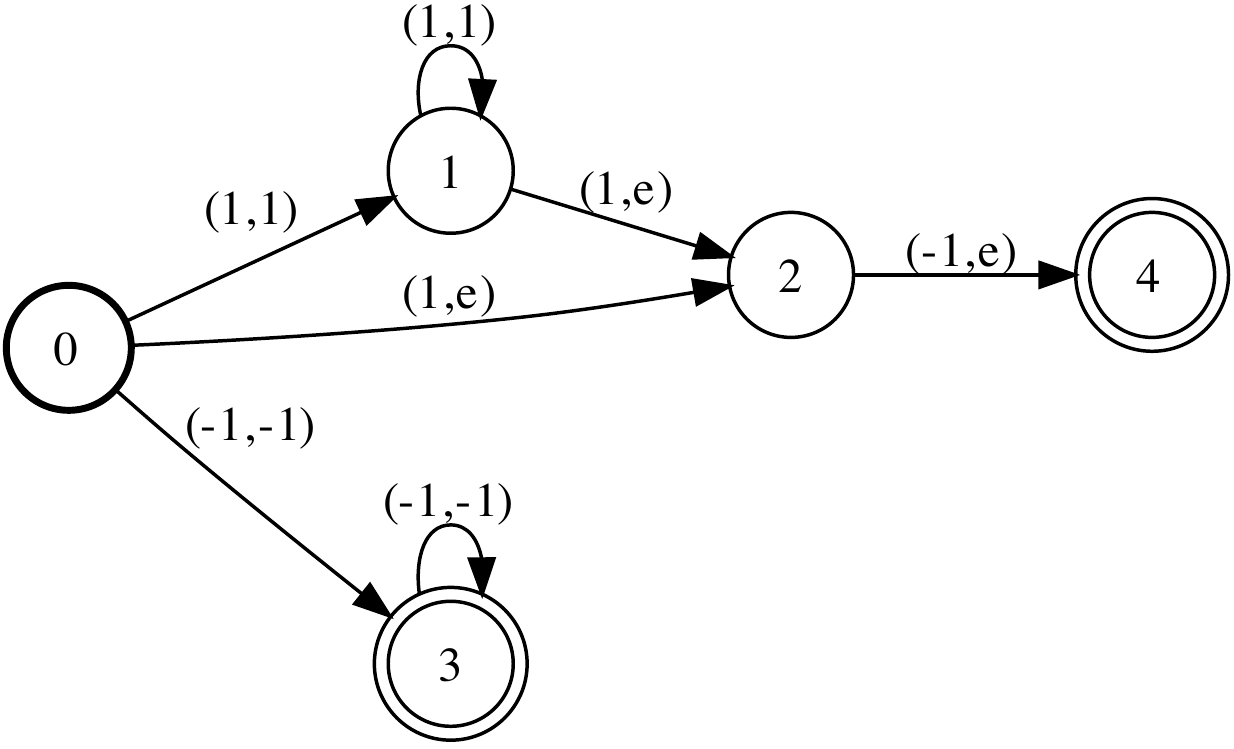} \\
	$\Ared$ & $\A^1$ & $\A^{-1}$
	\end{tabular}
	\end{figure}

\begin{figure}[H]
\centering
\caption{ }{\begin{minipage}{\linewidth} Automate $\A^{0}$ du semi-groupe engendré par les trois applications \begin{center} $\left\{ \begin{matrix} 0 :& x & \mapsto &  x/3, & \\ 1 :& x & \mapsto & x/3 & + & 1, \\ 3 :& x & \mapsto & x/3 & + & 3. \end{matrix} \right.$ \end{center} \end{minipage}} \label{fig_ex_mul0}
\includegraphics[scale=.65]{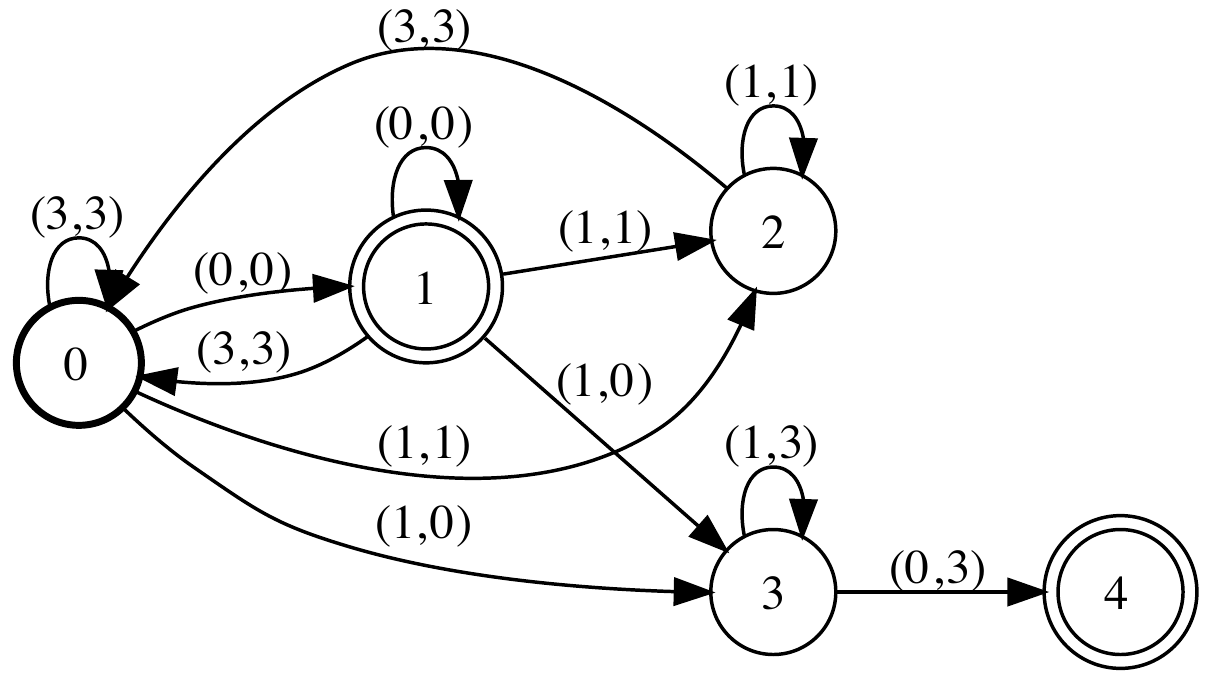}
\end{figure}

\begin{figure}[H]
\centering
\caption{ }{\begin{minipage}{\linewidth} Automate $\A^{1}$ du semi-groupe engendré par les trois applications \[ \left\{ \begin{matrix} 0 :& x & \mapsto &  x/3, & \\ 1 :& x & \mapsto & x/3 & + & 1, \\ 3 :& x & \mapsto & x/3 & + & 3. \end{matrix} \right. \] \end{minipage}} \label{fig_ex_mul1}
\includegraphics[scale=.65]{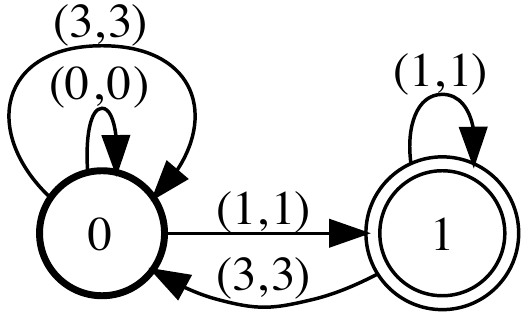}
\end{figure}

\begin{figure}[H]
\centering
\caption{ }{\begin{minipage}{\linewidth} Automate $\A^{3}$ du semi-groupe engendré par les trois applications \[ \left\{ \begin{matrix} 0 :& x & \mapsto &  x/3, & \\ 1 :& x & \mapsto & x/3 & + & 1, \\ 3 :& x & \mapsto & x/3 & + & 3. \end{matrix} \right.\] \end{minipage}} \label{fig_ex_mul3}
\includegraphics[scale=.65]{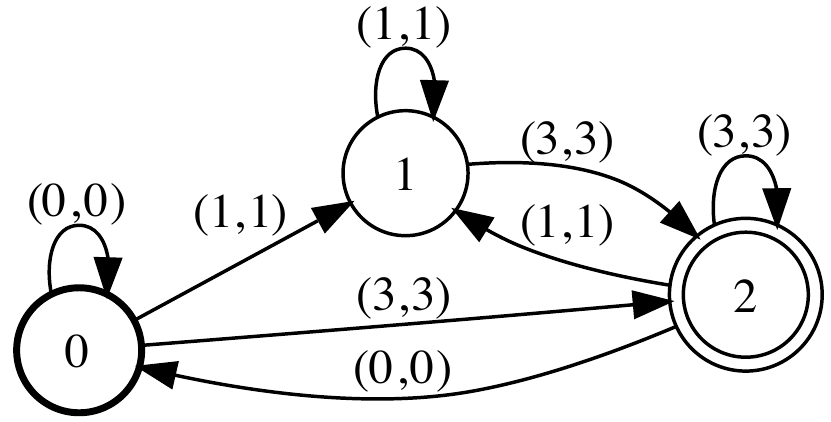}
\end{figure}

\begin{ex}
	Les semi-groupes suivants sont automatiques :
	\begin{itemize}
		\item $\Z$ : voir figure \ref{fig_ex_Z}.
		\item Le semi-groupe donné en introduction : voir figures \ref{fig_ex_red}, \ref{fig_ex_mul0}, \ref{fig_ex_mul1} et \ref{fig_ex_mul3}.
		\item Les groupes hyperboliques.
	\end{itemize}
	Il y a de nombreux autres exemples de groupes automatiques. Voir par exemple \cite{wpg}.
\end{ex}

La structure automatique d'un semi-groupe permet de manipuler celui-ci grâce à un système de mots réduits.
Cela permet par exemple d'effectuer des calculs dans le semi-groupe sur ordinateur.
Mais voici aussi plusieurs résultats donnant des informations sur le semi-groupe à partir de la structure automatique :

\begin{prop}
Si $\Gamma$ est automatique, alors il admet une présentation finie.
\end{prop}

Voir par exemple \cite{cain}.

\begin{prop}
Si $\Gamma$ est automatique, alors le nombre $c_n$ d'éléments de $\Gamma$ de longueur $n$ vérifie
$$ c_n = P(n) \alpha^n (1 + O_{n \rightarrow \infty}(e^{-\epsilon n})) $$
où $P$ est un polynôme, $\epsilon > 0$ est un réel, et $1 \leq \alpha \leq \# \Sigma$ est un réel.

Plus précisément, $c_n$ s'obtient comme la somme de coefficients de la puissance $n^{\text{ième}}$ de la matrice d'adjacence du graphe de l'automate $\Ared$.
\end{prop}

\emph{Preuve :}
Le nombre $c_n$ cherché est le nombre de chemins de longueur $n$ de l'automate $\Ared$, partant de l'état initial et aboutissant à un état final.
Ceci est donné par les puissances de la matrice d'adjacence du graphe.
Si l'automate est supposé émondé, le théorème de Perron-Frobenius nous donne l'existence d'une plus grande valeur propre $\alpha > 0$, pour laquelle on a l'asymptotique annoncée.
$\Box$

\begin{ex}
Pour $\Gamma = \Z$ et $\Sigma=\{-1, 1\}$ (voir figure \ref{fig_ex_Z}), la matrice d'adjacence du graphe de l'automate $\Ared$ est
$$
	\begin{pmatrix}
		0 & 1 & 1 \\
		0 & 1 & 0 \\
		0 & 0 & 1
	\end{pmatrix}
$$
qui est idempotente.
Le nombre $c_n$ d'éléments de $\Gamma$ de longueur $n$ est donc $c_0 = 1$ pour $n=0$ et $c_n = 2$ pour $n > 0$ (il s'agit en effet des éléments $n$ et $-n$).
\end{ex}

Les 4 figures qui suivent donnent une structure automatique complète de l'exemple de l'introduction : 

\begin{figure}[H]
\centering
\caption{ }{\begin{minipage}{\linewidth} Automate $\Ared$ du semi-groupe engendré par les trois applications \[\left\{ \begin{matrix} 0 :& x & \mapsto &  x/3, & \\ 1 :& x & \mapsto & x/3 & + & 1, \\ 3 :& x & \mapsto & x/3 & + & 3, \end{matrix} \right. \] pour l'ordre lexicographique (avec $0 < 1 < 3$). \end{minipage}} \label{fig_ex_red}
\includegraphics[scale=.65]{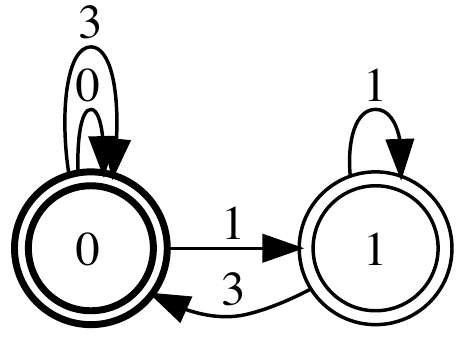}
\end{figure}

\noindent On voit que les mots réduits sont ici les mots ne contenant pas le sous-mot 10.

On voit sur la figure \ref{fig_ex_mul1} que le langage de $\A_1$ est formé des couples $(u1, u1)$ pour les mots $u$ réduits.
Sur la figure \ref{fig_ex_mul3}, on voit de même que le langage de $\A_3$ est formé des couples $(u3, u3)$ pour les mots $u$ réduits.
L'automate de la figure \ref{fig_ex_mul0} est plus compliqué puisque le mot $u0$ n'est pas réduit quand le mot réduit $u$ termine par un $1$ : son langage est l'ensemble des couples $(u1^n0, u03^n)$ pour les mots réduits $u$ et les entiers $n \in \N$.

\subsection{Fortement automatique implique automatique}

Voici un lien entre la structure fortement automatique et automatique :

\begin{prop} \label{afa}
Si le semi-groupe $\Gamma$ est fortement automatique, alors il est automatique.
\end{prop}

\emph{Preuve :}
Supposons que le semi-groupe $\Gamma$ soit fortement automatique : il existe donc un automate des relations $\Arel$, qui reconnait un langage rationnel $\Lrel$.
Par le point 1 des propriétés \ref{pptsfa}, le semi-groupe contient une relation entre des éléments de longueurs distinctes si et seulement si c'est un groupe fini.
Mais il est facile de voir qu'un groupe fini est automatique.
Ainsi, on supposera dans la suite de la preuve que les relations sont toujours entre des éléments de même longueur.
On a donc l'inclusion de l'ensemble des relations $\Lrel$ dans $(\Sigma \times \Sigma)^*$.

\paragraph{Construction de l'automate des mots réduits}
Considérons alors le langage rationnel
$$ L := L^{\operatorname{lex}} \cap \Lrel, $$
où $L^{\operatorname{lex}}$ est le langage de l'automate de la figure \ref{fig_lex}, en ayant muni l'alphabet $\Sigma$ d'un ordre total.
Posons alors
$$ \Lnonred := p_2(L), $$
le projecté du langage $L$ suivant la deuxième coordonnée.
Par la proposition \ref{proj_rat}, c'est un langage rationnel.
Alors $\Lnonred$ est l'ensemble des mots de $\Sigma^*$ qui sont non réduits pour l'ordre lexicographique.
En effet, le langage $L$ est l'ensemble des couples de mots équivalents $(u, v)$ de $\Sigma^*$ tels que $u$ est strictement inférieur à $v$ pour l'ordre lexicographique.
Ainsi, un mot $v$ est dans le langage $\Lnonred$ si et seulement s'il existe un mot $u$ équivalent et strictement inférieur dans l'ordre lexicographique.
En posant
$$ \Lred := \Sigma \backslash \Lnonred, $$
le langage $\Lred$ est donc un ensemble de mots réduits, et il est rationnel.

Ainsi, on a bien démontré l'existence de l'automate des mots réduits $\Ared$.

\paragraph{Construction des automates de multiplication}
Le paragraphe précédent a montré que le langage $\Lred$ des mots réduits est rationnel.
Pour $g \in \Sigma$, considérons alors le langage rationnel
$$ L^g := (\Lred g \times \Lred) \cap \Lrel. $$
C'est bien un langage rationnel par la proposition \ref{prod_rat}.
Pour obtenir un automate de multiplication $\A^g$, il suffit de considérer un automate $\A^g$ reconnaissant le langage $L^g$.

Ceci termine la preuve de la proposition \ref{afa}.
$\Box$

\subsection{Recherche du mot réduit}

Comment trouver le mot réduit correspondant à un mot donné ?
La structure automatique permet de faire cela :

\begin{prop}
Si le semi-groupe est automatique, alors il existe un algorithme quadratique qui prend en entrée un mot et rend le mot réduit correspondant.
\end{prop}

\emph{Preuve : }
Voir \cite{wpg}. \\

Ceci permet en particulier de résoudre le problème des mots (i.e. déterminer si deux mots donnés sont équivalents) en temps quadratique :

\begin{cor}
Si le semi-groupe est automatique, il existe un algorithme prenant en entrée deux mots et répondant en temps quadratique si les deux mots sont équivalents ou non.
\end{cor}

Lorsque le semi-groupe est fortement automatique, le problème des mots se résoud en temps linéaire et avec une mémoire constante puisqu'il est résolu par un automate.
Mais on peut aussi trouver le mot réduit correspondant à un mot donné rapidement :

\begin{prop}
Si le semi-groupe est fortement automatique, alors il existe un algorithme linéaire prenant en entrée un mot et rendant le mot réduit correspondant.
\end{prop}

\emph{Preuve :}
Définissons le langage
$$ L := \{ (a_1, b_1) ... (a_n, b_n) \in \Lrel | b_1 ... b_n \in \Lred e^* \} = (\Sigma^* \times \Lred e^*) \cap \Lrel. $$
Alors le langage $L$ est rationnel.

Soit $\A = (\Sigma \times (\Sigma \cup \{ e \}), Q, T, I, F)$ un automate déterministe reconnaissant le langage $L$.
Définissons alors l'automate $\A' = (\Sigma, Q', T', I', F')$ par :
\begin{itemize}
	\item $Q' := \Part(Q)$ (l'ensemble des parties de $Q$),
	\item $I' := \{I\}$,
	\item $F' := \{ P \in \Part(Q) | P \cap F \neq \emptyset \}$,
	\item $T'$ est défini par
	$$ (A, a, B) \in T' \text{ si et seulement si } B = \{ q \in Q | \text{Il existe } b \in (\Sigma \cup \{e\}), \quad{ et } p \in A \text{ tels que } (p, (a,b), q) \in T \}. $$
\end{itemize}
L'automate $\A'$ est clairement déterministe, et il reconnait le langage $\Sigma^*$, puisque pour tout mot $u \in \Sigma^*$, il existe un mot réduit $v \in \Sigma^*$ tel que l'on ait la relation $u = v$ dans le semi-groupe $\Gamma$, ce qui donne un mot $(u,ve^k)$ reconnu par l'automate $\A$ (en supposant que les mots réduits sont de longueur minimale).

Voici maintenant un algorithme permettant de trouver le mot réduit correspondant à un mot $u \in \Sigma^*$.
Considérons le chemin
$$A_0 \xrightarrow{u_1} A_1 \rightarrow ... \xrightarrow{u_n} A_n$$
dans l'automate $\A'$ étiqueté par le mot $u = u_1 ... u_n$.
Choisissons alors un état final $q_n \in A_n$ de l'automate $\A$.
Par définition, il existe alors une lettre $v_n \in (\Sigma \cup \{e\})$ et un état $q_{n-1} \in A_{n-1}$ tels que l'on ait la transition $q_{n-1} \xrightarrow{(u_n, v_n)} q_n$ dans l'automate $\A$. Et on peut trouver la lettre $v_n$ et l'état $q_{n-1}$ en temps constant. On peut alors continuer : on trouve une lettre $v_{n-1} \in (\Sigma \cup \{ e \})$ et un état $q_{n-2} \in A_{n-2}$ tels que l'on ait la transition $q_{n-2} \xrightarrow{(u_{n-1}, v_{n-1})} q_{n-1}$ dans l'automate $\A$, et ainsi de suite. On obtient finalement un chemin $q_0 \xrightarrow{(u_0, v_0)} q_1 \rightarrow ... \xrightarrow{(u_n, v_n)} q_n $ dans l'automate $\A$. L'état $q_0$ est un état initial de l'automate $\A$ puisque l'on a $q_0 \in A_0 \in I' = \{I\}$. On obtient donc un mot $v_1 ... v_n \in (\Sigma \cup \{e\})^*$ tel que le mot $(u_1, v_1) ... (u_n, v_n)$ est dans le langage $L$. Le mot $v = v_1 ... v_n$ est donc dans le langage $L_\Ared e^*$, et on a la relation $u_1 ... u_n = v_1 ... v_n$ dans le semi-groupe $\Gamma$. Ainsi, on obtient le mot réduit correspondant au mot $u$ en éliminant les lettres $e$ à la fin du mot $v$. Et tout ce calcul s'effectue en temps linéaire.
$\Box$

\section{Semi-groupes correspondant aux développements $\beta$-adique} \label{smgdb}

Soit $k$ un corps. Dans cette partie, on s'intéresse au semi-groupe $\Gamma$ engendré par les transformations affines
$$ x \mapsto \beta x + t $$
pout $t \in A \subset k$, où $A$ est une partie finie de $k$, et $\beta$ est un élément de $k$.

\begin{rem}
Si le corps $k$ est de caractéristique $0$, alors on peut supposer que l'on a $k=\C$. 
\end{rem}

Ce semi-groupe correspond au développement en base $\beta$, en utilisant l'ensemble de chiffres $A$.
Par exemple, l'exemple donné en introduction correspond au développement en base 3 en utilisant l'ensemble de chiffres $\{0, 1, 3\}$.

\subsection{Forte automaticité}

\begin{thm} \label{cm1}
Le semi-groupe $\Gamma$ est fortement automatique, sauf éventuellement dans le cas où le corps $k$ est de caractéristique nulle, et que le nombre complexe $\beta$ est algébrique, avec un conjugué de module 1.
\end{thm}

\emph{Preuve du théorème \ref{cm1} :}
Commençons par le cas où $\beta$ est une racine de l'unité. Par hypothèse, le corps $k$ est alors de caractéristique finie.
Dans ce cas, le semi-groupe $\Gamma$ est un groupe fini (et est donc fortement automatique par le point 1 des propriétés \ref{pptsfa}).
En effet, tous les générateurs sont d'ordre fini.

Supposons maintenant que $\beta$ n'est pas une racine de l'unité.
Par les propriétés \ref{pptsfa}, le semi-groupe est automatique si et seulement si il l'est pour la partie génératrice considérée ici, puisque les relations du semi-groupe $\Gamma$ sont entre des éléments de même longueur.
C'est-à-dire que si l'on a une égalité $a_1...a_n = b_1...b_n$ dans $\Gamma$ pour deux mots $a_1...a_n \in \Sigma^*$ et $b_1...b_n \in (\Sigma \cup \{e\})^*$, alors on a $b_1...b_n \in \Sigma^*$.
Ainsi, on peut considérer l'automate des relations $\Arel$ (à priori infini) sur l'alphabet $\Sigma \times \Sigma$ donné par la proposition \ref{osfa}, et on va montrer qu'il est fini (i.e. qu'il a un nombre fini d'états).

On peut représenter les états de l'automate $\Arel$ par des éléments du corps $k$, puisque les seuls éléments de $\Gamma^{-1} (\Gamma \cup \{e\})$ qui sont des états de l'automate $\Arel$ sont de la forme $x \mapsto x + t$ pour $t \in k$ (les autres éléments ne sont pas dans l'automate émondé).

Dans ce cadre particulier, quitte à remplacer $\Gamma$ par son inverse (ce qui ne change pas le fait qu'il soit fortement automatique), l'automate $\Arel$ est donc l'émondé de l'automate $\A = (\Sigma_{\A}, Q_{\A}, T_{\A}, I_{\A}, F_{\A})$ définit par :
\begin{enumerate}
	\item $\Sigma_{\A} = \Sigma \times \Sigma$,
	\item $Q_{\A} = k$,
	\item	$I_{\A} = \{ 0 \}$,
	\item $F_{\A} = \{ 0\}$,
	\item $T_{\A}$ est définit par :
		$(p, (g,h), q) \in T_{\A}  \text{ si et seulement si } q = \beta p + g-h$.
\end{enumerate}

\begin{rem}
J'aurait pu considérer l'automate $\A$ qui correspond directement à celui du semi-groupe $\Gamma$ et non pas à son inverse. Cela aurait donné les transitions $(p, (g,h), q) \in T_{\A}  \eq q = (p - g+h)/\beta$.
Mais la formule donnant les transitions de l'automate $\A$ pour le semi-groupe inverse me semblait plus agréable.
\end{rem}

Le lemme suivant fournit une critère algébrique d'appartenance à l'ensemble des sommets de l'automate des relations $\Arel$ :

\begin{lemme} \label{care}
	Un élément $x \in k$ est un état de l'automate des relations $\Arel$ (c'est-à-dire de l'émondé de l'automate $\A$) si et seulement si il existe deux polynômes $P, Q \in (A-A)[X]$ à coefficients dans $A-A$ tels que l'on ait $x = P(\beta) = \beta^{-1} Q(\beta^{-1})$.
\end{lemme}

\emph{Preuve du lemme :}
Un élément $x \in k$ est un état de l'automate des relations $\Arel$ si et seulement si il existe un chemin
$$0 \xrightarrow{(a_1, b_1)} ... \xrightarrow{(a_n, b_n)} x \xrightarrow{(a_k', b_k')} ... \xrightarrow{(a_0', b_0')} 0$$
dans l'automate $\A$, avec $a_i, b_i, a_i'$ et $b_i' \in A$.
On a alors $x = \sum_{i = 0}^{n-1} (a_{n-i} - b_{n-i}) \beta^{i}$ et $\beta^{k} x + \sum_{i=0}^{k-1} (a_i' - b_i') \beta^i = 0$, d'où $x = P(\beta) = \beta^{-1} Q(\beta^{-1})$, avec
$$P(X) = \sum_{i = 0}^{n-1} (a_{n-i} - b_{n-i}) X^{i} \quad \text{  et  } \quad X^{-1} Q(X^{-1}) = \sum_{i=0}^{k-1} (a_i' - b_i') X^{i-k}.$$
Et réciproquement, deux tels polynômes nous donnent un chemin de $0$ à $0$ passant par $x$.
$\Box$ \\

Notons $\F_p := \Frac(\Z/p\Z)$ le corps de fractions de $\Z/p\Z$ où $p$ est la caractéristique du corps $k$.
On a ainsi $\F_0 = \Q$ si le corps $k$ est de caractéristique nulle, et sinon $\F_p = \Z/p\Z$ est le corps fini à $p$ éléments.

Considérons $D$ le $\F_p(\beta)$-espace vectoriel engendré par la partie $A$, et soit $C$ une base de cet espace vectoriel.
Considérons un système de projecteurs $(p_c)_{c \in C}$ de $D$ dans $\F_p(\beta)$ tels que l'on ait
$$ \sum_{c \in C} p_c c = id. $$
Pour chaque $c \in C$, considérons alors le semi-groupe $\Gamma_c$ engendré par les transformations
$$ x \mapsto \beta x + t, \text{ pour } t \in p_c(A).$$
On a alors le lemme :

\begin{lemme} \label{proj}
Le semi-groupe $\Gamma$ est fortement automatique si et seulement si les semi-groupes $\Gamma_c$, $c \in C$, sont tous fortement automatiques.
\end{lemme}

\emph{Preuve :}
Chaque semi-groupe $\Gamma_c$ est fortement automatique si et seulement si il l'est pour la partie génératrice naturelle, de même que pour $\Gamma$, d'après les propriétés \ref{pptsfa}.

Montrons qu'un élement $x \in k$ est un état de l'automate $\Arel_\Gamma$ si et seulement si chaque projeté $p_c(x) \in \F_p(\beta)$ est dans l'automate $\Arel_{\Gamma_c}$.
Si $x \in k$ est un état de $\Arel_\Gamma$, alors par le lemme \ref{care} il existe deux polynômes $P$ et $Q$ à coefficients dans $A-A$ tels que $x = P(\beta) = \beta^{-1}Q(\beta^{-1})$.
En projetant, on obtient deux polynômes $p_c(P)$ et $p_c(Q)$ à coefficients dans $p_c(A-A) = p_c(A)-p_c(A)$ tels que $p_c(x) = p_c(P)(\beta) = \beta^{-1} p_c(Q)(\beta^{-1})$, donc $p_c(x)$ est un état de l'automate $\Arel_{p_c}$ par le lemme \ref{care}.

Et réciproquement, si l'on a $x = \sum_{c \in C} x_c c$ pour $x_c \in \F_p(\beta)$, avec $x_c$ état de l'automate $\Arel_{\Gamma_{c}}$, alors par le lemme \ref{care}, il existe des polynômes $P_c$ et $Q_c$ à coefficients dans $p_c(A) - p_c(A)$ tels que $x_c = P_c(\beta) = \beta^{-1} Q_c(\beta^{-1})$.
On a alors $x = P(\beta) = \beta^{-1}Q(\beta^{-1})$ en posant $P = \sum_{c \in C} P_c c$ et $Q = \sum_{c \in C} Q_c c$. Donc $x$ est un état de l'automate $\Arel$ par le lemme \ref{care}.
$\Box$ \\

Grâce au lemme \ref{proj}, on se ramène à ce que l'on ait $A \subseteq \F_p(\beta)$.
Quitte à multiplier la partie $A$ par un élément du corps $k$ (ce qui ne change pas le semi-groupe), on peut supposer que l'on a même $A \subseteq (\Z/p\Z)[\beta]$.
Ainsi, on peut supposer que l'ensemble des états de l'automate $\A$ est $(\Z/p\Z)[\beta]$. \\
Il y a maintenant deux cas à considérer :
\begin{itemize}
	\item \underline{Si $\beta$ est transcendant :}
		Dans ce cas, on a $(\Z/p\Z)[\beta] \simeq (\Z/p\Z)[X]$, et on peut voir $A$ comme une partie de $(\Z/p\Z)[X]$.
		Montrons alors que le degré des états de l'automate $\Arel$ (vus comme des polynômes) est strictement majoré par
		$$\max_{P \in A-A} \deg P.$$
		Soit $Q$ un état non nul de l'automate $\A$ de degré $\deg Q \geq \max_{P \in A-A} \deg P$.
		Pour toute transition $(Q, (U,V), R)$, on a alors
		$$\deg R = \deg (XQ + U - V) = 1 + \deg Q \geq \max_{P \in A-A} \deg P,$$
		car $\deg(U-V) \leq \max_{P \in A-A} \deg P < \deg (XQ)$.
		Ainsi, par récurrence, il ne peut pas exister de chemin de l'état $Q$ vers l'état final $0$, donc l'état $Q$ n'est pas dans l'automate émondé $\Arel$.
		
		Si le corps $k$ est de caractéristique non nulle, cela prouve donc que l'ensemble des états de l'automate $\Arel$ est fini.
		
		Supposons maintenant que le corps $k$ est de caractéristique nulle.
		Montrons alors que les polynômes $P \in \Z[X]$ qui sont des états de l'automate $\Arel$ ont leur $i^{\text{ème}}$ coefficient borné par
		$$ \sum_{j \geq i} \max_{\sum_{j} p_j X^j \in A-A} \abs{p_j}. $$
		Soit $Q = \sum_{j} q_j X^j$ un état de l'automate $\A$ ayant le $i^{\text{ème}}$ coefficient trop grand : $q_{i} > \sum_{j \geq i} \max_{\sum_{j} p_j X^j \in A-A} \abs{p_j}$, et soit $(Q, (U,V), R)$ une transition de l'automate $\A$. Alors le $(i+1)^{\text{ème}}$ coefficient du polynôme $R = XQ + U-V$ est trop grand :
		\begin{eqnarray}
			\abs{q_i + u_{i+1} - v_{i+1}}	&>& (\sum_{j \geq i} \max_{\sum_{j} p_j X^j \in A-A} \abs{p_j}) - \abs{u_{i+1} - v_{i+1}} \nonumber \\
									&\geq& (\sum_{j \geq i} \max_{\sum_{j} p_j X^j \in A-A} \abs{p_j}) - \max_{\sum_{j} p_j X^j \in A-A} \abs{p_{i+1}} \nonumber \\
									&=& \sum_{j \geq i+1} \max_{\sum_{j} p_j X^j \in A-A} \abs{p_j}. \nonumber
		\end{eqnarray}
		De la même façon que précédemment, il ne peut donc pas exister de chemin de l'état $Q$ vers l'état final $0$, donc l'état $Q$ n'est pas dans l'automate émondé $\Arel$.
		
		On a montré que les polynômes $P$ qui sont des états de l'automate $\Arel$ ont tous leur coefficients bornés et ont leur degré borné.
		On conclut donc que l'automate $\Arel$ est fini, et donc le semi-groupe $\Gamma$ est fortement automatique.
		
	\item \underline{Si $\beta$ est algébrique :}
		Alors dans le cas où le corps $k$ est de caractéristique non nulle, l'ensemble $(\Z/p\Z)[\beta]$ des états de $\A$ est fini, et donc l'automate $\Arel$ est aussi fini.
		
		Supposons donc que le corps $k$ est de caractéristique nulle. Sans perte de généralité, on supposera que l'on a $k = \Q(\beta)$. \\
		
		\begin{define}
			Soit $\Part$ l'ensemble (fini) des valeurs absolues $\abs{.}$ du corps $k = \Q(\beta)$ qui sont telles que $\abs{\beta} \neq 1$, et des valeurs absolues archimédiennes.
			On définit un anneau $\Ring$ stable par multiplication par $\beta$ et par $\beta^{-1}$, par
			$$ \Ring = \{ x \in k | \text{ Pour toute valeur absolue } v \not\in \Part, \abs{x}_v \leq 1 \}. $$
		\end{define}
				
		\begin{prop} \label{dccpt}
			L'anneau $\Ring$ est un réseau dans l'espace
			$$ E := \prod_{\abs{.} \in \Part} k_{\abs{.}}, $$
			dans lequel il est plongé diagonalement, où $k_{\abs{.}}$ est le complété du corps $k$ pour la valeur absolue $\abs{.}$. \\
		\end{prop}
		
		\emph{Preuve de la proposition \ref{dccpt} :}
		
		La discrétude de $\Ring$ dans $E$ est une conséquence de la formule du produit :
		
		\begin{prop}[Formule du produit] \text{ } \\
			Pour tout $x \in k \backslash \{0\}$, on a
			$$ \prod_{\abs{.} \in \Part_k} \abs{x} = 1, $$
			où $\Part_k$ est l'ensemble des valeurs absolues du corps $k$ (à équivalence près).
		\end{prop}
		\begin{rem}
			Dans la proposition précédente, on a choisi les valeurs absolues "canoniques" dans chaque classe d'équivalence.
			Voir \cite{SL}, V.1. pour plus de détails.
		\end{rem}
		
		Etant donné un point $x_0 \in \Ring$, et un ensemble non vide de valeurs absolues $\Part_0 \subseteq \Part$ contenant les valeurs absolues archimédiennes, la formule du produit nous donne que le voisinage
		$$ V := \{ x \in E | \text{ Pour toute valeur absolue } v \in \Part \backslash \Part_0, \abs{x - x_0}_v \leq 1, \text{ et pour toute valeur absolue } v \in \Part_0, \abs{x - x_0}_v < 1 \} $$
		a pour intersection $\{ x_0 \}$ avec $\Ring$, ce qui donne bien la discrétude de $\Ring$ dans $E$.
		
		Montrons maintenant la co-compacité de $\Ring$ dans $E$.
		
		Pour cela, on utilise le théorème :
		\begin{thm}  \text{ }\\
			Le corps $k$ est discret et co-compact dans l'ensemble des adèles $\Adele_k$. \\
		\end{thm}
		
		Voir \cite{SL} pour la définition des adèles et une preuve du résultat. \\
		L'espace $E' := E \times \prod_{\abs{.} \in \Part_k \backslash \Part} \ane_{k_{\abs{.}}}$ est une partie ouverte de l'ensemble des adèles $\Adele_k$, donc son image dans le quotient $\Adele_k/k$ est aussi ouverte.
		On en déduit qu'elle est aussi fermée, puisque l'on peut écrire l'orbite de $k$ sous l'action du groupe additif $E'$ comme l'union des autres $E'$-orbites.
		Ainsi, l'image de $E'$ dans le quotient $\Adele_k/k$ est compacte, d'où la co-compacité de $\Ring = k \cap E'$ dans $E'$ et donc dans $E$.
		
		$\Box$
		
		L'espace $E$ est le produit d'un nombre fini de corps $p$-adiques, et de copies de $\R$ et $\C$.
		
		\begin{ex}
			Pour $\beta = \frac{1+\sqrt{-14}}{5}$, l'espace $E$ est
			$$ E = \C \times \Q_3 \times \Q_5. $$
			
			Pour $\beta = \frac{\sqrt{-14}}{5}$, l'espace $E$ est
			$$ E = \C \times \Q_5 \times \Q_5 \times E_2 \times E_7, $$
			où $E_2$ et $E_7$ sont respectivement des extensions de degré 2 de $\Q_2$ et de $\Q_7$.
		\end{ex}
		
		\hfill
		
		On va montrer que l'ensemble des états de l'automate $\Arel$ est inclu dans une partie compacte de $E$, ce qui prouvera sa finitude.
		Soit $\abs{.} \in \Part$ une des valeurs absolues, et soit $\sigma : E \rightarrow k_{\abs{.}}$ la projection correspondante.
		Soit $\gamma = \sigma(\beta)$ le conjugué de $\beta$ correspondant.
		Montrons que les états $x$ de l'automate $\Arel$ vérifient
		$$ \abs{\sigma(x)} < \frac{1}{\abs{1-\abs{\gamma}}} \max_{a \in A-A} \abs{\sigma(a)}. $$
		Par hypothèse, on a $\abs{\gamma} \neq 1$.
		On a alors deux cas :
		\begin{enumerate}
			\item	$\abs{\gamma} < 1$ \\
				D'après le lemme \ref{care}, il existe un polynôme $P \in (A-A)[X]$ tel que l'on ait $x = P(\beta)$.
				On a alors
				$$\abs{\sigma(x)} = \abs{\sigma(P)(\gamma)} < \max_{a \in A-A} \abs{\sigma(a)} \sum_{i=0}^\infty \abs{\gamma}^i = \frac{1}{1 - \abs{\gamma}} \max_{a \in A-A} \abs{\sigma(a)}.  $$

			\item	$\abs{\gamma} > 1$ \\
				D'après le lemme \ref{care}, il existe un polynôme $Q \in (A-A)[X]$ tel que l'on ait $x = \beta^{-1}Q(\beta^{-1})$.
				On a alors
				$$\abs{\sigma(x)} = \abs{\gamma^{-1}\sigma(Q)(\gamma^{-1})} < \frac{1}{\abs{\gamma} - 1} \max_{a \in A-A} \abs{\sigma(a)}.  $$
		\end{enumerate}
		Le domaine de l'espace $E$ délimité par ces inégalités est relativement compact.
		Ainsi, la discrétude de l'anneau $\Ring$ dans l'espace $E$ entraîne que l'automate $\Arel$ n'a qu'un nombre fini d'états. Donc le semi-groupe $\Gamma$ est fortement automatique.
\end{itemize}
Ceci termine la preuve du théorème \ref{cm1}.
$\Box$

\begin{rem}
Dans les directions $p$-adiques, le fait que les valeurs absolues soient ultra-métriques permet d'obtenir les inégalités plus précises suivantes :
$$ \abs{x}_p \leq \max_{a \in A-A} \abs{a}_p \abs{\beta^{-1}}_p \text{ si } \abs{\beta}_p > 1,$$
$$ \abs{x}_p \leq \max_{a \in A-A} \abs{a}_p \text{ si } \abs{\beta}_p < 1,$$ 
pour tout état $x$ non nul de l'automate $\Arel$. \\
\end{rem}

\begin{rem}
La condition pour le nombre algébrique $\beta$ d'être sans conjugué de module 1 est nécessaire : voir exemple \ref{exs} et proposition \ref{rcm1}.
Cependant, il existe tout de même des nombres algébriques ayant au moins un conjugué de module 1 et pour lequels le semi-groupe est automatique.
Par exemple, le semi-groupe engendré par les deux applications
$$
\left\{
	\begin{array}{ccl}
		x &\mapsto& \beta x \\
		x &\mapsto& \beta x + 1
	\end{array}
\right.
$$
est libre (et donc fortement automatique) dès que le nombre $\beta$ a un conjugué de module strictement supérieur à 2.
Ainsi, par exemple, il est libre pour le nombre de Salem qui est racine du polynôme $X^4 - 3X^3 - 3X^2 - 3X + 1$.
\end{rem}

\subsection{Réciproque} 
Voici une réciproque au théorème \ref{cm1} :

\begin{prop} \label{rcm1}
En caractéristique nulle, si $\beta$ est un nombre algébrique ayant un conjugué de module $1$, alors il existe une partie finie $A \subset \Ring$ telle que le semi-groupe $\Gamma$ n'est pas fortement automatique.
\end{prop}

\begin{rem} \label{rrcm1}
D'après le lemme \ref{care}, la proposition \ref{rcm1} revient à dire que si $\beta$ a un conjugué de module 1, alors il existe une partie finie $A \subset \Ring$ telle que l'ensemble
$$ \{ x \in \Ring | \text{Il existe } P, Q \in (A-A)[X] \text{ tels que } x = P(\beta) = \beta^{-1}Q(\beta^{-1}) \} $$
est infini.
\end{rem}

\begin{rem}
	Dans la proposition \ref{cm1}, on peut même choisir $A \subset \Z[\beta]$.
\end{rem}

\begin{rem} \label{cinv}
Sous les hypothèses de la proposition, pour tout $\gamma$ conjugué de $\beta$, l'inverse $1/\gamma$ est aussi un conjugué de $\beta$.
\end{rem}
En effet, si $\gamma \in \C$ est de module 1, alors $1/\gamma$ est son conjugué complexe. \\ 

L'idée de la preuve de la proposition \ref{rcm1} est la suivante : grâce au lemme \ref{l1} on se ramène à seulement montrer l'existence de chemins jusqu'à 0, plutôt que dans les deux sens.
On choisit alors une partie finie $A$ et un domaine infini $\D$ qui soient tels que l'on puisse trouver des transitions des points de $\D$ vers d'autres points de $\D$ plus proche de 0, jusqu'à tomber dans une partie compacte. Il suffira alors de rajouter à la partie $A$ le bon ensemble fini de points pour obtenir des transitions de tous les points du compact vers 0. \\

\emph{Preuve de la proposition \ref{rcm1} :}
Soit $\beta \in \C$ un nombre algébrique ayant un conjugué de module 1.
Pour toute partie $A$ finie, on considèrera l'automate $\A$ définit dans la preuve du théorème \ref{cm1}.
D'après la proposition \ref{dccpt}, on peut plonger l'anneau $\Ring$ dans un espace $E$ qui est un produit de corps $p$-adiques et de copies des corps $\R$ et $\C$, de façon à ce que l'anneau $\Ring$ soit un réseau dans l'espace $E$.
On peut alors écrire l'espace $E$ comme un produit de trois espaces :
$$E = E_- \times E_0 \times E_+, $$
où \begin{itemize}
	\item l'espace $E_+$ est le produit des corps dans lesquels le nombre $\beta$ est de valeur absolue strictement supérieure à $1$,
	\item l'espace $E_0$ est le produit des corps archimédiens où le nombre $\beta$ est de module $1$,
	\item l'espace $E_-$ est le produit des corps dans lesquels le nombre $\beta$ est de valeur absolue strictement inférieure à $1$.
\end{itemize}
On notera respectivement $\Part_-$, $\Part_0$ et $\Part_+$ les ensembles de valeurs absolues des corps des espaces $E_-$, $E_0$ et $E_+$.
On notera aussi $\norm{.}_0$ la norme infinie sur $E_0$.

\begin{notation}
Etant donné $x = P(\beta) \in \Ring$, on note $\overline{x} := P(\beta^{-1}) \in \Ring$.
\end{notation}
L'application $x \mapsto \overline{x}$ est un élément de $\operatorname{Gal}(\Q(\beta) / \Q)$ d'après la remarque \ref{cinv}.

\begin{lemme} \label{l1}
S'il existe un chemin de $x$ à $0$ dans l'automate $\A$, alors il existe aussi un chemin de $0$ à $\overline{x}/\beta$.
\end{lemme}

\emph{Preuve du lemme :}
De même que dans la preuve du lemme \ref{care}, l'existence d'un chemin de $x$ à $0$ est équivalente à l'existence d'un polynôme $Q \in (A-A)[X]$ tel que $x = \beta^{-1}P(\beta^{-1})$.
On obtient alors un chemin de $0$ à $\overline{x}/\beta = Q(\beta)$ dans l'automate $\A$.
$\Box$ \\

On va montrer qu'il existe une partie finie $A \subset \Ring$ et un domaine $\D$ de l'espace $E$ tels que pour tout point $x$ de $\D \cap \Ring$, il existe un chemin de $x$ à $0$ dans l'automate $\A$.

\begin{rem} \label{car_co}
Ceci revient à démontrer qu'il existe une partie $A \subset \Ring$ finie et un domaine $\D$ de l'espace $E$ tels que
$$ \text{Pour tout } x \in \D \cap \Ring, \text{ il existe } Q \in (A-A)[X] \text{ tel que } x = \beta^{-1} Q(\beta^{-1}). $$
\end{rem}

\paragraph{Construction de la partie $A$} \text{ }\\
Pour $R > 0$, on définit une partie $A_R \subseteq \Ring$ par $x \in  A_R$ si et seulement si l'on a
$$ \abs{\sigma(x)} < R, $$
pour tout plongement $\sigma$. \\
On va maintenant fixer un réel $R$ assez grand de la façon suivante :

Soit $r$ le diamètre d'une maille du réseau $\Ring$ dans $E$.
Définissons une partie $K \subseteq E_0 \times E_+$ par $x \in K$ si et seulement si
$$ \text{pour toute valeur absolue } v \in \Part_0, \quad \abs{x}_v < 3r,$$
$$ \text{pour toute valeur absolue } v \in \Part_+, \quad \abs{x}_v < \abs{\beta}_vr.$$

En choisissant le rayon $R = 3 r \max_{\abs{.} \in \Part}{ \abs{\beta} }$, la partie $A_R$ est $r$-couvrante dans l'ensemble $K$.
C'est-à-dire tel que l'on a :
$$ K \subseteq \bigcup_{x \in A_R} B(x, r) \subseteq E_0 \times E_+. $$

On définit maintenant une partie $A' \subseteq \Ring$ par $x \in A'$ si et seulement si l'on a
$$ \text{pour toute valeur absolue } v \in \Part_+, \quad \abs{x}_v < \abs{\beta}_v r, $$
$$ \text{pour toute valeur absolue } v \in \Part_0, \quad \abs{x}_v < 3 r, $$
$$ \text{pour toute valeur absolue } v \in \Part_-, \quad \abs{x}_v < \frac{\max_{x \in A_R} \abs{x}_v}{1 - \abs{\beta}_v}, $$

La partie $A$ finie de $\Ring$ que l'on considère est
$$ A := A' \cup A_R. $$

\paragraph{Construction du domaine $\D$} \text{ }\\
On considère le domaine $\D \subseteq E$ définit par $x \in \D$ si et seulement si
$$ \abs{x}_v < r, \quad \text{pour toute valeur absolue } v \in \Part_+, $$
$$ \abs{x}_v < \frac{\max_{x \in A_R} \abs{x}_v}{1 - \abs{\beta}_v}, \quad \text{pour toute valeur absolue } v \in \Part_-. $$

\paragraph{Preuve de la non finitude de $\Arel$}

\begin{lemme} \label{l1r}
	Pour tout point $x \in \D$ tel que $\norm{x}_0 \geq 3 r$, il existe un point $y \in \D$ tel que
	\begin{itemize}
		\item Il existe une transition dans l'automate $\A$ de $x$ à $y$,
		\item On a l'inégalité $\norm{y}_0 < \norm{x}_0$.
	\end{itemize}
\end{lemme}
\emph{Preuve :}
Soit $x$ un point de $\D$ tel que $\norm{x}_0 \geq 3 r$.
Dans l'espace $E_0$, considérons une boule $B_0$ de centre $c$ sur le segment $[0,\beta x]$, de rayon $r$, et qui soit incluse dans la boule de centre 0 et de rayon $3r$.
Le point de coordonnées $b$ dans l'espace $E_0$ et $\beta x$ dans l'espace $E_+$ est dans le compact $K$ de l'espace $E_0 \times E_+$.
La partie $A_R$ étant $r$-couvrante, on peut alors trouver un point $t \in A_R$ qui est dans la boule $B_0$ dans l'espace $E_0$ et à distance au plus r de $x$ dans l'espace $E_+$.
Alors le point $y := \beta x - t$ convient car
\begin{enumerate}
	\item Il y a bien une transition de $x$ à $y$ dans l'automate $\A$ puisque l'on a
		$$-t \in -A \subseteq A-A.$$
	\item Dans l'espace $E_+$, on a bien $ \text{pour toute valeur absolue } v \in \Part_+, \abs{\beta x - t}_v < r $, puisque le point $\beta x$ est dans la boule de centre $t$ et de rayon $r$.
	\item Dans l'espace $E_0$, le point $\beta x - t$ est dans la boule de centre $\beta x - c$ et de rayon $r$. Le point $\beta x$ étant en dehors de la boule de centre 0 et de rayon $3r$, et le point $t$ étant dans la boule $B_0$ de centre $c$ et de rayon $r$ qui est incluse dans la boule de centre 0 et de rayon $3r$, on en déduit que la norme de l'élément $\beta x - t$ est strictement inférieure à celle de $x$ :
	$$ \norm{\beta x - t}_0 < \norm{\beta x}_0 = \norm{x}_0. $$
	\item Dans l'espace $E_-$, on a bien $\text{pour toute valeur absolue } v \in \Part_-$,
		$$
			\abs{\beta x - t}_v	\leq \abs{\beta}_v \abs{x}_v + \abs{t}_v
						< \frac{ \abs{\beta}_v \max_{a \in A_R} \abs{a}_v}{1 - \abs{\beta}_v} + \max_{a \in A_R} \abs{a}_v
						\leq \frac{ \max_{a \in A_R} \abs{a}_v}{1 - \abs{\beta}_v}.
		$$
\end{enumerate}
$\Box$

\begin{lemme} \label{l2r}
	Pour tout point $x \in \D \cap \Ring$, il existe un chemin de $x$ à $0$ dans l'automate $\A$.
\end{lemme}
\emph{Preuve :}
Soit $x \in \D \cap \Ring$.

Supposons d'abord que l'on ait $\norm{x}_0 < 3$.
Dans ce cas l'élément $\beta x$ est dans l'ensemble $A'$ (puisque l'on a $\norm{\beta x}_0 = \norm{x}_0 < 3$).
Il existe donc un élément $t \in -A' \subseteq -A \subseteq A-A$ tel que $\beta x + t = 0$, ce qui prouve l'existence d'une transition de $x$ vers $0$.

On est donc ramené à supposer que l'on ait $\norm{x}_0 \geq 3$.
Le lemme \ref{l1r} permet alors d'obtenir une transition vers un état de norme $\norm{.}_0$ strictement inférieure.
Par récurrence, et par discrétude de $\Ring$ dans $E$, on est donc ramené au premier cas.
$\Box$ \\

Pour achever la preuve de la proposition \ref{rcm1}, il suffit de remarquer que les lemmes \ref{l2r} et \ref{l1} entraînent que les éléments de l'ensemble infini
$$\D \cap \beta^{-1} \overline{\D} \cap \Ring$$
sont des états de l'automate émondé $\Arel$.
$\Box$

\subsection{Un exemple non fortement automatique}

La proposition \ref{rcm1} nous donne des parties $A$ finies pour lesquelles le semi-groupe n'est pas fortement automatique. Voici un exemple de semi-groupe non fortement automatique pour une partie $A = \{0, 1\}$ fixée.

\begin{prop} \label{exs}
Soit le nombre de Salem $\beta = \frac{1+\sqrt{2}+\sqrt{2\sqrt{2}-1}}{2} \simeq 1.8832035059$ (qui est une racine du polynôme $X^4-2X^3+X^2-2X+1$).
Alors, le semi-groupe engendré par les deux applications
$$
\left\{
	\begin{array}{ccl}
		x &\mapsto& \beta x \\
		x &\mapsto& \beta x + 1
	\end{array}
\right.
$$
n'est pas fortement automatique.
\end{prop}

\begin{cor}
Pour le nombre de Salem $\beta$ de la proposition précédente, il n'existe aucune partie $A \subset \C$ finie et de cardinal au moins 2, telle que le semi-groupe $\Gamma$ soit fortement automatique.
\end{cor}

\emph{Preuve du corollaire :}
Soit $A_0 \subseteq A$ une partie de $A$ de cardinal 2, et
soit $\Gamma_0$ le semi-groupe engendré par les deux applications
$$
	\begin{array}{ccl}
		x &\mapsto& \beta x + t, \text{ pour } t \in A_0.
	\end{array}
$$
Alors les sommets de l'automate $\Arel_{\Gamma_0}$ sont aussi des sommets de l'automate $\Arel_\Gamma$.
En effet, tout chemin de l'automate $\Arel_{\Gamma_0}$ est aussi un chemin de l'automate $\Arel_\Gamma$.
Ainsi, la forte automaticité du semi-groupe $\Gamma$ entraîne celle du semi-groupe $\Gamma_0$.
Or, le semi-groupe $\Gamma_0$ est le même que le semi-groupe de la proposition \ref{exs} modulo similitude (et donc le même d'un point de vue combinatoire).
Ainsi, la proposition \ref{exs} implique que le semi-groupe $\Gamma$ n'est pas fortement automatique.
$\Box$ \\

L'idée de la preuve de la proposition est la suivante : De même que dans la preuve de la proposition \ref{rcm1}, on se ramène à montrer l'existence de chemins des points d'un domaine infini $\D$ vers 0.
Et de la même façon, on commence par ramener tout point dans une partie compacte, puis on montre qu'il existe un chemin vers 0 pour chaque point de la partie compacte.

Ici il n'est pas possible de choisir la partie $A$ pour pouvoir rapprocher les points de $0$ en suivant une seule transition. Ce que l'on fera est de donner, suivant l'endroit où se trouve le point dans la partie contractante, des suites d'arêtes qui permettent de se rapprocher de 0 tout en restant bien dans le domaine $\D$. \\

\emph{Preuve de la proposition :}
Le réel $\beta$ est strictement supérieur à 1, et possède un conjugué réel de module strictement inférieur à 1, ainsi que deux conjugués complexes conjugués de module 1.
L'anneau des entiers $\Z[\beta] = \Ring$ de $\beta$ se plonge donc dans $\R^2 \times \C$, et on a deux plongements réels $\sigma_+$ et $\sigma_-$ respectivement dilatant et contractant et un plongement complexe $\sigma_0$ qui préserve le module.
Notons $\abs{.}_+$, $\abs{.}_-$ et $\abs{.}_0$ les valeurs absolues correspondantes.
La preuve du théorème \ref{cm1} nous permet de voir que les sommets de l'automate $\Arel$ sont dans un domaine de $\R^2 \times \C$ délimité par les inégalités
$$ \abs{x}_i <  \frac{1}{\abs{\abs{\beta}_i - 1}} , \text{ pour } i \in \{+, -\}. $$
Définissons alors le domaine $\D \subset \R^2 \times \C$ délimité par les inégalités
$$ x \in \D \text{ si et seulement si } \abs{x}_+ <  \frac{c}{\abs{\beta}_+ - 1} \text{ et } \abs{x}_- <  \frac{1}{1 - \abs{\beta}_-}, $$
où $0 < c < 1$ est une constante qui sera fixée ultérieurement. \\

De même que dans la preuve de la proposition \ref{rcm1}, on souhaite démontrer le lemme :
\begin{lemme}
	Pour tout point $x$ de $\D \cap \Z[\beta]$, il existe un chemin dans l'automate restreint au domaine $\D$ qui part de $x$ et aboutit en $0$.
\end{lemme}

\begin{rem}
Ceci revient à démontrer le résultat suivant :
$$ \text{Pour tout } x \in \D \cap \Z[\beta], \text{ il existe } Q \in \{-1, 0, 1\}[X] \text{ tel que } x = \beta^{-1} Q(\beta^{-1}). $$
\end{rem}

Pour cela, nous allons donner une stratégie qui permet, partant d'un point $x \in \D$, d'aboutir à un sommet $y \in \D$ tel que $\abs{y}_0 < \abs{x}_0$, en suivant des transitions de l'automate $\A$ données par une suite d'étiquettes dans $\{-1, 0, 1\}$. Nous allons donner cette stratégie de la façon suivante : étant donné un intervalle dans lequel se situe le point $(\beta-1) x$ dans la direction dilatante $E_+ = \R$, nous donnerons trois suites d'éléments de $\{-1, 0, 1\}$ qui donnent des chemins vers trois états dont l'un au moins sera de module strictement inférieur à $x$ dans la direction correspondant aux complexes conjugués.

Pour obtenir cela, nous avons découpé l'intervalle de la direction dilatante correspondant au domaine $\D$ en intervalles vérifiant le lemme :

\begin{lemme}
	Soit $I$ un intervalle de $E_+ = \R$. S'il existe trois suites $(a^1_i)_{1 \leq i \leq n_1}$, $(a^2_i)_{1 \leq i \leq n_2}$ et $(a^3_i)_{1 \leq i \leq n_3}$ de $\{-1, 0, 1\}^{(\N)}$ vérifiant :
	\begin{enumerate}
		\item Dans l'espace $E_0 = \C$, les trois nombres complexes
			$$c_j := \sum_{i = 1}^{n_j} a^j_i \beta^{-i} \quad \text{ pour } j \in \{1, 2, 3\}, $$
			forment un triangle contenant 0 dans son intérieur.
		\item Pour tout $j \in \{1, 2, 3\}$, on a l'inclusion
			$$ \beta^{n_j} I + \sum_{i = 1}^{n_j} a^j_i \beta^{n_j-i} \subseteq \D$$
	\end{enumerate}
	alors pour tout point $x \in \D$ tel que $\abs{x}_0$ est assez grand, il existe un chemin $x \xrightarrow{a^j_1} ... \xrightarrow{a^j_{n_j}} y $ pour un $j \in \{1,2,3\}$, vers un point $y \in \D$ tel que $\abs{y}_0 < \abs{x}_0$.
\end{lemme}

\emph{Preuve :}
La première condition permet de trouver un $j \in \{1, 2, 3\}$ tel que l'on ait l'inégalité
$$ \abs{x + c_j}_0 < \abs{x}_0, $$
dès que $\abs{x}_0$ est assez grand. \\
En effet, il suffit que dans l'espace $E_0 = \C$ le point $x$ soit en dehors du triangle formé par les médiatrices des segments $[0, -c_j], j=1,2,3$.

La deuxième condition assure que le point $y := \beta^{n_j} (x + c_j)$ est dans le domaine $\D$, et on a $\abs{y}_0 = \abs{x + c_j}_0$. Et pour finir, la définition de $y$ donne l'existence du chemin
$$x \xrightarrow{a^j_1} ... \xrightarrow{a^j_{n_j}} y$$
dans l'automate $\A$.
$\Box$

Voici cette stratégie, pour $c = 0.883204$, en supposant que l'on parte d'un point $x \in \D$ tel que $\abs{x}_0 > 4.978154$ et $\sigma_+(x) \geq 0$ :

\begin{tabular}{c | c | c}

%
%

%


[0.468990, 0.601232]	& [0.601232, 0.647807]	& [0.647807, 0.671454] \\
-1 0 -1 1				& -1 -1 0 1 1			& -1 -1 0 0 1 1   \\
-1 0 -1				& -1 -1 1 -1 1			& -1 -1 0 1 0 -1   \\
-1 0 0				& -1 0 -1 -1 1			& -1 0 -1 -1 1 -1   \\
& & \\

[0.671454, 0.685095]	& [0.685095, 0.708742]	& [0.708742, 0.718028] \\
-1 0 -1 -1 -1 1			& -1 -1 0 0 0 1			& -1 -1 0 0 1 -1   \\
-1 -1 0 1 0 -1			& -1 -1 0 0 1 -1			& -1 0 -1 -1 -1 0   \\
-1 0 -1 -1 0 -1			& -1 0 -1 -1 0 -1			& -1 0 -1 -1   \\
& & \\

[0.718028, 0.780048]	& [0.780048, 0.812982]	& [0.812982, 0.832782] \\
-1 -1 -1 1 1			& -1 -1 -1 0 1 1			& -1 -1 -1 0   \\
-1 -1 -1				& -1 -1 -1 0 1			& -1 -1 -1 0 1 0 -1   \\
-1 -1 0 -1 1			& -1 -1 -1 1 -1 1			& -1   \\
& & \\

[0.832782, 0.850270]	& [0.850270, 0.883204]	& [-0.468990, 0.468990] \\
-1 -1 -1 0 0 1 -1			& -1 -1 -1 -1 1			& $0^n$ \\
-1 -1 -1 0 0 0			& -1 -1 -1 -1			& \\
-1 -1 -1 0 1 -1 -1		& -1 -1 -1 0 0			& \\

\end{tabular} \\

Si l'on part d'un point $x \in \D$ tel que $\abs{x}_0 > 4.978154$ et $\sigma_+(x) < 0$, alors on déduit la stratégie de celle donnée ci-dessus : il suffit de faire l'opposé de la suite de coups de l'intervalle opposé. \\
Quand on arrive dans l'intervalle [-0.468990, 0.468990], il suffit de suivre suffisamment d'arêtes étiquetées par $0$ pour retomber dans l'un des intervalles ci-dessus ou son opposé.


Il reste ensuite à donner la stratégie pour $\abs{x}_0 \leq 4.978154$. Le nombre de points est alors fini : il y en a 76.
Voici une stratégie pour 38 de de ces points :

\begin{tabular}{l | l}

\noindent
$1 \rightarrow$ -1 -1 -1 -1 1  &
$\beta-2 \rightarrow$ 1 -1  \\
$\beta-1 \rightarrow$ -1 -1 -1  &
$2\beta-3 \rightarrow$ -1 -1  \\
$\beta^2-3\beta+1 \rightarrow$ 1 1 1  &
$\beta^2-3\beta+2 \rightarrow$ 0 1  \\
$\beta^2-3\beta+3 \rightarrow$ -1  &
$\beta^2-2\beta \rightarrow$ 0 0 1  \\
$\beta^2-2\beta+1 \rightarrow$ -1  &
$\beta^2-\beta-1 \rightarrow$ -1  \\
$2\beta^2-4\beta \rightarrow$ 0 1 & 
$2\beta^2-4\beta+1 \rightarrow$ -1  \\
$2\beta^2-3\beta-1 \rightarrow$ -1  &
$\beta^3-3\beta^2+\beta+1 \rightarrow$ 1  \\
$\beta^3-3\beta^2+2\beta \rightarrow$ 1  &
$\beta^3-3\beta^2+2\beta+1 \rightarrow$ -1 -1 -1 1 1  \\
$\beta^3-3\beta^2+3\beta-2 \rightarrow$ 1  &
$\beta^3-3\beta^2+3\beta-1 \rightarrow$ -1 -1 0 1 1 -1  \\
$\beta^3-3\beta^2+4\beta-3 \rightarrow$ -1 0 -1  &
$\beta^3-2\beta^2-\beta+2 \rightarrow$ 0 1  \\
$\beta^3-2\beta^2 \rightarrow$ 0 1  &
$\beta^3-2\beta^2+1 \rightarrow$ -1 0 0 -1 1 \\
$\beta^3-2\beta^2+\beta-2 \rightarrow$ 1  &
$\beta^3-2\beta^2+\beta-1 \rightarrow$ 0 -1 -1 -1 1  \\
$\beta^3-\beta^2-2\beta \rightarrow$ 1 0 1 -1  &
$\beta^3-\beta^2-2\beta+1 \rightarrow$ -1  \\

\end{tabular}

\begin{tabular}{l | l}

$\beta^3-\beta^2-\beta-1 \rightarrow$ -1  &
$\beta^3-4\beta \rightarrow$ 1 1 1  \\
$\beta^3-3\beta \rightarrow$ -1  &
$2\beta^3-5\beta^2+3\beta-2 \rightarrow$ 1  \\
$2\beta^3-5\beta^2+4\beta-3 \rightarrow$ -1 1 0 1 1  &
$2\beta^3-4\beta^2 \rightarrow$ 1  \\
$2\beta^3-4\beta^2+\beta-1 \rightarrow$ 1 -1 -1 -1  &
$2\beta^3-4\beta^2+\beta \rightarrow$ -1 -1 -1 -1 -1 1 0 1 1  \\
$2\beta^3-4\beta^2+2\beta-2 \rightarrow$ -1 -1 -1 0 1  &
$2\beta^3-3\beta^2-2\beta \rightarrow$ 1 1  \\
$2\beta^3-3\beta^2-\beta-1 \rightarrow$ 0 1  &
$2\beta^3-3\beta^2-\beta \rightarrow$ -1 -1 0 -1  \\

\end{tabular} \\

De la même façon que précédemment, on déduit la stratégie pour l'autre moitié des points en considérant la suite de coups opposée de l'élément de $\Z[\beta]$ opposé.

En suivant cette stratégie, on aboutit à l'état $0$ en partant de n'importe quel point du domaine $\D$.
Il suffit donc de vérifier chacun des cas ci-dessus pour obtenir une preuve de l'exemple.
Tout ceci a été vérifié par ordinateur.
$\Box$

\begin{ex}
En suivant la stratégie ci-dessus, en partant de $x = \beta - 2$, on obtient le chemin
$$ x = \beta - 2 \xrightarrow{\text{ 1 -1 }} \beta^3 - 2 \beta^2 + \beta -1 \xrightarrow{\text{ 0 -1 -1 -1 1 }} 0. $$
\end{ex}
\begin{ex}
En partant de $x = 2 \beta^3 - 5 \beta^2 + \beta + 2$, on a $(\beta - 1)x \in [-0.468990, 0.468990]$ (et on a bien $x \in \D$).
Et on a $\abs{x}_0 \simeq 7.6615$.
On obtient le chemin
$$ x = 2 \beta^3 - 5 \beta^2 + \beta + 2 \xrightarrow{\text{ 0 }} -\beta^3 - \beta^2 + 6 \beta - 2, $$
ce qui nous fait tomber dans l'intervalle $(\beta - 1)(-\beta^3 - \beta^2 + 6 \beta - 2) \in [-0.832782, -0.812982]$.
On a alors le chemin
$$ -\beta^3 - \beta^2 + 6 \beta - 2 \xrightarrow{\text{ 1 }} -3 \beta^3 + 7 \beta^2 - 4 \beta + 2 , $$
et on tombe dans l'intervalle $(\beta - 1)(-3 \beta^3 + 7 \beta^2 - 4 \beta + 2) \in [-0.667607, -0.634165]$,
avec $\abs{-3 \beta^3 + 7 \beta^2 - 4 \beta + 2}_0 \simeq 6.8587$.
On a ensuite
$$ -3 \beta^3 + 7 \beta^2 - 4 \beta + 2 \xrightarrow{\text{ 1 0 1 1 -1 1 }} -\beta^3 + 4 \beta - 1, $$
avec $\abs{ -\beta^3 + 4 \beta - 1}_0 \simeq 5.3049$. Et l'on tombe dans l'intervalle $[-0.468990, 0.468990]$.
On continue avec
$$ -\beta^3 + 4 \beta - 1 \xrightarrow{\text{ 0 0 0 }} \beta^3 - 4 \beta^2 + 4 \beta - 1, $$
qui nous met dans l'intervalle $[-0.883204, -0.850270]$.
Puis
$$ \beta^3 - 4 \beta^2 + 4 \beta - 1 \xrightarrow{\text{ 1 1 1 1 }} -\beta^3 + 4 \beta, $$
avec $\abs{-\beta^3 + 4 \beta}_0 \simeq 4.9309$.
Comme on a $\abs{-\beta^3 + 4 \beta}_0 \leq 4.978154$, l'élément $-\beta^3 + 4 \beta$ est dans le tableau.
Il suffit donc de suivre le chemin indiqué :
$$ -\beta^3 + 4 \beta \xrightarrow{\text{ -1 -1 -1 }} 2 \beta^3 - 5 \beta^2 + 3 \beta - 2 \xrightarrow{\text{ 1 }} - \beta^3 + \beta^2 + 2 \beta - 1 $$
$$ \xrightarrow{\text{ 1 }} -\beta^3 + 3 \beta^2 - 3 \beta + 2 \xrightarrow{\text{ -1 }} \beta^3 - 2 \beta^2 \xrightarrow{\text{ 0 1 }} - \beta^3 + 2 \beta^2 - \beta + 1 \xrightarrow{\text{ 0 1 1 1 -1}} 0. $$ 

On vérifie que l'on a en effet bien la relation
$$ \beta^{-2}+\beta^{-3}+\beta^{-5}+\beta^{-6}-\beta^{-7}+\beta^{-8}+\beta^{-12}+\beta^{-13}+\beta^{-14}+\beta^{-15}-\beta^{-16}-\beta^{-17}-\beta^{-18}$$
$$+\beta^{-19}+\beta^{-20}-\beta^{-21}+\beta^{-23}+\beta^{-25}+\beta^{-26}+\beta^{-27}-\beta^{-28} =  2 \beta^3 - 5\beta^2 + \beta + 2. $$

\end{ex}

\begin{quest}
Le semi-groupe de l'exemple \ref{exs} est-il automatique ?
A t'il un ensemble de mots réduits qui soit un langage rationnel ? 
A t'il une présentation finie ?
Je conjecture une réponse négative à ces questions.
\end{quest}

\section{Exemples} \label{Exemples}

\subsection{Le cas où $\beta$ est un nombre de Pisot}
On appelle \defi{nombre de Pisot} un réel algébrique $\beta > 1$ qui a tous ses conjugués de modules strictement inférieurs à 1.
Dans son article \cite{Lal}, Lalley s'intéresse aux semi-groupes  engendrés par les transformations affines
$$ x \mapsto \beta x + t $$
pout $t \in A \subset \Z[\beta]$, où $1/\beta$ est un nombre de Pisot et $A$ est une partie finie de $\Z[\beta]$.
D'après le théorème \ref{cm1}, le semi-groupe est fortement automatique. En particulier, il existe un automate des mots réduits pour l'ordre lexicographique, par la proposition \ref{afa}.
Lalley donne une construction de l'automate des mots réduits (mais sans parler d'automates), mais la construction que je propose dans cet article est différente. Il obtient ainsi la vitesse de croissance du semi-groupe, qu'il arrive à relier à la dimension de Hausdorff de l'ensemble limite.

\subsection{L'exemple de Kenyon}
R. Kenyon étudie en détails dans son article \cite{ken} le semi-groupe engendré par les trois transformations
$$
\left\{
	\begin{array}{cccc}
		0: x &\mapsto& x/3& \\ 
		b: x &\mapsto& x/3& + t \\ 
		1: x &\mapsto& x/3& + 1 \\ 
	\end{array}
\right.
$$
où $0 < t < 1$ est un réel.

D'après le théorème \ref{cm1}, ce semi-groupe est fortement automatique pour tout réel $t$. 
La proposition suivante permet de savoir si le semi-groupe $\Gamma$ est libre ou non (et donc de savoir si l'automate des relations $\Arel$ est trivial ou non).
\begin{prop}[Kenyon] \label{pken}
Le semi-groupe $\Gamma$ est libre si et seulement si le réel $t$ n'est pas un rationnel de la forme $\frac{p}{q}$ avec $p+q \not\equiv 0 \mod 3$, pour $p \wedge q = 1$.
\end{prop}

Dans son article, Kenyon propose un construction d'automates qui est un cas particulier de la construction que je propose dans cet article. Il s'intéresse à la dimension de Hausdorff de l'ensemble limite du semi-groupe $\Gamma$. Il montre que celle-ci est reliée à la vitesse exponentielle de croissance du semi-groupe (que l'on peut calculer avec l'automate des mots réduits) quand le paramètre de translation $t$ est rationnel. C'est encore une conjecture (connue sous le nom de conjecture de Furstenberg) que l'ensemble limite est de dimension de Hausdorff 1 quand le paramètre de translation $t$ est irrationnel. Voir \cite{ken}.

Voici quelques exemples d'automates des relations que l'on obtient pour le semi-groupe engendré par les 3 transformations
$$ \left\{ \begin{matrix} 0 :& x & \mapsto &  x/3, & \\ a :& x & \mapsto & x/3 & + & a, \\ b :& x & \mapsto & x/3 & + & b. \end{matrix} \right. $$
(Ce qui revient à considérer le semi-groupe de Kenyon avec $t = a/b$.)
Cela permet de voir à quoi ressemble les premiers exemples d'automates correspondant au semi-groupe de Kenyon.

\begin{center}
	\begin{tabular}{cc}
		\includegraphics[scale=.5]{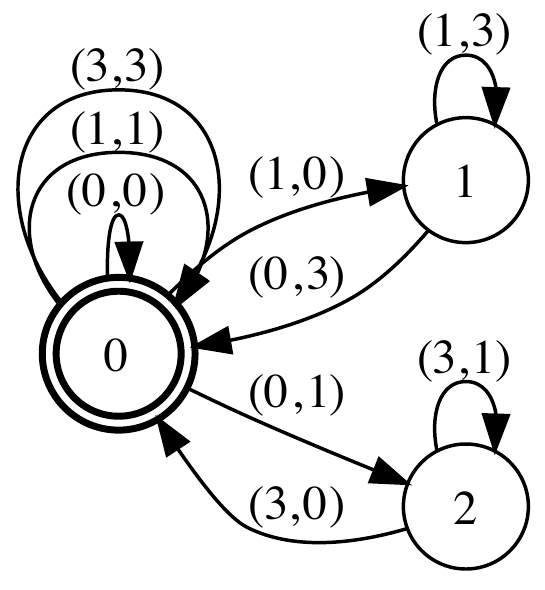} &
		\includegraphics[scale=.5]{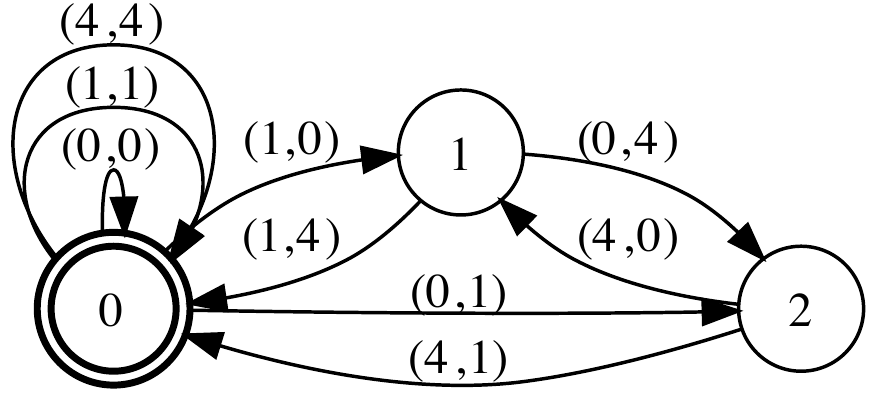} \\
		a=1 et b=3 & a=1 et b=4 \\ \\
	\end{tabular}
	\begin{tabular}{cc}
		\includegraphics[scale=.5]{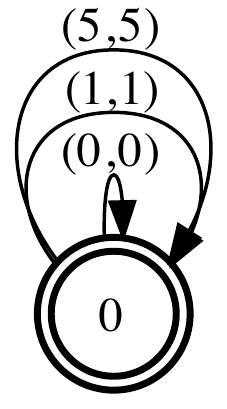} &
		\includegraphics[scale=.4]{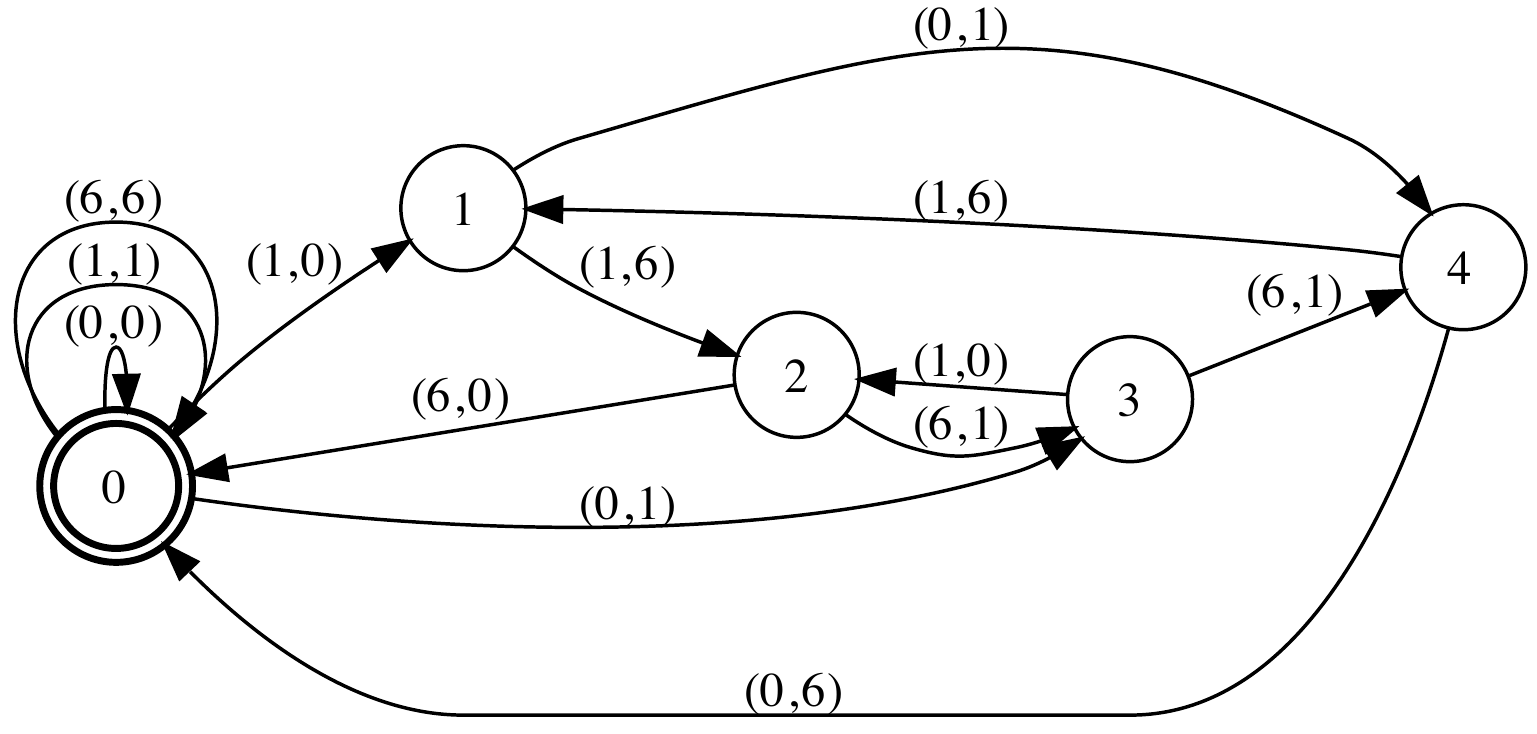} \\
		a=1 et b=5 & a=1 et b=6 \\ \\
		
		\includegraphics[scale=.3]{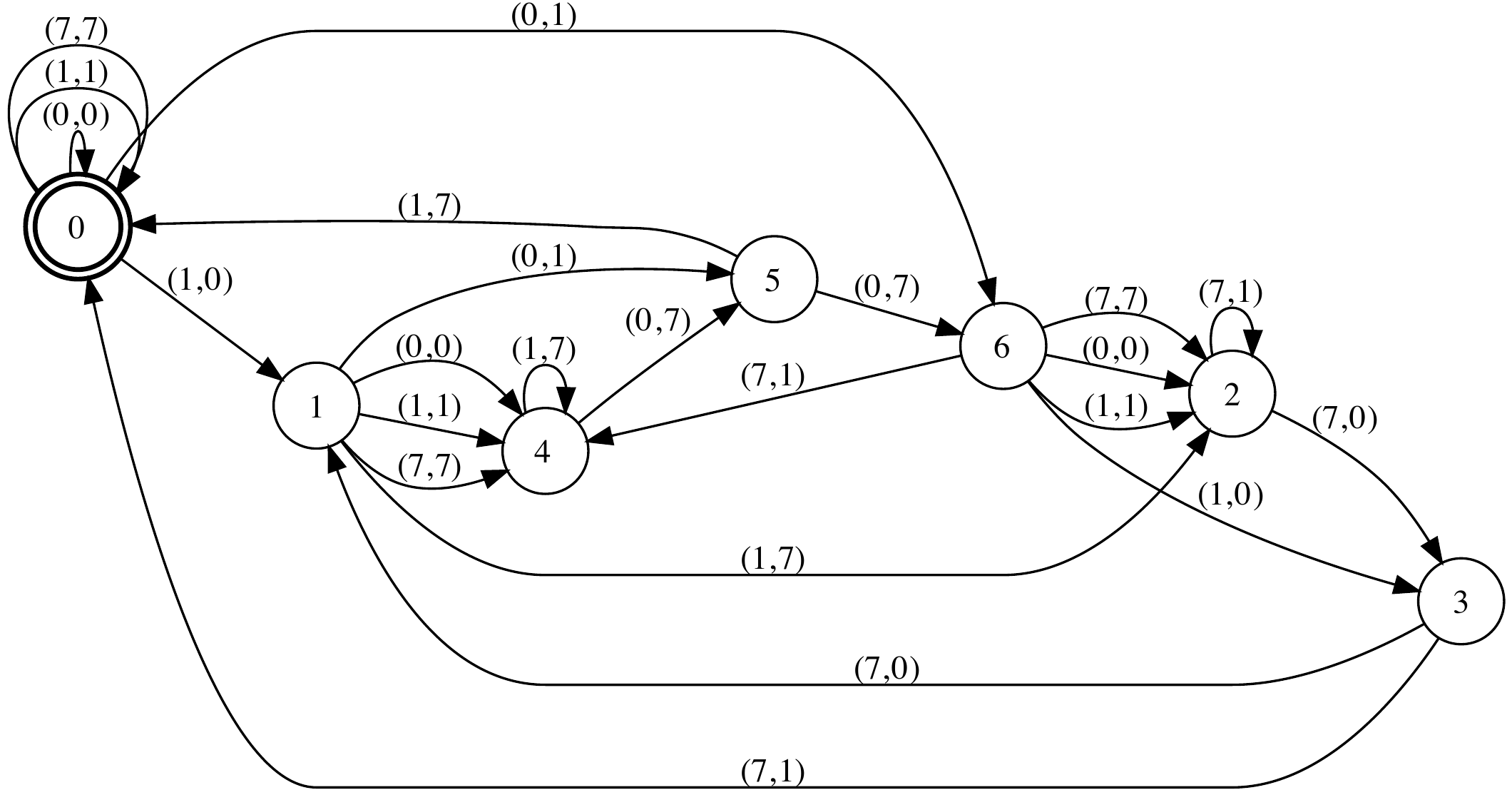} &
		\includegraphics[scale=.3]{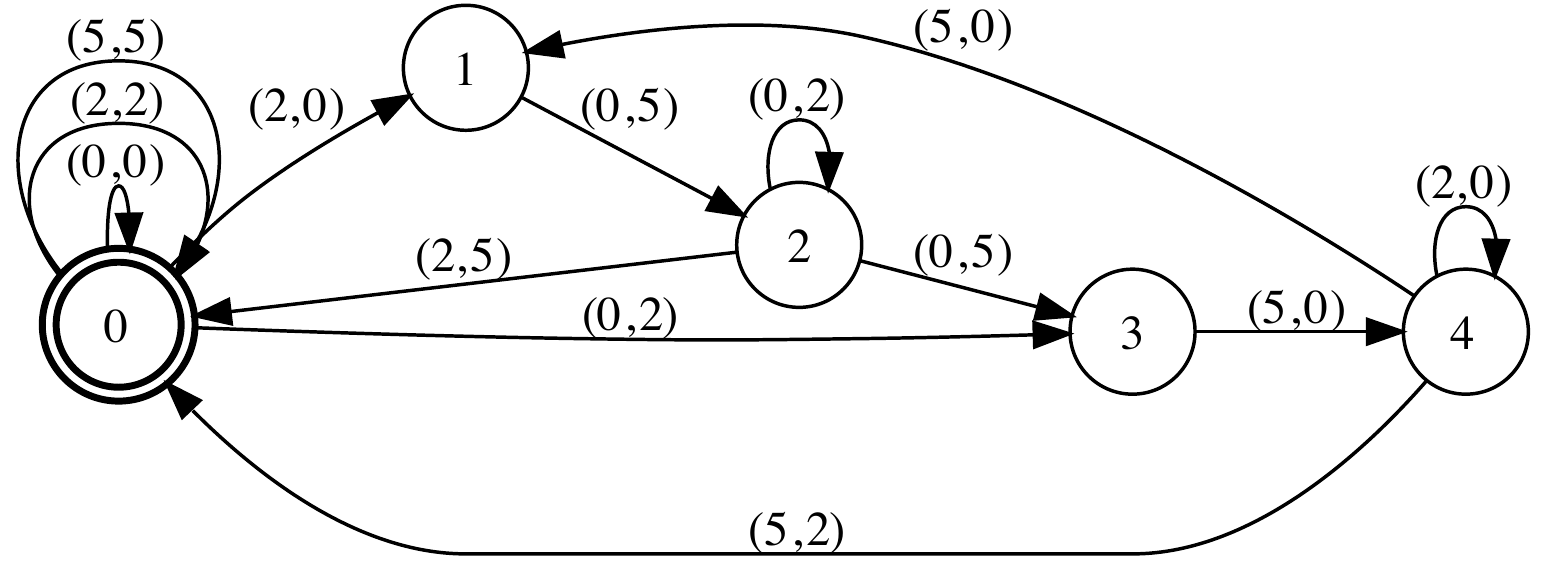} \\
		a=1 et b=7 & a=2 et b=5 \\ \\
	\end{tabular}
\end{center}

\begin{center}
	\begin{tabular}{cc}
		\includegraphics[scale=.5]{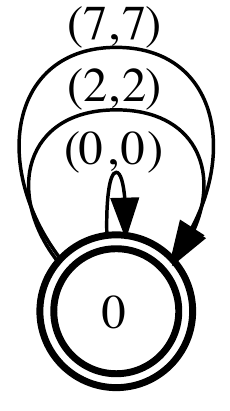} &
		\includegraphics[scale=.3]{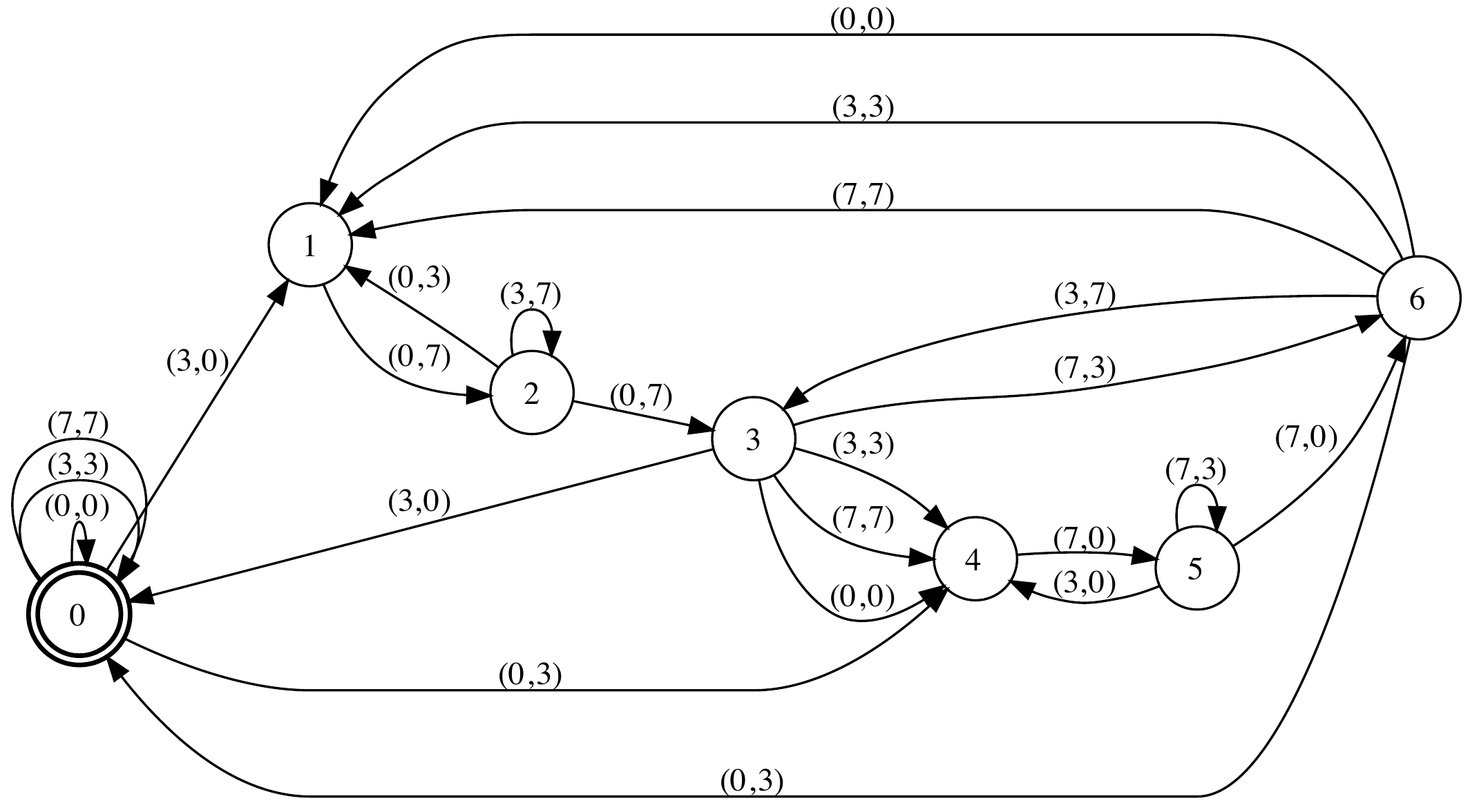} \\
		a=2 et b=7 & a=3 et b=7
	\end{tabular}
\end{center} \text{ } \\ \\
	
\noindent Et voici les automates minimaux des mots non réduits pour les même exemples (sauf pour $a/b = 1/5$ et $a/b = 2/7$ puisque les semi-groupes correspondant sont libres) :
\begin{center}
	\begin{tabular}{cc}
		\includegraphics[scale=.5]{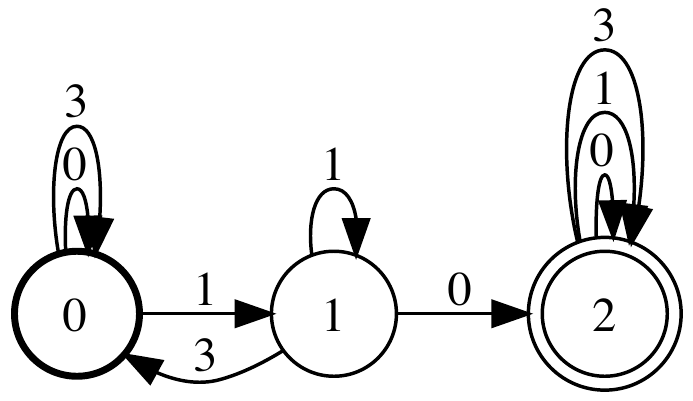} &
		\includegraphics[scale=.5]{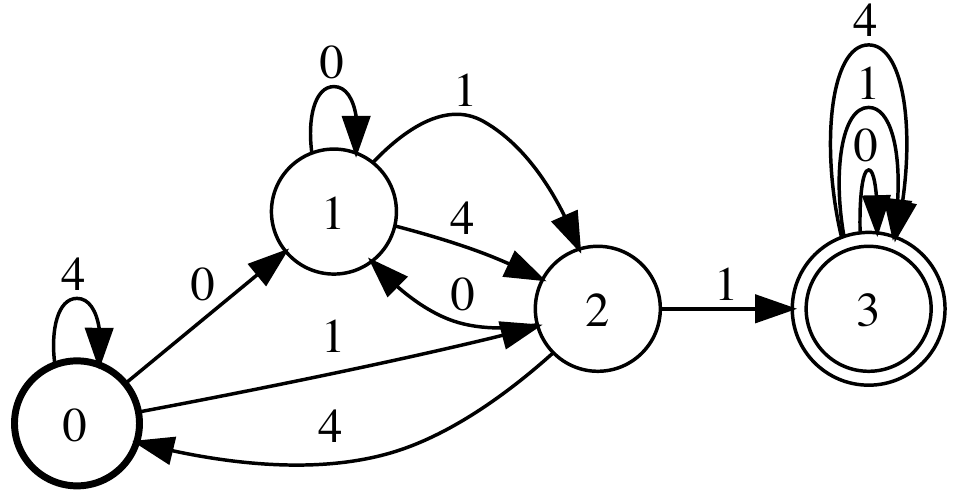} \\
		a=1 et b=3 & a=1 et b=4 \\ \\
	\end{tabular}
	\begin{tabular}{cc}
		\includegraphics[scale=.4]{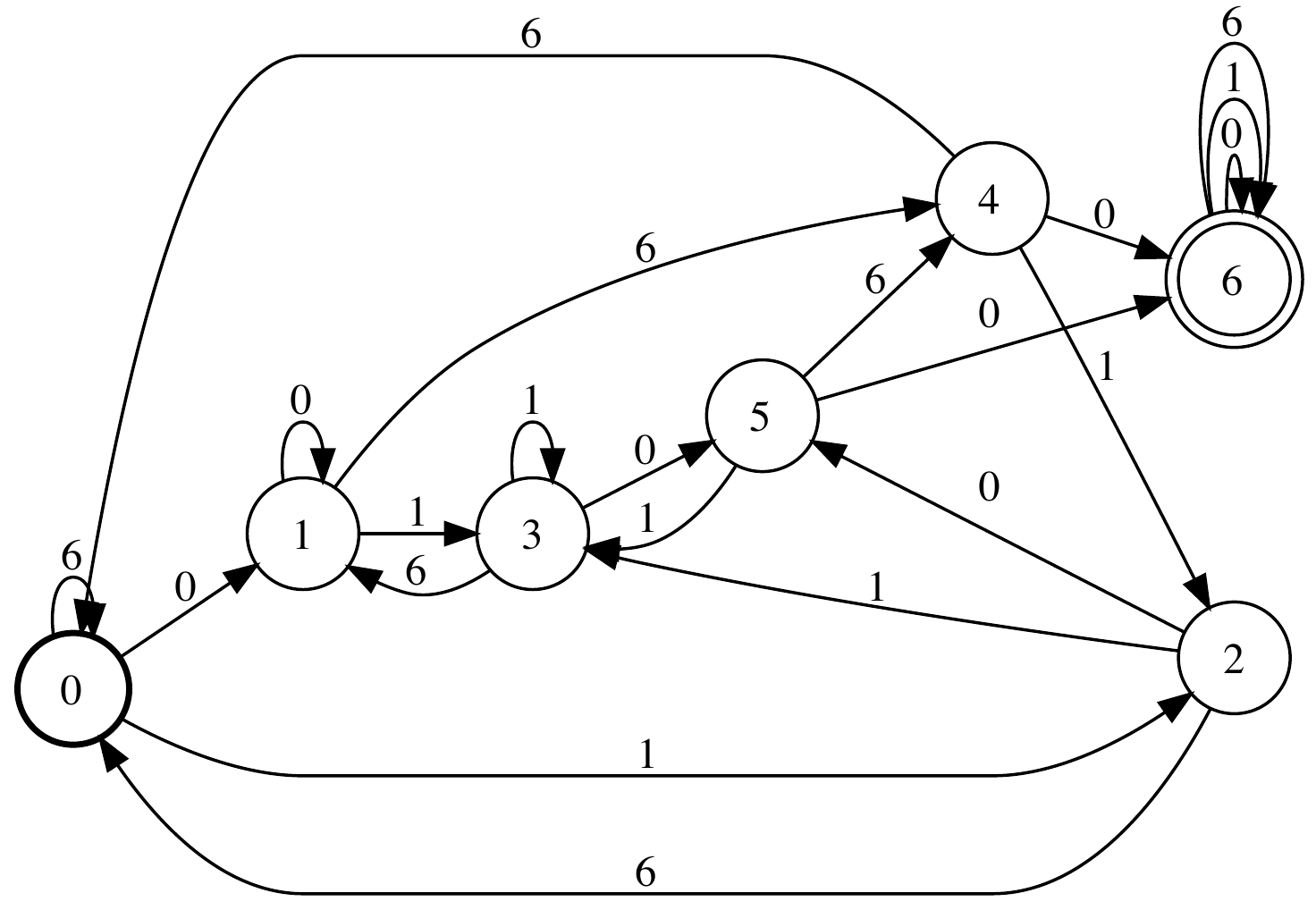} &
		\includegraphics[scale=.3]{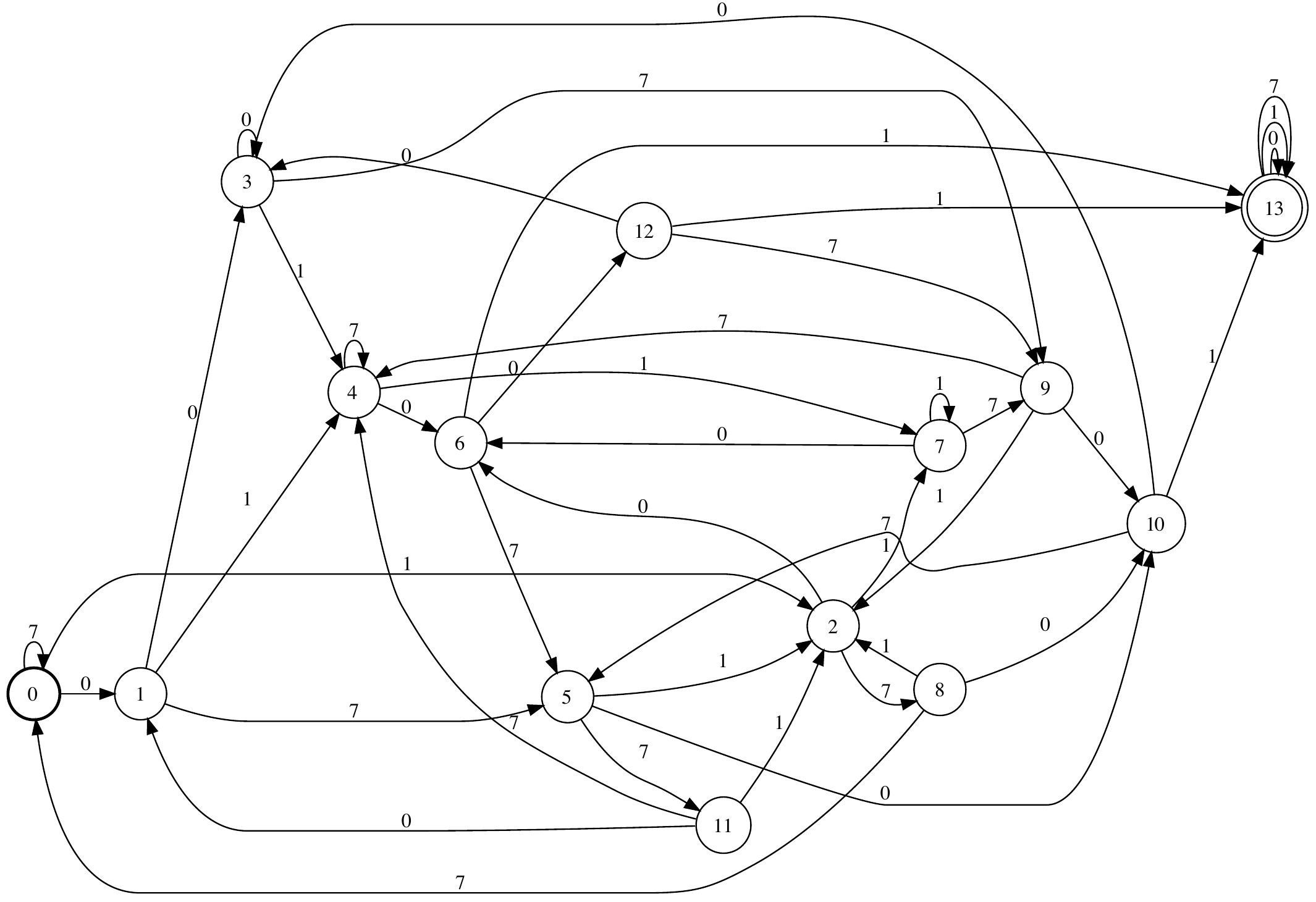} \\
		a=1 et b=6 & a=1 et b=7 \\ \\
	\end{tabular}
\end{center}

\begin{center}
	\begin{tabular}{cc}
		\includegraphics[scale=.3]{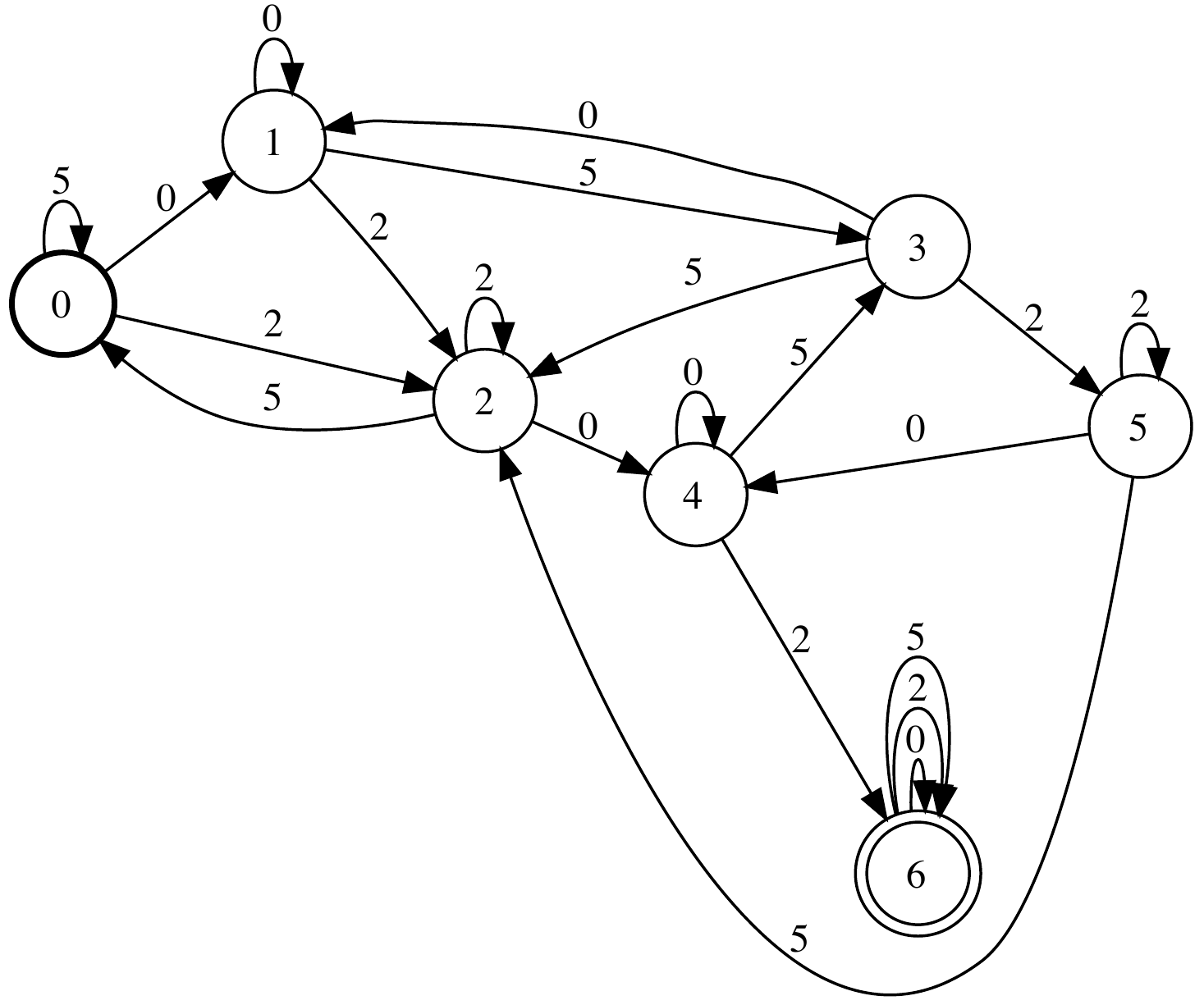} &
		\includegraphics[scale=.3]{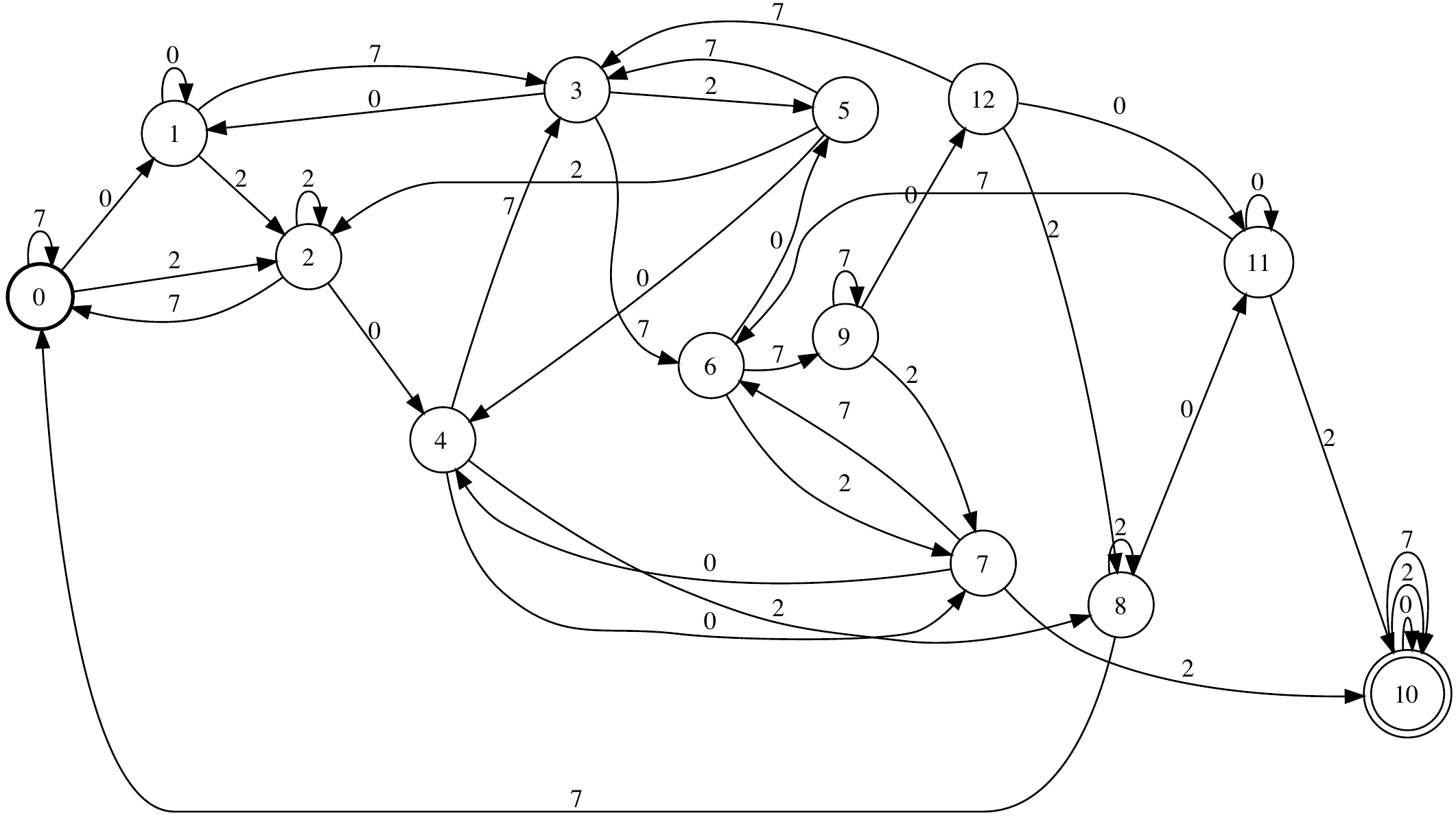} \\
		 a=2 et b=5 & a=3 et b=7		
		
	\end{tabular}
\end{center}

Voici la vitesse exponentielle de croissance $\lambda$ du semi-groupe de Kenyon pour quelques valeurs du paramètre de translation $t$.
On note $\pi_\lambda$ le polynôme minimal de $\lambda$.

\begin{center}
	\begin{tabular}{| c | c | c |}
		\hline
		t & $\lambda$ & $\pi_\lambda$ \\
		\hline
		1/3 & 2.6180	& $x^2-3x+1$ \\
		1/4 & 2.6180	& $x^2-3x+1$ \\
		2/5 & 2.8019	& $x^3-4x^2+3x+1$ \\
		1/6 & 2.7321	& $x^2-2x-2$ \\
		1/7 & 2.7383	& $x^5-3x^4+x^2+3x-1$ \\
		3/7 & 2.8794	& $x^3-3x^2+1$ \\
		3/8 & 2.8136	& $x^3-2x^2-3x+2$ \\
		1/9 & 2.6180	& $x^2-3x+1$ \\
		2/9 & 2.7233	& $x^6-3x^5+x^3+3x^2-1$ \\
		4/9 & 2.8794	& $x^3-3x^2+1$ \\
		1/10 & 2.6180	& $x^2-3x+1$ \\
		3/10 & 2.7699	& $x^6-2x^5-4x^4+x^3+9x^2+6x+3$ \\
		2/11 & 2.7421	& $x^5-4x^4+3x^3+x^2+x-1$ \\
		3/11 & 2.8073	& $x^9-4x^8+x^7+7x^6-9x^3-x^2+3x-1$ \\
		5/11 & 2.9242	& $x^{11}-5x^{10}+6x^9-7x^8+32x^7-32x^6+15x^5-49x^4+$ \\ && $20x^3-13x^2+3x-1$ \\
		\hline
	\end{tabular}
\end{center}

J'ai donné ici tous les rationnels $a/b$ dans l'intervalle $]0, 1/2[$ ayant un dénominateur inférieur ou égal à 11, avec $a+b \not\equiv 0 \mod 3$. On en déduit facilement les valeurs des vitesses exponentielles de croissance pour tous les rationnels $a/b$ ayant un dénominateur inférieur ou égal à 11.
On constate que ces vitesses de croissance sont difficiles à prévoir. Il existe cependant des suites de valeurs de $a/b$ pour lesquelles on connait la vitesse de croissance, comme par exemple les $1/3^n$ et $1/(3^n + 1)$.

\subsection{Développement $\beta$-adique avec ensemble de chiffres $\{0, 1\}$}

Considérons le semi-groupe engendré par les deux transformations affines
$$
\left\{
	\begin{array}{ccl}
		0: x &\mapsto& \beta x, \\
		1: x &\mapsto& \beta x + 1,
	\end{array}
\right.
$$
où $\beta > 1$ est un réel.
Si $\beta$ est transcendant, le semi-groupe est libre, et donc sa structure automatique est triviale.
Supposons que $\beta$ est algébrique.

\begin{define}
On appelle \defi{produit de Malher} (ou mesure de Malher) d'un nombre algébrique $\beta$ le produit
$$ m_\beta := \prod_{\begin{array}{c} v \text{ place de } \Q(\beta) \\ \abs{\gamma}_v > 1 \end{array}} \abs{\gamma}_v. $$
\end{define}
Autrement dit, le produit de Malher est le produit des modules des conjugués strictements supérieurs à $1$ et du coefficient dominant du polynôme minimal de $\beta$.

On sait que le semi-groupe est libre quand le nombre algébrique $\beta$ a un conjugué de module supérieur ou égal à 2.
Voici une réciproque :

\begin{prop}
Si le produit de Malher $m_\beta$ de $\beta$ est strictement inférieur à 2, alors le semi-groupe n'est pas libre.
\end{prop}

\emph{Preuve :}
Tous les conjugués $p$-adiques de $\beta$ sont de module 1, puisque sinon le produit de Malher serait trop grand. En particulier, le réel $\beta$ est un entier algébrique.
On peut donc plonger l'anneau des entiers de $\beta$ dans $\R^n$, de tel façon qu'il y soit un réseau.
Si l'on considère les $2^{n+1}$ polynômes en $\beta$ de degré $n$ et à coefficients dans $\{0, 1\}$, on remarque qu'ils sont inclus dans un domaine de $\R^n$ qui est de volume majoré par $P(n) m_\beta^n$, pour un certain polynôme $P$. Comme l'anneau des entiers est un réseau de $\R^n$, ceci nous donne que le nombre d'éléments du semi-groupe de longueur $n$ est majoré par $c P(n) m_\beta^n$, pour une constante $c>0$, ce qui est strictement inférieur à $2^{n+1}$ pour $n$ assez grand. D'où l'existence d'une relation non triviale dans le semi-groupe.
$\Box$ \\


La preuve ci-dessus montre que, sous les hypothèses de la proposition, l'exposant de croissance du semi-groupe est majoré par le produit de Malher.
Dans tous les exemples vérifiant les hypothèses de la proposition que j'ai pu voir, l'exposant de croissance du semi-groupe est même égal au produit de Malher.

\begin{rem}
Dans le cas particulier où le réel $\beta$ est strictement inférieur à 2, l'étude du semi-groupe est liée à celle du système dynamique
$$
	\begin{array}{ccccc}
		T :	& \R/\Z	&\rightarrow	& \R/\Z & \\
			& x		&\mapsto		& \beta x & \mod 1.
	\end{array}
$$
\end{rem}

Voici quelques exemples d'automates des relations pour le semi-groupe engendré par les applications
$$
\left\{
	\begin{array}{ccl}
		0: x &\mapsto& \beta x, \\
		1: x &\mapsto& \beta x + 1,
	\end{array}
\right.
$$
où l'étiquette $0$ signifie $(0,1)$, l'étiquette $1$ signifie $(1,0)$ et l'étiquette $*$ signifie que l'on a deux arêtes étiquetées respectivement par $(0,0)$ et $(1,1)$.

\begin{center}
	\begin{tabular}{cc}
		\includegraphics[scale=.5]{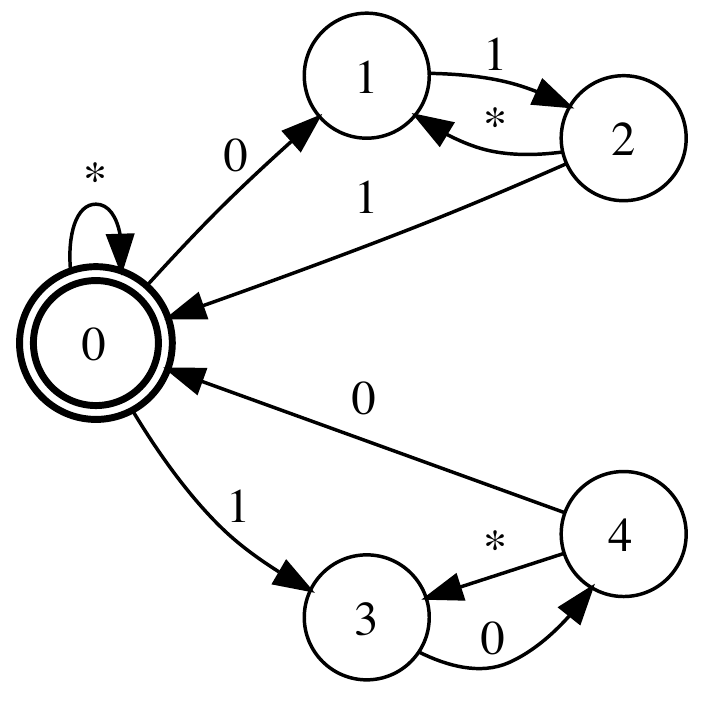} &
		\includegraphics[scale=.5]{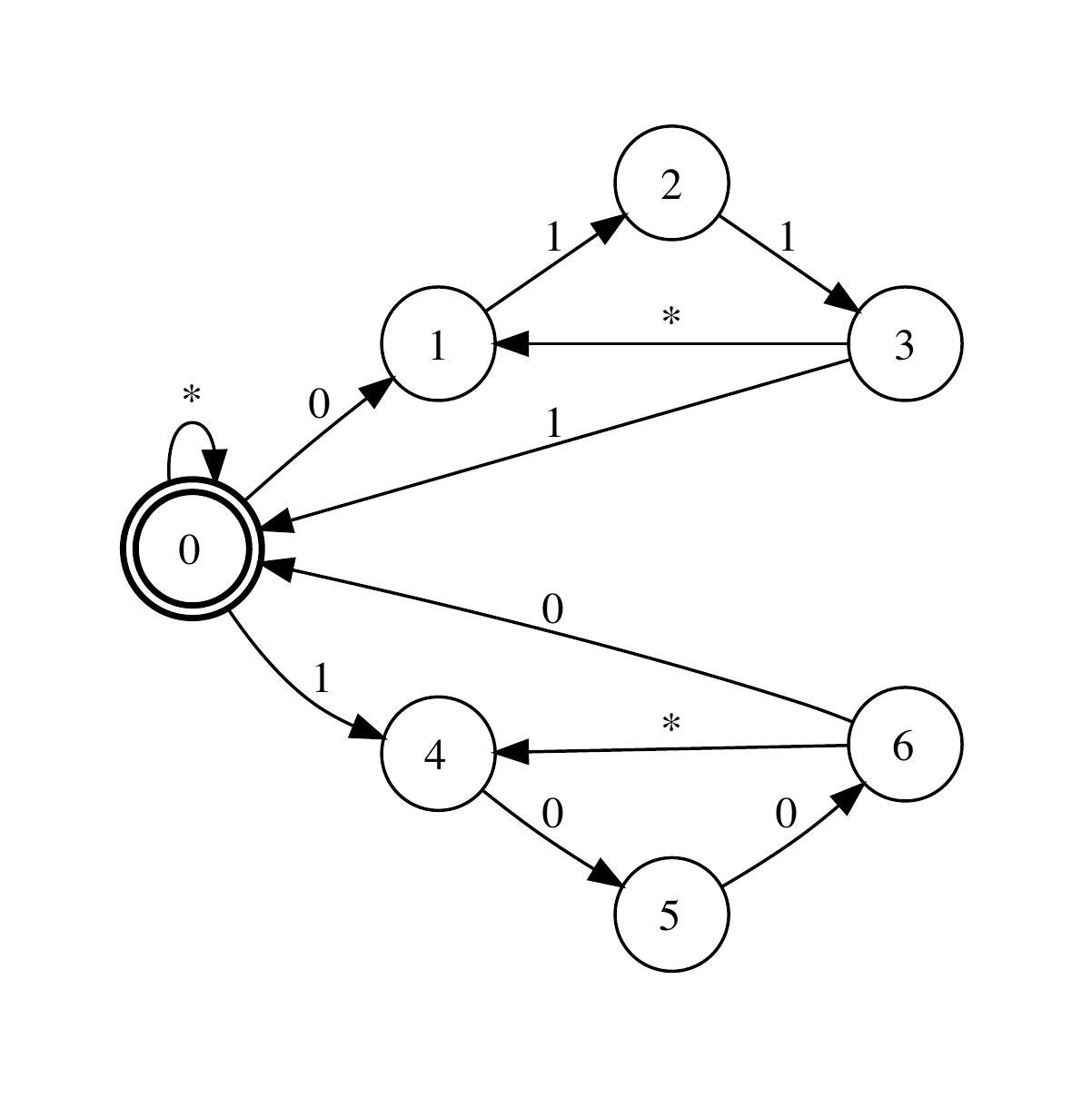} \\
		$\pi_\beta = x^2-x-1$ & $\pi_\beta = x^3-x^2-x-1$ \\ \\
	\end{tabular}
\end{center}

\begin{center}
	\includegraphics[scale=.6]{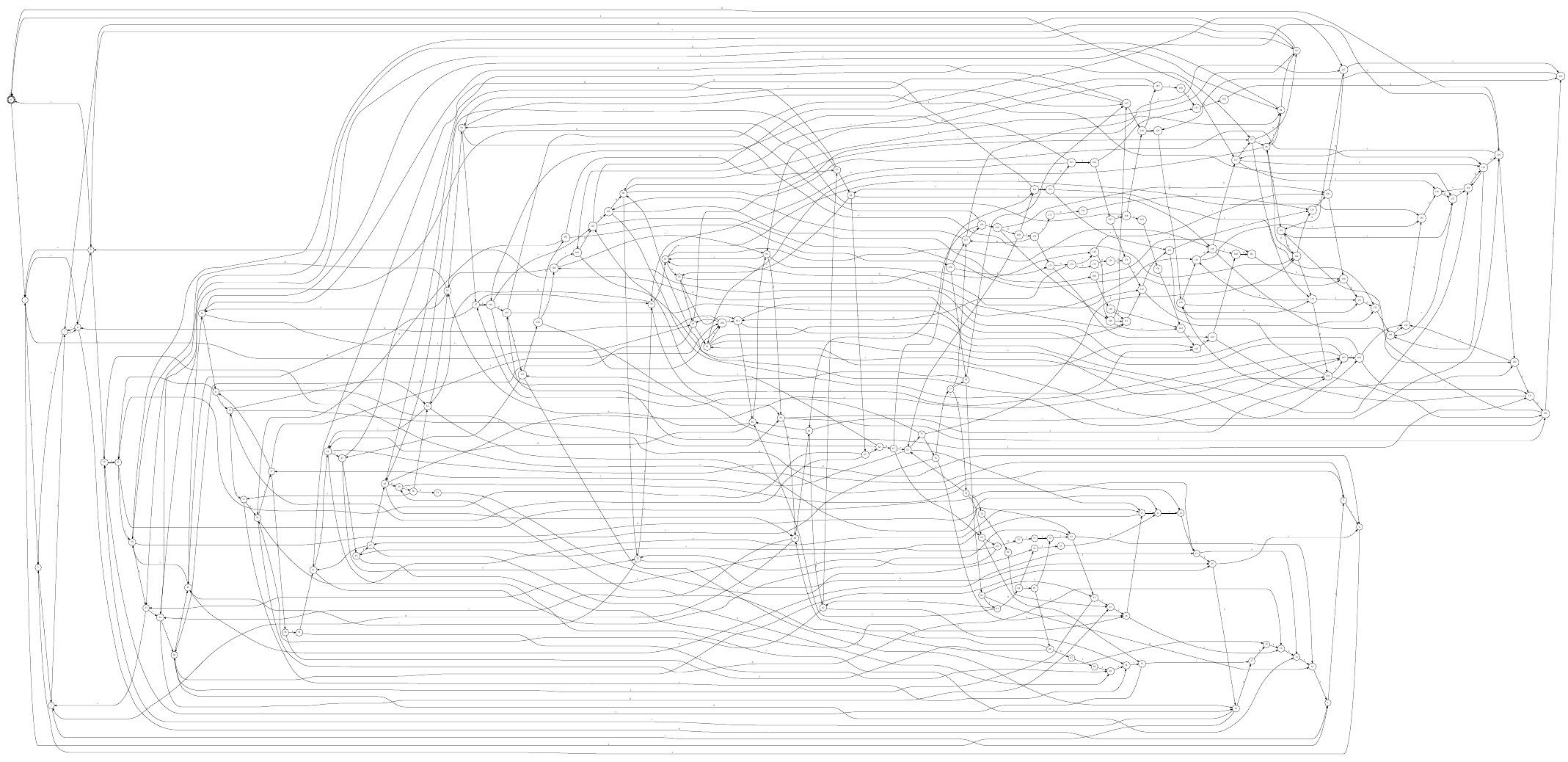} \\
	$\pi_\beta = x^3-x-1$
\end{center}

\begin{center}
	\includegraphics[scale=.5]{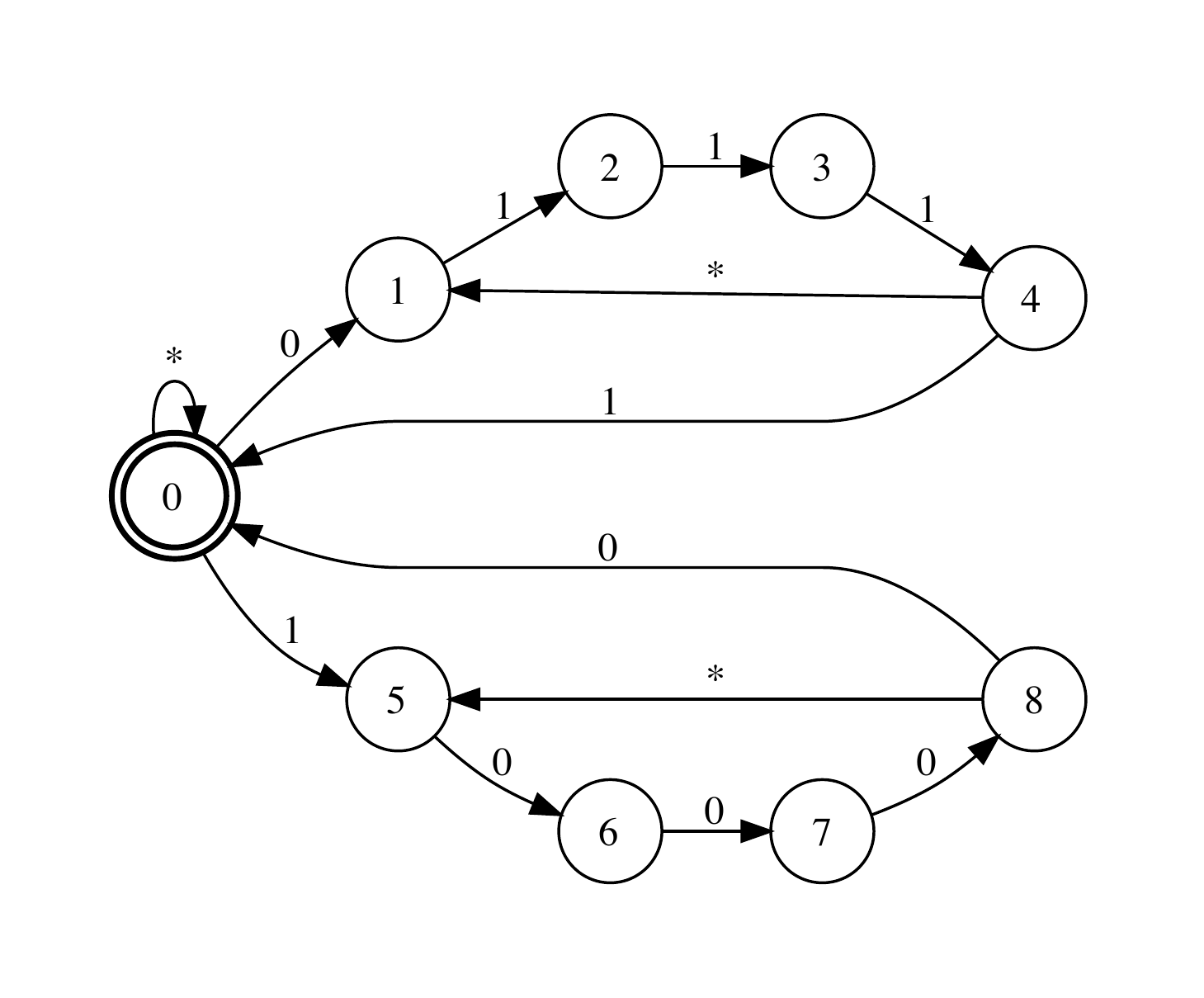} \\
	$\pi_\beta = x^4-x^3-x^2-x-1$
\end{center}


\begin{center}
	\includegraphics[scale=.6]{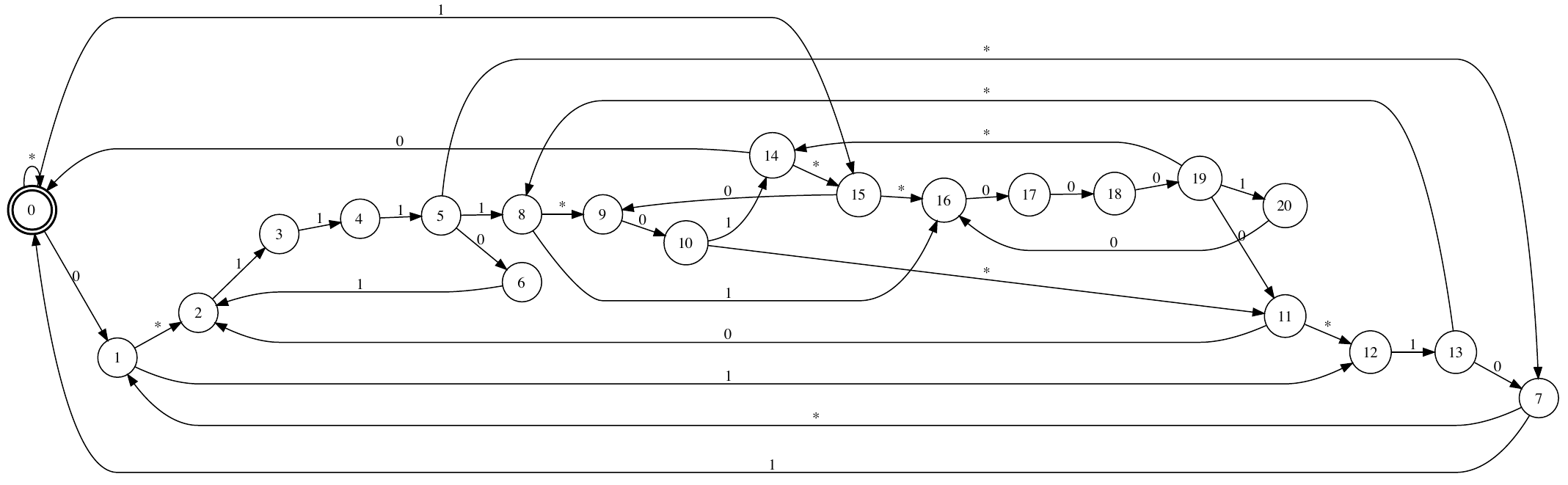} \\
	$\pi_\beta = x^4-x^3-x^2+x-1$
\end{center}

Les trois premiers exemples d'automates des relations sont pour des nombres de Pisot, tandis que ce dernier exemple est pour un nombre qui n'est pas de Pisot.
On peut voir que ces automates peuvent être très simples, mais aussi très compliqués même pour un des plus simple exemple de nombre de Pisot.

\subsection{Un cas où $\beta$ est un nombre transcendant}
Considérons le semi-groupe engendré par les 3 transformations
$$
\left\{
	\begin{array}{ccl}
		0: x &\mapsto& \beta x, \\
		p : x & \mapsto& \beta x + P(\beta), \\
		q: x &\mapsto& \beta x + Q(\beta),
	\end{array}
\right.
$$
où $\beta$ est un nombre transcendant, et $P$ et $Q$ sont deux polynômes à coefficients entiers.

\begin{rem}
On peut supposer que l'on a $\deg{P} < \deg{Q}$ et $\pgcd(P, Q) = 1$.
\end{rem}
En effet, si l'on a $\deg{Q} < \deg{P}$, alors il suffit d'échanger $P$ et $Q$, et si $\deg{P} = \deg{Q}$, alors le semi-groupe engendré par les transformations
$$
\left\{
	\begin{array}{ccl}
		0: x &\mapsto& \beta x, \\
		p : x & \mapsto& \beta x + (Q - P)(\beta), \\
		q: x &\mapsto& \beta x + Q(\beta),
	\end{array}
\right.
$$
est le même.
Si maintenant on a $\deg{P} = \deg{Q} = \deg{(Q - P)}$, alors le semi-groupe est libre, puisque l'on ne peut avoir d'égalité non triviale
$$ \sum_{i = 0}^{n} \epsilon_i \beta^i = 0,$$
avec $\epsilon_i \in \{ 0, P, Q \} - \{ 0, P, Q \} = \{ 0, P, Q, -P, -Q, P-Q, Q-P \}$.
On peut toujours supposer que les polynômes $P$ et $Q$ sont premiers entre eux quitte à tout diviser par le pgcd.

D'après la proposition \ref{pken} de Kenyon, le semi-groupe est libre dès que l'on a $P(3)+Q(3) \equiv 0 \mod 3$. Et d'après la preuve du théorème \ref{cm1}, déterminer si le semi-groupe est libre est toujours décidable, puisque la structure fortement automatique du semi-groupe est calculable. La remarque \ref{rem_lib} donne un critère pour déterminer si le semi-groupe est libre. Mais existe-t'il un critère simple pour déterminer si le semi-groupe est libre à partir des polynômes $P$ et $Q$ ?

Voici quelques exemples d'automates des relations pour ce semi-groupe

\begin{center}
	\begin{tabular}{cc}
		\includegraphics[scale=.6]{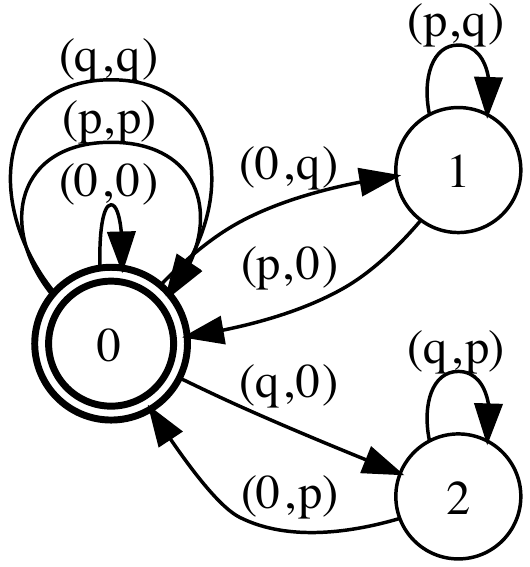} &
		\includegraphics[scale=.6]{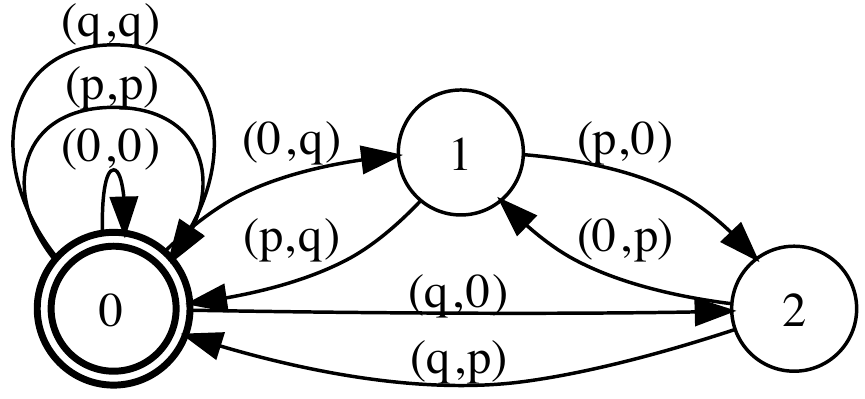} \\
		P=1 et Q=X & P=1 et Q=X+1 \\ \\
%
%
	\end{tabular}
\end{center}
\begin{center}
	\includegraphics[scale=.5]{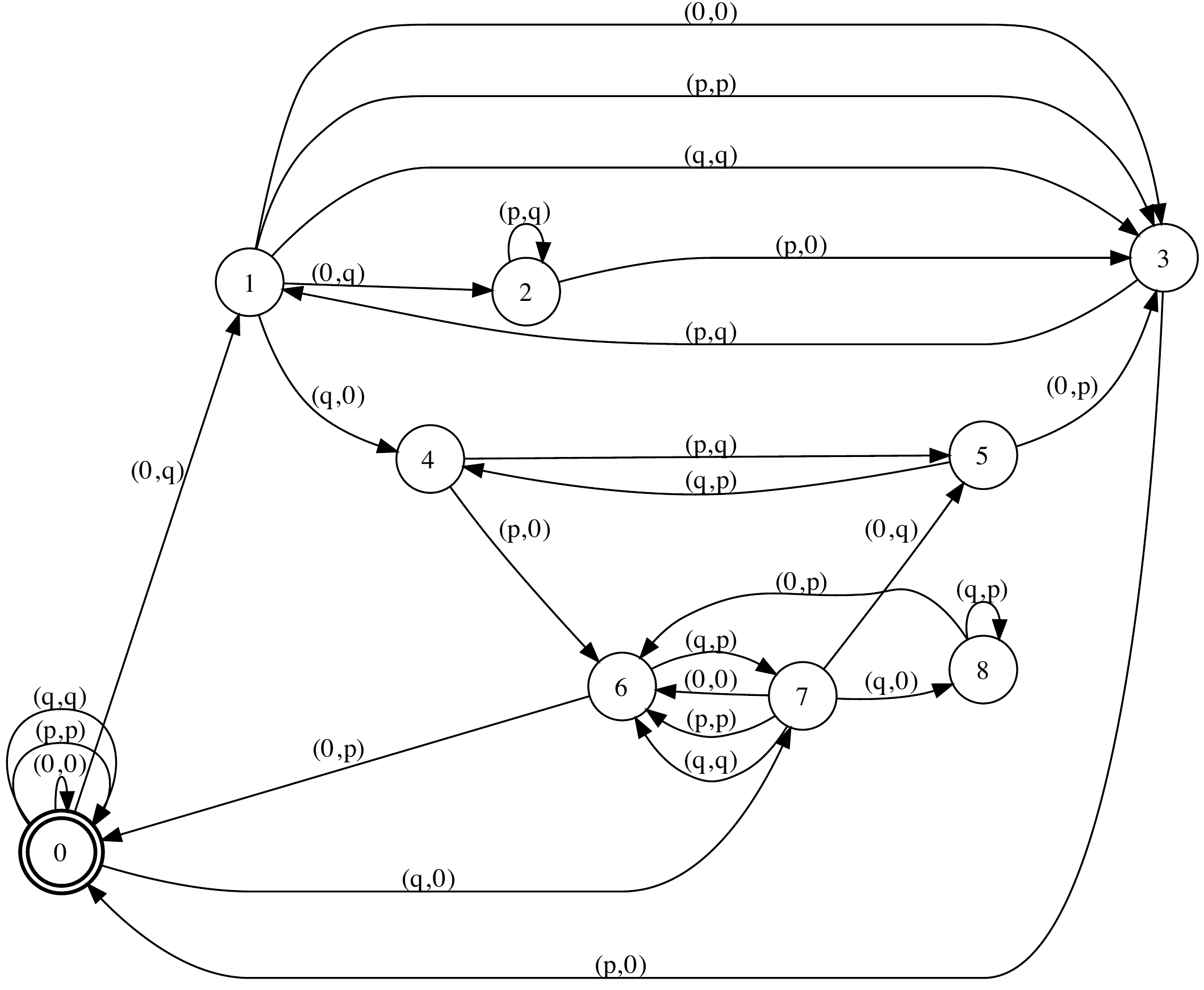}
\end{center}
\begin{center}
	$P=1$ et $Q = X^2$ \\
\end{center}

\begin{center}
	\includegraphics[scale=.5]{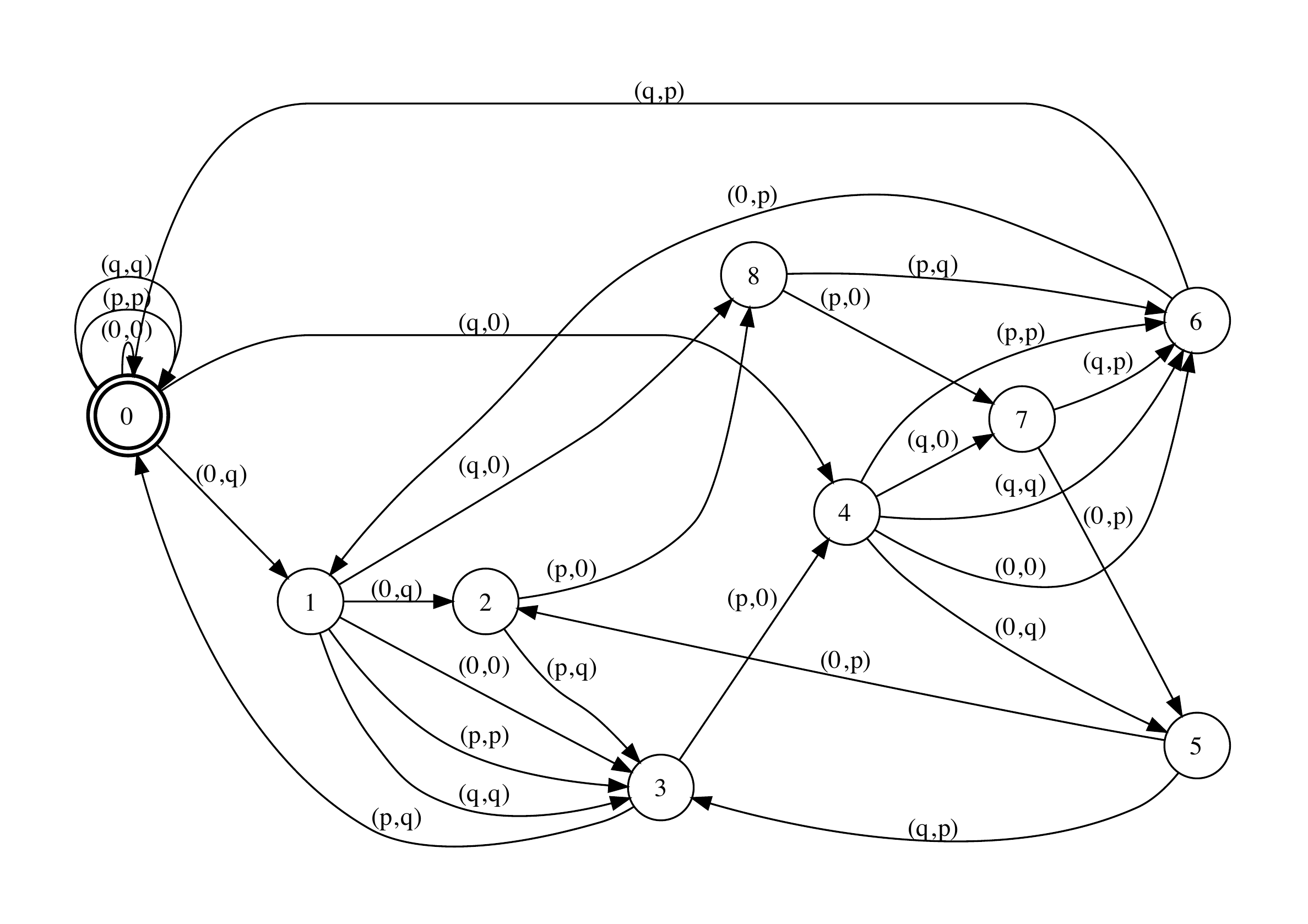}
\end{center}
\begin{center}
	$P=1$ et $Q = X^2+1$ \\
\end{center}

\begin{center}
	\includegraphics[scale=.5]{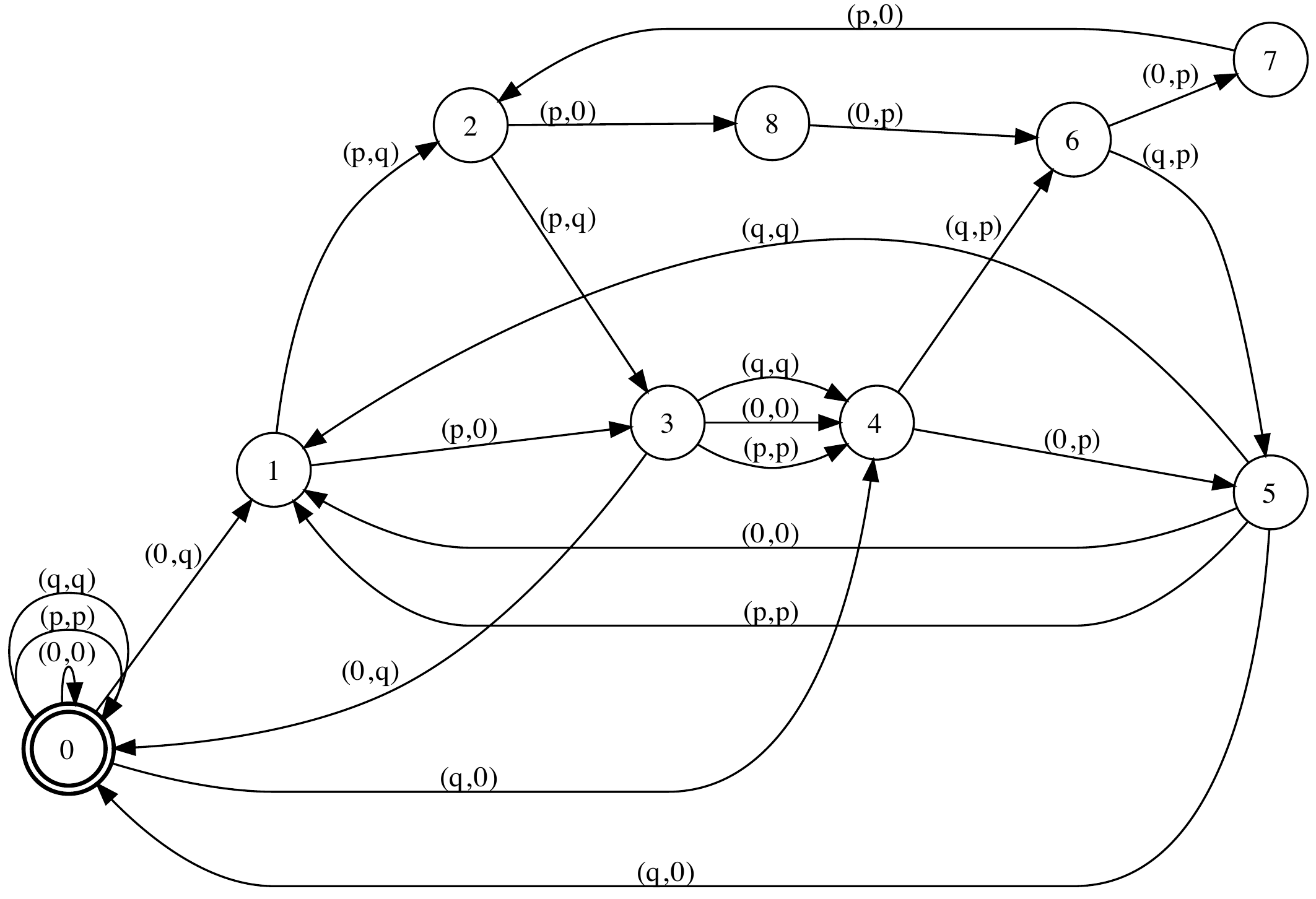}
\end{center}
\begin{center}
	$P=X$ et $Q = X^2+1$ \\
\end{center}

\begin{center}
	\includegraphics[scale=.5]{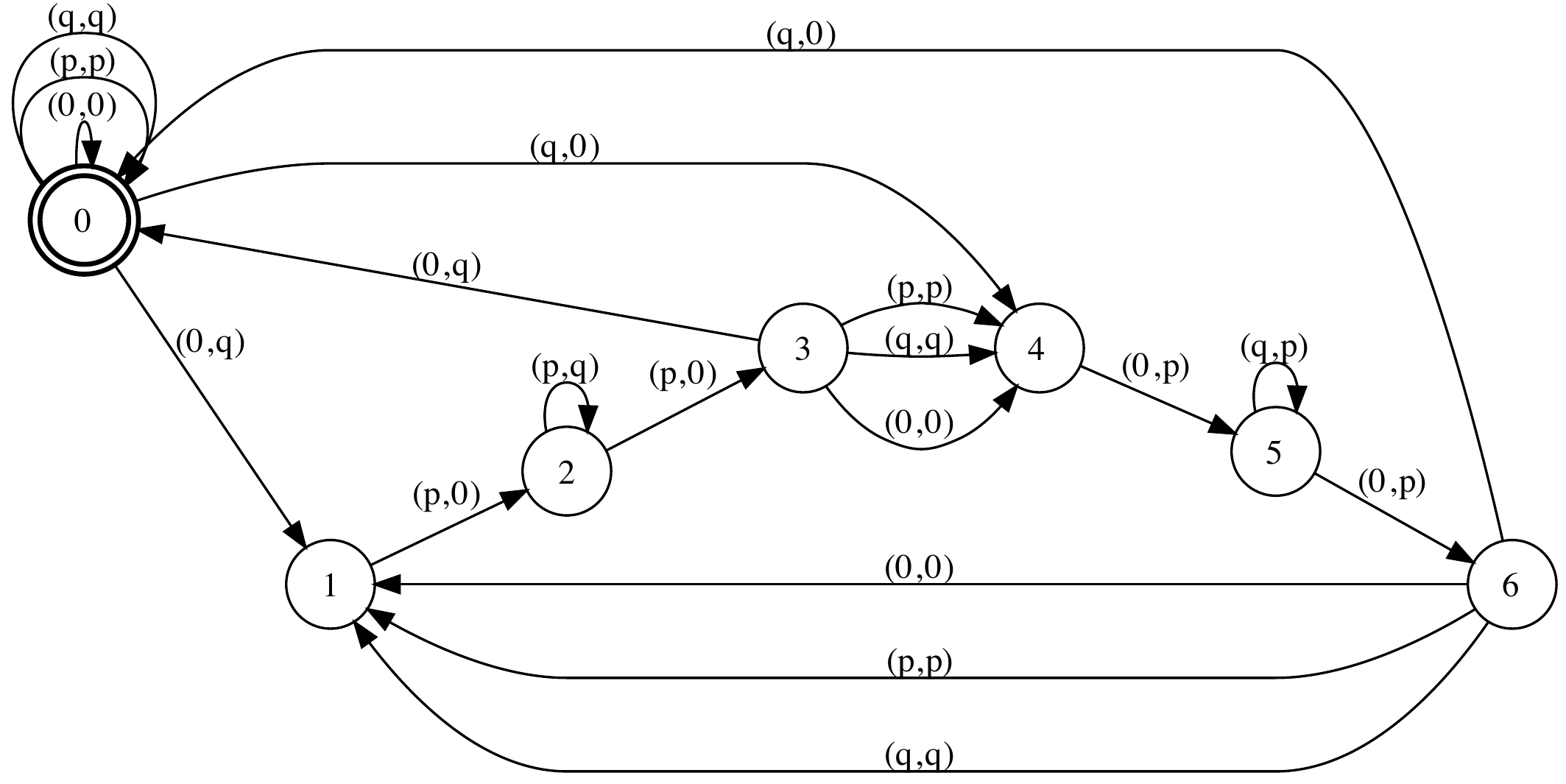}
\end{center}
\begin{center}
	$P=X$ et $Q = X^2-X+1$ \\
\end{center}

Et voici les automates des mots réduits pour les mêmes exemples :

\begin{center}
	\begin{tabular}{cc}
		\includegraphics[scale=.6]{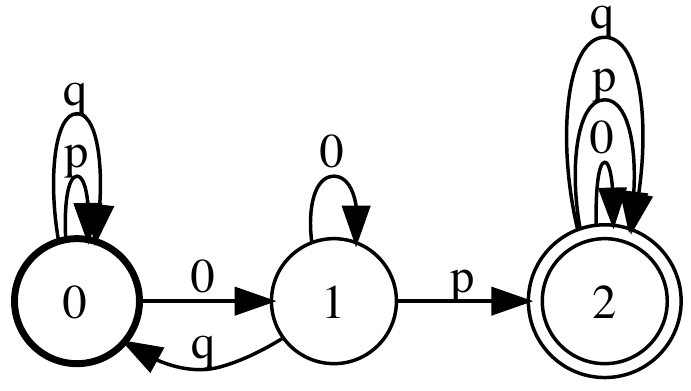} &
		\includegraphics[scale=.6]{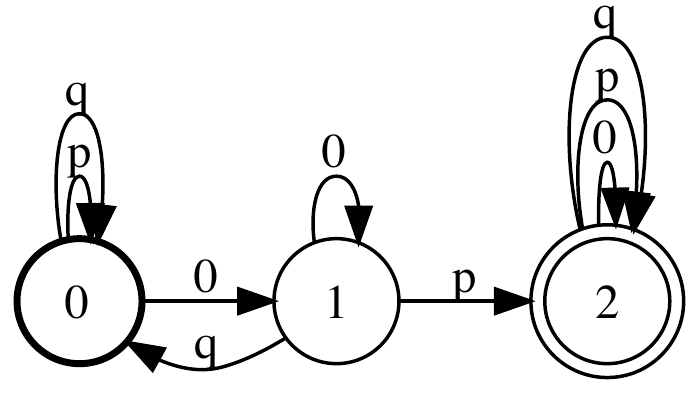} \\
		$P=1$ et $Q=X$ & $P=1$ et $Q=X+1$ \\ \\
	\end{tabular}
\end{center}

\begin{center}
	\begin{tabular}{cc}
		\includegraphics[scale=.4]{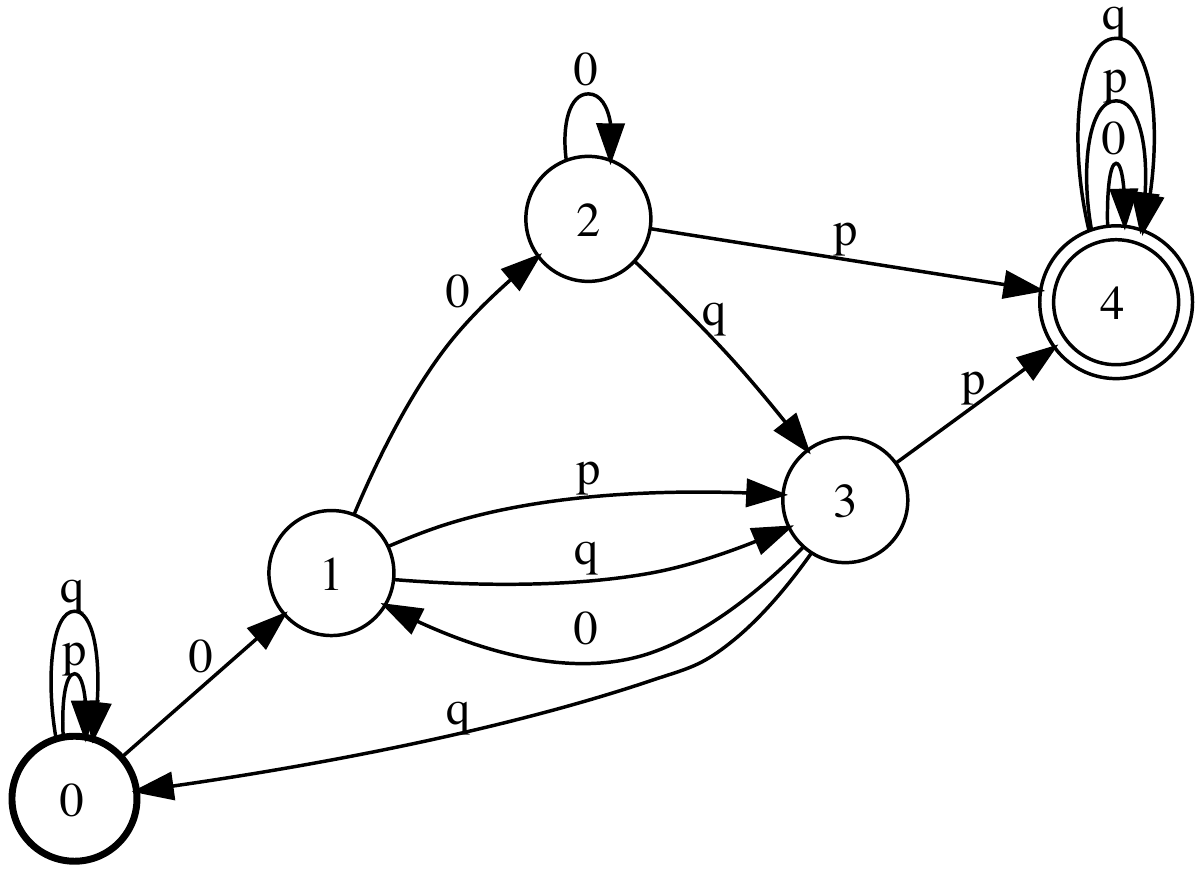} &
		\includegraphics[scale=.4]{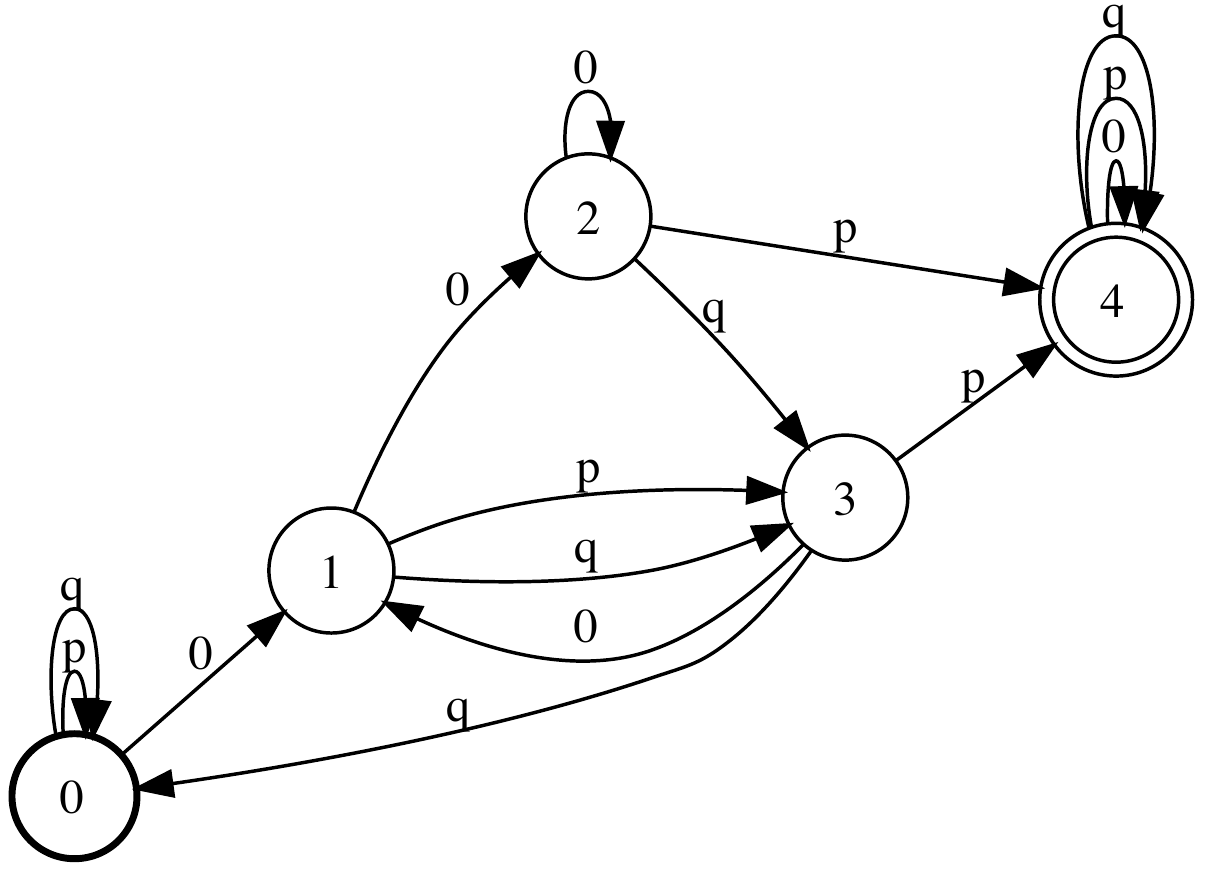} \\
		$P=1$ et $Q=X^2$ & $P=1$ et $Q=X^2+1$ \\ \\
	\end{tabular}
\end{center}

%

\begin{center}
	\includegraphics[scale=.5]{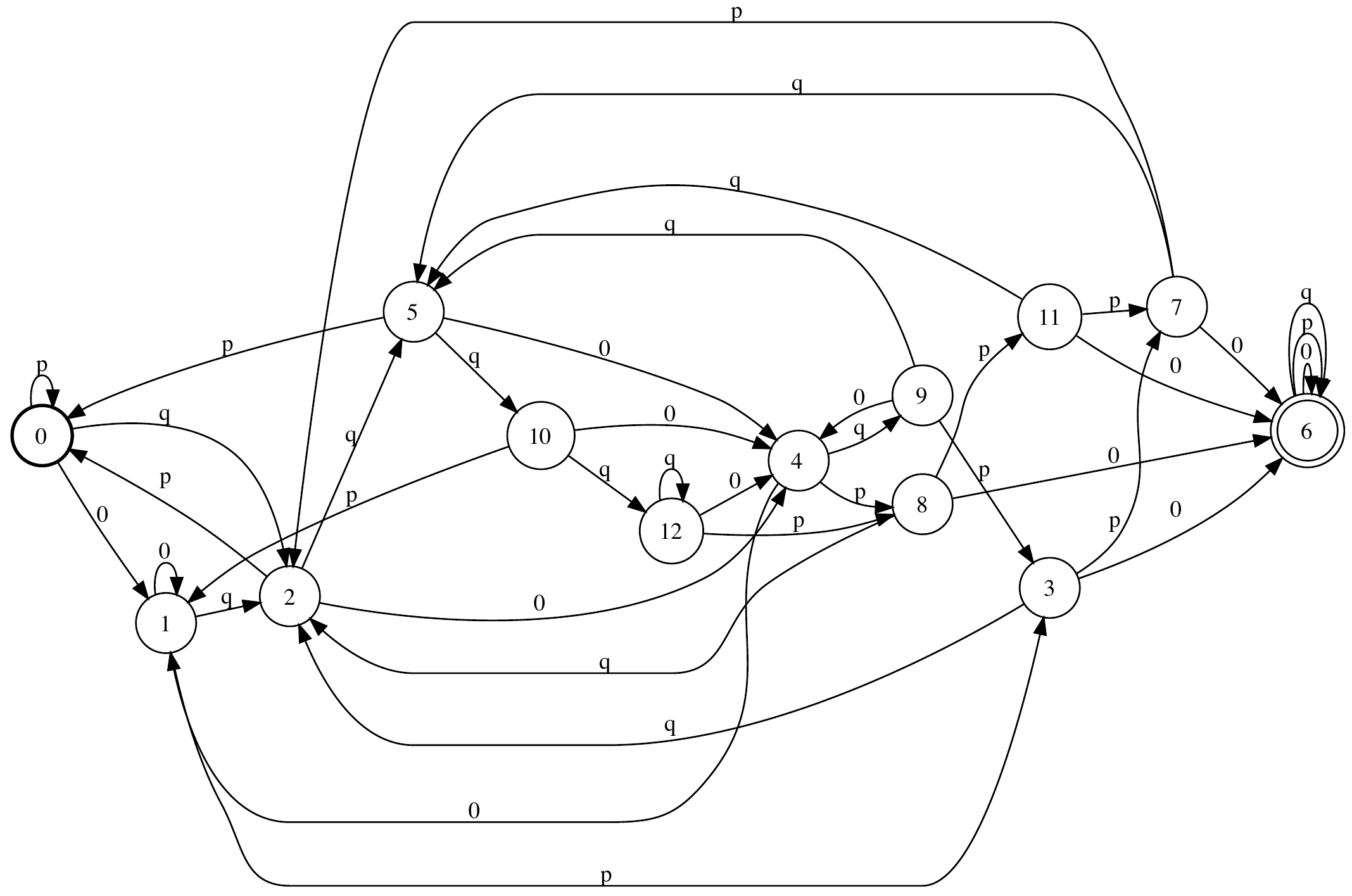}
\end{center}
\begin{center}
	$P=X$ et $Q = X^2+1$ \\
\end{center}

\begin{center}
	\includegraphics[scale=.5]{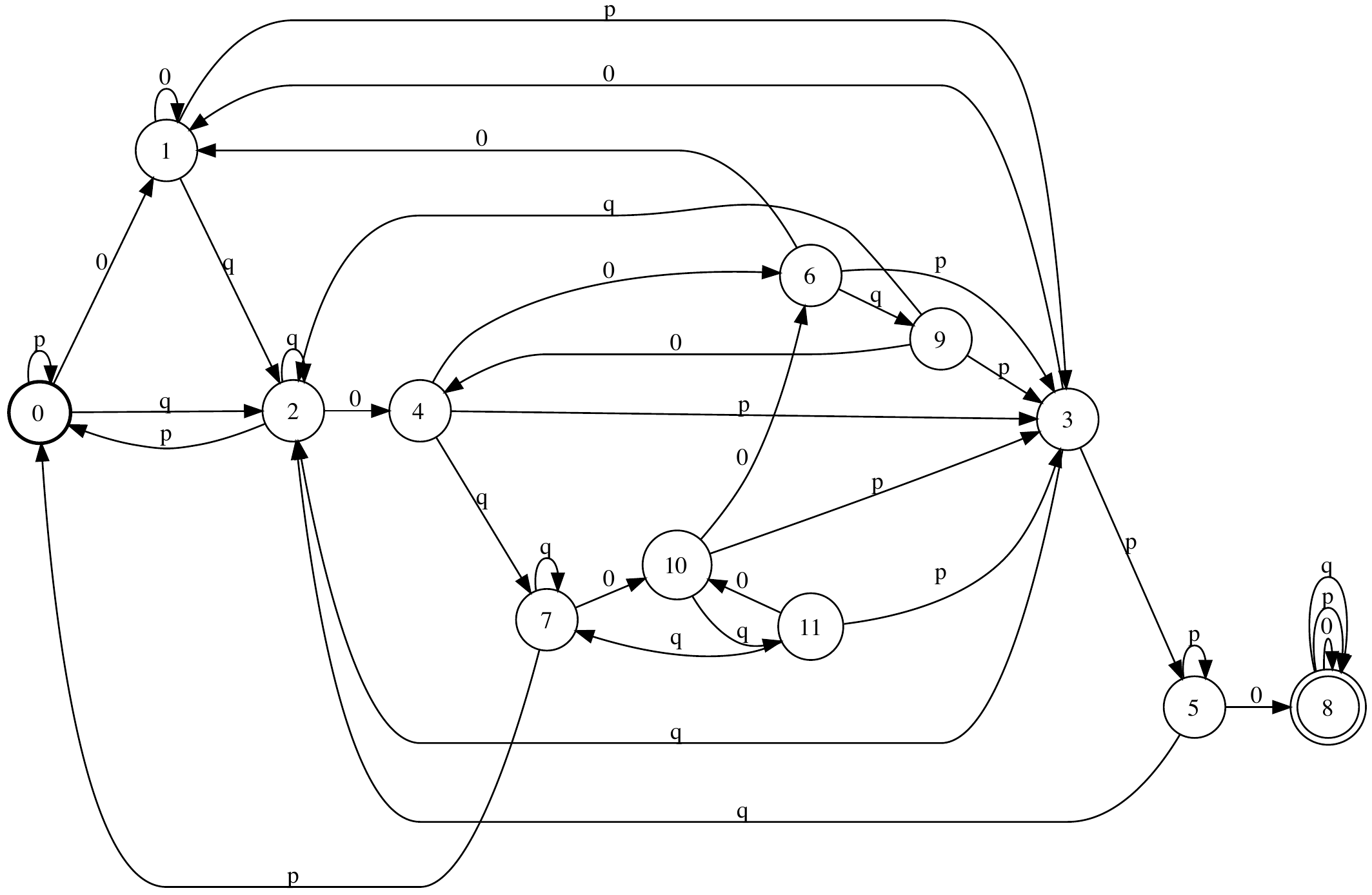}
\end{center}
\begin{center}
	$P=X$ et $Q = X^2-X+1$ \\
\end{center}

Voici les vitesses exponentielles de croissance $\lambda$ pour quelques exemples, où $\pi_\lambda$ est le polynôme minimal de $\lambda$ :

\begin{center}
	\begin{tabular}{| c | c | c |}
		\hline
		$P/Q$ & $\lambda$ & $\pi_\lambda$ \\
		\hline
		$1/X$		& 2.6180	& $x^2-3x+1$ \\
		$1/(X+1)$		& 2.6180	& $x^2-3x+1$ \\
		$1/(X^2-X)$	& 2.8794	& $x^3 - 3 x^2 + 1$ \\
		$1/(X^2-X+1)$	& 2.7971	& $x^4-2x^3-2x^2-x+1$ \\
		$1/X^2$		& 2.6180	& $x^2-3x+1$ \\
		$1/(X^2+1)$	& 2.6180	& $x^2-3x+1$ \\
		$1/(X^2+X)$	& 2.8794	& $x^3-3x^2+1$ \\
		$1/(X^2+X+1)$	& 2.7693	& $x^3-3x^2+x-1$ \\
		$(X-1)/X^2$	& 2.7971	& $x^4-2x^3-2x^2-x+1$ \\
		$(X-1)/(X^2+X-1)$ & 2.8794	& $x^3-3x^2+1$ \\
		$1/(X^3-X^2-X)$	& 2.9615	& $x^4-3x^3+1$ \\
		$1/(X^3-X^2)$ & 2.8584	& $x^7-3x^6+3x^3+x^2-1$ \\
		$1/(X^3-X^2+1)$ & 2.8396	& $x^{10}-3x^9+3x^6+x^5+4x^4-3x^3-3x^2+1$ \\
		$1/(X^3-X^2+X)$	& 2.8444	& $x^{13}-3x^{12}-2x^{11}+7x^{10}-2x^9+7x^8-$\\
						&		& $16x^6+ 6x^5-6x^3+8x^2+x-2$ \\
		
		\hline
	\end{tabular}
\end{center}

On remarque qu'à nouveau ces vitesses exponentielles de croissance sont difficiles à prévoir, mais qu'il y a tout de même des valeurs particulières pour lesquelles on les connait (par exemple les $1/X^n$).

\begin{rem}
La vitesse exponentielle de croissance du semi-groupe pour $\beta$ transcendant majore celle du semi-groupe pour $\beta$ algébrique.
Pour l'exemple de l'introduction (qui correspond aux polynômes $P=1$ et $Q=X$), le caractère algébrique de $1/3$ n'a aucun rôle : le semi-groupe est le même (d'un point de vue combinatoire) en prenant $\beta$ transcendant plutôt que $\beta = 1/3$. Cependant, les semi-groupes diffèrent quand on prend par exemple $P=X$ et $Q=X^2-X+1$ suivant que $\beta=1/3$ ou que $\beta$ est transcendant.
\end{rem}

\begin{rem} \label{rem_lib}
Le semi-groupe n'est pas libre si et seulement si l'on a $P/Q = A/B$ pour deux polynômes $A, B \in \{-1, 0, 1\}[X]$ et avec $A-B \in \{-1, 0, 1\}[X]$.
C'est pourquoi tous les exemples considérés ci-dessus sont de cette forme.
\end{rem}

\end{document}